%% file: waveletextension.tex
\documentclass[a4paper,10pt]{article}

\usepackage{mymath}
{
	\theoremstyle{plain}
	\newtheorem{assumption}{Assumption}
}

\usepackage[normalem]{ulem}
\ifdefined\C
\renewcommand{\C}{\mathbb{C}}
\else
\newcommand{\C}{\mathbb{C}}
\fi
\title{Efficient function approximation on general bounded domains using wavelets on a cartesian grid}
\author{Vincent Copp\'e\footnote{Email: \texttt{vincent.coppe@cs.kuleuven.be}. Website: \texttt{https://people.cs.kuleuven.be/\textasciitilde vincent.coppe}.} \and Daan Huybrechs\footnote{Email: \texttt{daan.huybrechs@cs.kuleuven.be}. Website: \texttt{https://people.cs.kuleuven.be/\textasciitilde daan.huybrechs}.}}

\date{KU Leuven \\ Department of Computer Science \\ Celestijnenlaan 200A \\ 3001 Leuven, Belgium \\ \vskip10pt \today}

\newcommand*{\SPAN}{\mathrm{span}}

\usepackage{tikz}
\usepackage{pgfplots}
\usepgfplotslibrary{groupplots}
\def\osi{q}

\def\K{{\mathcal K}}
\def\L{{\mathcal L}}
\def\M{{\mathcal M}}
\begin{document}
\maketitle
\begin{abstract}
Fourier extension is an approximation method that alleviates the periodicity requirements of Fourier series and avoids the Gibbs phenomenon when approximating functions. We describe a similar extension approach using regular wavelet bases on a hypercube to approximate functions on subsets of that cube. These subsets may have a general shape. This construction is inherently associated with redundancy which leads to severe ill-conditioning, but recent theory shows that nevertheless high accuracy and numerical stability can be achieved using regularization and oversampling. Regularized least squares solvers, such as the truncated singular value decomposition,  that are suited to solve the resulting ill-conditioned and skinny linear system generally have cubic computational cost. We compare several algorithms that improve on this complexity. The improvements benefit from the sparsity in and the structure of the discrete wavelet transform. We present a method that requires $\mathcal O(N)$ operations in 1-D and $\mathcal O(N^{3(d-1)/d})$ in $d$-D, $d>1$. We experimentally show that direct sparse QR solvers appear to be more time-efficient, but yield larger expansion coefficients.
\end{abstract}

\noindent \textbf{Keywords} \, Fourier extension, wavelets, efficient algorithms, frames, ill-conditioning, function approximation, oversampling  \\

\noindent \textbf{Mathematics Subject Classification (2010)} \,65D15, 65T60, 65Y20
\section{Introduction}
Wavelets have many applications in signal processing. Their most common uses are in compression, edge detection, denoising and other signal enchancements. The wide applicability of wavelets is mostly due to the localization properties of wavelets in time and frequency, such that many signals can be sparsely represented, as well as to the existence of the (bi)orthogonal Discrete Wavelet Transform that can be applied very efficiently. Wavelets are also applied in solution methods for partial differential equations or integral equations using wavelet-based discretizations \cite{beylkin1991wavelet,dahmen1997opeq,cohen2003adaptivewaveletschemes,stevenson04compressibility,dahmen06compression}. Here, too, the localization properties of wavelets are of interest. In this setting, wavelets with suitable level-dependent scalings can generate stable bases for a range of function spaces, and the existence of dual bases with varying smoothness is frequently useful.
%\GREEN{In welke settings is het hebben van duale basissen met verschillende smoothness nuttig? Is dit omdat met het gebruik van quadratuur coefficienten gevonden worden zonder het oplossen van een stelsel? Waarom is het hebben van veschillende smoothness hierbij interessant? Niet per se om hier uit te leggen, maar uit interesse. In welke paper wordt gebruik gemaakt van deze zaken?}
% Zie Dahmen, 1997, hoofdstuk 5, bv. (5.3.11) en (5.3.12)
% Voor het belang van de orde van de duale MRA: zie ook hoofdstuk 6 Preconditioning

However, it is in general difficult to create a wavelet basis on a complex geometry, i.e., to create a basis for a domain that is not a hypercube. Several methods have been proposed to enable the use of wavelet on general domains. One method is based on solving a Dirichlet problem with a fictitious domain method \cite{Wells1995}. Other methods employ adaptive finite element techniques \cite{Cohen2000,Barinka2001,Stevenson2003}. However, in \cite{Stevenson2003} it is stated that the required smoothness conditions for efficient adaptive wavelet methods are difficult to impose on domains that do not have product structure. A wavelet frame method is constructed in \cite{Stevenson2003} that can cope with domains that are overlapping unions of subdomains, each of them being the image under a smooth parametrization of a hypercube. Frames generalize a basis in the sense that they allow for redundancy \cite{christensen2016}.

As in \cite{Stevenson2003}, in this paper we resort to using a wavelet frame rather than a wavelet basis. However, we use a different type of frame and we restrict ourselves to the problem of function approximation rather than the solution of operator equations. The approximation problem we consider is the same as that considered in \cite{matthysen2018function,coppe2019splines} based on Fourier series and splines respectively. We aim for a fast algorithm for the approximation of a function $f$ on a compact domain $\Omega$ that can have an arbitrary shape. Without loss of generality, this bounded domain can be scaled such that $\Omega\subset\Xi$ with $\Xi=[0,1]^d$. With Fourier extensions, (tensor products of) Fourier series are used on $\Xi$, while in \cite{coppe2019splines} an analogous spline extension is introduced based on a periodic spline basis on $\Xi$. Here, we extend the idea further to wavelets.
%The main advantage is that we can use a standard wavelet basis on the cube, with an associated efficient (discrete) wavelet transform. The main disadvantage is that the restriction of that basis to the subset $\Omega$ constitutes a frame and this complicates the computation of function approximations. That is the topic of this paper.

While it is difficult to create a wavelet basis on $\Omega$, it is easy to create one on $\Xi$. Consider for example tensor products of Daubechies or CDF wavelets, periodized to the interval $[0,1]$ \cite{Daubechies1988,Cohen1992}. If we restrict the basis to $\Omega$, we naturally arrive at a frame that we will call a \emph{wavelet extension frame}. In this paper we focus on algorithms, rather than on the properties of a frame. Nevertheless, we recall its definition and the motivation for its use in function approximation.

A family of functions $\Phi = \{\phi_k\}_{k=1}^\infty$ is a frame for a Hilbert space $\mathcal H$ if \cite[Def. 5.1.1]{christensen2016}
\begin{equation*}\label{eq:frame}
	A\|f\|^2\leq\sum_{k=1}^\infty |\langle f,\phi_k\rangle|^2  \leq B\|f\|^2,\qquad \forall f\in\mathcal H
\end{equation*}
for constants $A,B>0$. It is more general than a basis, as demonstrated by the construction above. In particular frames may be redundant. In the setting of this paper, redundancy arises since our approximant can take any value in $\Xi\setminus\Omega$. This leads to apparent ill-conditioning of the approximation problem. However, recent theory indicates that the ill-conditioning of the linear systems to be solved does not prevent stable and highly accurate function approximation if one uses regularization techniques in combination with oversampling \cite{Adcock2019,adcock2017frames}. For that reason we consider least squares approximations and develop an efficient regularizing solver.

Wavelets are by their nature adaptive. It is possible to extend or refine a wavelet basis by adding basis functions on a finer scale. This is not possible in a spline basis. A translation-invariant spline basis $\Phi_{N}=\{\phi(\cdot-hk)\}_{k=-\infty}^\infty$ with $h > 0$ can be refined by dilating the basis functions, but all basis functions change as a result. It is the possibility of adaptivity of wavelets that motivates their study in this paper. However, we will not (yet) fully take advantage of the possibilities. We do use one form of adaptivity at the end of the paper to arrive at a wavelet extension approximation with a smooth extension by choosing level-dependent weights, which is not possible in the context of spline extensions.

Though the methods of the paper are general, we consider in our examples the Daubechies and Cohen--Daubechies--Feauveau (CDF) family of wavelets, since they are widely used and have compact support. The duals of these bases are a key ingredient in the construction of efficient solvers. For Daubechies and CDF wavelets, dual bases in $L^2(\R)$ are well studied. They can be used for function approximation using a Galerkin-type approach, i.e., based on (bi)orthogonal projections using inner products. However, inner products with wavelets on general domains are not easily computed, especially not in the multivariate case, since they require the numerical evaluation of integrals on domains of general (and possibly irregular) shape. Instead, we focus in our experiments on a collocation approach based on discrete function samples. Collocation and oversampling necessitate the construction of bases that are dual with respect to a discrete oversampled equispaced grid. We provide such a construction on the bounding box using cartesian grids, taking advantage of their regular structure, and demonstrate how this construction can be used for the efficient solution on the subdomain of general shape.

The structure of the paper is as follows. In \S\ref{wavs:wavelets}, wavelets are introduced along with the discrete wavelet transform. The structure of the latter can be used to create efficient matrix-vector products. We recall these basics in order to modify them later on. In \S\ref{wavs:periodicwavelet}, bases biorthogonal to  periodic wavelet bases on the interval are discussed. We describe the construction of discrete dual bases. In \S\ref{wavs:approximationproblem} we discretize the function approximation and arrive at a matrix system. Next, we compare several algorithms to solve this system in \S\ref{wavs:az}.  Finally, we use the adaptive nature of wavelets to construct a smooth extension in \S\ref{wavs:smooth} and end with some concluding remarks in \S\ref{wavs:conclusion}.

\section{Wavelets}\label{wavs:wavelets}

Wavelets may be created by dilating and translating a given function.  For particular choices of $\psi(t)\in L^2(\R)$, the family
\begin{equation}\label{eq:motherfunction}
	\psi_{jk}(t) = {2^{j/2}}\psi(2^{j}t-k), \qquad j,k\in\Z
\end{equation}
forms a basis for $L^2(\R)$. This family is a wavelet basis and $\psi(t)$ is called the \emph{mother function}. There exists a great variety of other types of wavelets, some giving rise to a frame rather than a basis \cite{christensen2016,daubechies1992ten}. However, we limit ourselves here to orthogonal and biorthogonal wavelet bases that are translation invariant as above, with compact support, and that can be constructed using a multiresolution analysis \cite{Mallat1989,mallat1989multiresolution,Daubechies1988,Cohen1992}. We can take advantage of their regular structure to implement efficient operations.

\subsection{Multiresolution analysis}

A multiresolution analysis in the context of wavelets was introduced in \cite{mallat1989multiresolution} and can be defined as follows.
\begin{definition}\cite[Definition 7.1]{Mallat2009}\label{def:multiresolution}
	A multiresolution analysis of $L^2(\R)$ is a nested sequence
	\begin{equation}\label{eq:nestedsubspaces}
	\cdots\subset V_{-2}\subset V_{-1}\subset V_{0}  \subset V_{1}\subset V_{2}\subset \cdots
	\end{equation}
	 of closed subspaces of $L^2(\R)$ such that
	\begin{enumerate}
		\item $\lim_{j\rightarrow\infty}V_{j}=\overline{\bigcup_{j\in\Z}V_j}=L^2(\R)$ and $\lim_{j\rightarrow-\infty}V_{j}=\bigcap_{j\in\Z}V_j=\emptyset$.
		\item $f(t)\in V_{j}\Leftrightarrow f(2t)\in V_{j+1}$, $\forall j\in\Z$.
		\item $f(t)\in V_0\Leftrightarrow f(t-k)\in V_0$, $\forall k\in\Z$.
		\item there exists a $\phi(t)\in V_0$ such that $\{\phi(\cdot-k)\}_{k\in\Z}$ forms a Riesz basis for $V_0$.
	\end{enumerate}
\end{definition}

The first condition states that the sequence of subspaces is a non-redundant approximation of $L^2(\R)$. The second and third condition introduce scale and translation invariance. The last one demands the existence of a translation invariant basis for $V_0$. Similar to the wavelet mother function~\eqref{eq:motherfunction} we call the function $\phi(t)$ introduced in Definition~\ref{def:multiresolution} the father function. Analogously to the mother function it generates a family of functions:
\begin{equation}\label{eq:scalingbases}
	\phi_{jk}(t)={2^{j/2}}\phi(2^{j}t-k), \qquad \forall k\in\Z, \quad j\in\Z.
\end{equation}
For every $j$, $\{\phi_{jk}\}_{k\in\Z}$ forms a Riesz basis of $V_j$. We call $\{\phi_{jk}\}_{k\in\Z}$ a scaling basis of $V_j$. More specifically, if  $\{\phi(\cdot-k)\}_{k\in\Z}$ forms an orthonormal basis for $V_0$, $\{\phi_{jk}\}_{k\in\Z}$ forms an orthonormal basis of $V_j$ for every $j\in\Z$. In that case, we can also define the sequence of orthogonal projections
\begin{equation}\label{eq:projection}
 {\mathcal P}_j f = \sum_{k \in \Z} \langle f, \phi_{jk} \rangle \, \phi_{jk}.
\end{equation}

The scale invariance of the multiresolution analysis in Definition~\ref{def:multiresolution} implies the existence of a two-scale relation
\begin{equation}\label{eq:twoscalerelation}
	\phi(t) = \sqrt 2\sum_{k\in\Z}h_k\phi(2t-k)
\end{equation}
in which $h_k$ is a sequence. If we require that $\int_{\R}\phi(t)\d t\neq 0$, we have that
\begin{equation*}
	\sum_{k\in\Z}h_k=1.
\end{equation*}
If we further require that $\{\phi(\cdot-k)\}_{k\in\Z}$ forms an orthonormal basis for $V_0$, then $h$ satisfies so-called double shift orthogonality conditions:
\begin{equation}\label{eq:doubleshiftorthonormality}
	\sum_{k\in\Z}h_k\overline h_{k+2n}=\delta_{0n},\qquad \forall n\in\Z.
\end{equation}

A wavelet basis that follows from an orthonormal multiresolution analysis is
\begin{equation}\label{eq:orthonormalwaveletbasis}
	\psi_{jk}(t)=\sqrt 2\sum_{l\in\Z} g_k \phi_{jl}(2t-k),\qquad \forall k,j\in\Z
\end{equation}
with $g_k=(-1)^k\overline h_{-k+1}$ and $h$ as in~\eqref{eq:twoscalerelation}. It forms an orthonormal basis for~$L^2(\R)$.

\subsection{Biorthogonal multiresolution analysis}

Orthogonality is a rather restrictive requirement. A compactly supported and symmetric sequence $h$ that satisfies double shift orthogonality~\eqref{eq:doubleshiftorthonormality} can only have two non-zero coefficients \cite[Proposition 4.1.]{Daubechies1988}. That restriction is lifted using biorthogonal wavelets.

To construct biorthogonal wavelets we create, as in~\cite{Cohen1992}, a biorthogonal multiresolution analysis. To that end, next to the first (primal) multiresolution~\eqref{eq:nestedsubspaces}, we define a second (dual) one
\begin{equation*}\label{eq:dualnestedsubspaces}
	\cdots\subset \tilde V_{-2}\subset \tilde V_{-1}\subset \tilde V_{0}  \subset \tilde V_{1}\subset \tilde V_{2}\subset \cdots
\end{equation*}
for which a dual scaling function~$\tilde\phi(t)\in\tilde V_0$ exists such that $\{\tilde\phi(\cdot-k)\}_{k\in\Z}$ forms a Riesz basis of~$\tilde V_0$. The dual scaling function satisfies the two-scale relation
\begin{equation*}\label{eq:dualtwoscalerelation}
	\tilde\phi(t) = \sqrt 2\sum_{k\in\Z}\tilde h_k\tilde\phi(2t-k).
\end{equation*}
We call this basis a dual scaling basis, while the scaling bases in \eqref{eq:scalingbases} are primal scaling bases.  If
\begin{equation*}\label{eq:biorthogonalscalingbases}
	\left\langle\phi(\cdot-k),\tilde\phi(\cdot-l)\right\rangle_{L^2(\R)}=\delta_{kl},
\end{equation*}
i.e., the primal scaling and dual scaling bases are biorthogonal to each other, both multiresolution analyses together form a biorthogonal multiresolution analysis. The orthogonal projection of \eqref{eq:projection} becomes a more general oblique projection,
\begin{equation}\label{eq:dualprojection}
 {\mathcal P}_j f = \sum_{k \in \Z} \langle f, \tilde \phi_{jk} \rangle \, \phi_{jk}.
\end{equation}
Alternatively, with the roles of primal and dual scaling functions interchanged, we also have
\begin{equation*}\label{eq:dualdualprojection}
 \tilde{\mathcal P}_j f = \sum_{k \in \Z} \langle f, \phi_{jk} \rangle \, \tilde \phi_{jk}.
\end{equation*}
All further analysis in this paper is based on the biorthogonal setting. The orthogonal setting corresponds to $V_j=\tilde V_j$ and $\phi=\tilde\phi$.

Next to the primal wavelet basis~\eqref{eq:orthonormalwaveletbasis} we also define a dual wavelet basis
\begin{equation*}\label{eq:dualwavelets}
	\tilde\psi_{jk}(t)=\sqrt 2\sum_{k\in\Z} \tilde  g_k \tilde \phi(2t-k),\qquad \forall k,j\in\Z.
\end{equation*}
To obtain biorthogonal wavelet bases we require the mixed conditions
\begin{gather*}
	\left\langle \phi_{jk},\tilde\phi_{jl}\right\rangle=\delta_{kl},\quad \left\langle \phi_{jk},\tilde\psi_{jl}\right\rangle=0,\quad \left\langle \psi_{jk},\tilde\phi_{jl}\right\rangle=0,\qquad \forall {j,k,l}\in\Z\\
	\left\langle \psi_{ik},\tilde\psi_{jl}\right\rangle=\delta_{kl}\delta_{ij},\qquad\forall {i,j,k,l}\in\Z.
\end{gather*}
One can verify that these conditions follow from a dual double-shift orthogonality and two mixed alternating flip-relations,
\begin{gather*}
	\sum_{k\in\Z}\overline h_k\tilde h_{k+2n}=\delta_n,\qquad \forall n\in\Z\\
	 \overline{g}_k=(-1)^k\tilde h_{1-k},\quad  \tilde g_k=(-1)^k\overline{h}_{1-k},\qquad \forall k\in\Z.
\end{gather*}
In the remainder of the text we will use compactly supported sequences, i.e., sequences $a$ for which there exist $K_1,K_2\in\Z$ such that $a_k=0$ if $k<K_1$ or $k>K_2$. As a result, all associated scaling functions and wavelets have compact support as well.

\subsection{Discrete wavelet transform}\label{ss:dwt}

The sequences $h,\tilde h, g, \tilde g$ that describe the wavelet and scaling bases in the previous section can be used to define the discrete wavelet transform (DWT). We revisit its definition in order to motivate the statements in the complexity analysis of the numerical methods later on. The DWT transforms scaling coefficients of a given function $f(t)\in L^2(\R)$:  $v_{jk}=\left\langle f,\tilde\phi_{jk}\right\rangle$, $j,k\in\Z$ to its wavelet coefficients: $w_{jk}=\left\langle f,\tilde\psi_{jk}\right\rangle$, $j,k\in\Z$.  The inverse discrete wavelet transform (iDWT) transforms wavelet coefficients back into scaling coefficients.

Both the DWT and iDWT are recursive algorithms. In every step, the DWT transforms scaling coefficients at a given level $j+1$ to wavelet and scaling coefficients at a coarser level $j$, while the iDWT recovers in each step the scaling coefficients of the fine level $j+1$ using wavelet and scaling coefficients at level $j$:
\begin{gather}
	v_{jk} = \sum_{l\in\Z}\overline{\tilde h}_{l-2k}v_{j+1,l},\quad
	w_{jk} = \sum_{l\in\Z}\overline{\tilde g}_{l-2k}v_{j+1,l},\qquad \forall j,k\in\Z\label{eq:dwtstep}\\
	v_{j+1,k} = \sum_{l\in\Z} h_{k-2l}v_{jl}+g_{k-2l}w_{jl},\qquad \forall j,k\in\Z\label{eq:idwtstep}.\nonumber
\end{gather}
Usually, the DWT is implemented to transform a finite vector of length $N=2^J$, $\mathbf v_J=\{v_{Jk}\}_{k=0}^{N-1}\in\C^{N}$, to a vector $\mathbf w_J\in\C^{N}$:
\begin{align}
 \mathbf w_J &= [v_{00}, w_{00}, \underbrace{w_{10}, w_{11}}_{\text{2 elements}},\dots,\underbrace{w_{l,0},\dots,w_{l,2^l-1}}_{2^{l}\text{ elements}},\dots, \underbrace{w_{J-1,0},\dots,w_{J-1,2^{J-1}-1}}_{2^{J-1}\text{ elements}} ]\label{eq:waveletvector}\\
	&= [\mathbf{v}_{0}^T,\hat{\mathbf{w}}_{0}^T,\hat{\mathbf{w}}_{1}^T,\dots,\hat{\mathbf{w}}_{J-1}^T]^T\nonumber
\end{align}
with $\hat{\mathbf{w}}_{j} = \{w_{j,k}\}_{k=0}^{2^j-1}$. Boundary conditions  deal with the finite nature of the vectors. We will assume a periodic boundary condition\footnote{Since we intend to employ wavelets on a bounding box $\Xi$ to approximate functions on a subset $\Omega \subset \Xi$, the periodicity of the basis on $\Xi$ is not actually a restriction on $\Omega$, as long as the boundaries of $\Omega$ and $\Xi$ do not touch. One can use other boundary conditions on $\Xi$, but periodicity is the simplest one to implement and manipulate.}, i.e.,
\begin{equation*}\label{eq:periodicboundarycondition}
	v_{jk} = v_{j,k+2^j} ,\quad w_{jk} = w_{j,k+2^j}\qquad \forall j\in\Z^+,\forall k\in\Z.
\end{equation*}

To transform $\mathbf v_J$ into $\mathbf w_J$ the DWT performs $J$ steps like~\eqref{eq:dwtstep}. If we use $[A]_{\downarrow q}$ to denote the down-sampling of a matrix, i.e., the selection of every $q$th row $\left([A]_{\downarrow q}\right)(k,l) = A(qk,l)$ and use $A^*$ to denote the adjoint of $A$, one step of the DWT can be represented in matrix notation as
\begin{equation*}
	\begin{bmatrix}
	\mathbf v_{j-1}\\\hat{\mathbf{w}}_{j-1}
	\end{bmatrix} =
	\begin{bmatrix}
	{\tilde{H}_j^*}\\{\tilde{G}_j^*}
	\end{bmatrix}_{\downarrow 2} \mathbf v_{j}
\end{equation*}
with matrix $\tilde{H}_j\in\C^{2^j\times 2^j}$
\begin{align*}
	\tilde H_j(k,l) = \sum_{m\in\Z}\tilde h_{k-l + m2^j}, \qquad\forall k,l=0,\dots,2^j-1,
\end{align*}
such that $[\tilde{H}_j^*]_{\downarrow 2}\in\C^{2^{j-1}\times 2^j }$ and $[\tilde{H}_j^*]_{\downarrow 2}(k,l)= \sum_{m\in\Z}\tilde h_{l-2k + m2^j}$.
Note that the summation over $m$ here is used to incorporate the periodic boundary conditions. The matrices $H_j,G_j$~and~$\tilde G_j$ are defined analogously.

The full DWT in matrix notation is $\mathbf w_J = W_J\mathbf v_J$, with $W_J\in\C^{N\times N}$ and
\begin{equation*}\label{eq:fulldwt}
	W_J=\underbrace{
	\begin{bmatrix}
	\begin{bmatrix}
	{\tilde{H}_1^*}\\{\tilde{G}_1^*}
	\end{bmatrix}_{\downarrow 2}&0\\
	0&I_{2^{J}-2}
	\end{bmatrix}
	\begin{bmatrix}
	\begin{bmatrix}
	{\tilde{H}_2^*}\\{\tilde{G}_2^*}
	\end{bmatrix}_{\downarrow 2}&0\\
	0&I_{2^{J}-2^{2}}
	\end{bmatrix}
	\cdots
	\begin{bmatrix}
	\begin{bmatrix}
	{\tilde{H}_{J-1}^*}\\{\tilde{G}_{J-1}^*}
	\end{bmatrix}_{\downarrow 2}&0\\
	0&I_{2^{J-1}}
	\end{bmatrix}
	\begin{bmatrix}
	{\tilde{H}_J^*}\\{\tilde{G}_J^*}
	\end{bmatrix}_{\downarrow 2}
}_{J\text{ terms}}.
\end{equation*}
Similarly, the full iDWT can be written as $\mathbf v_J = W^{-1}_J\mathbf w_J$. The iDWT matrix is the inverse of $W_J$, $W^{-1}_J$. This inverse can be decomposed as
\begin{equation*}\label{eq:fullidwt}
W^{-1}_J=\underbrace{
	\begin{bmatrix}
	{{H}_J^*}\\{{G}_J^*}
	\end{bmatrix}_{\downarrow 2}^*
	\begin{bmatrix}
	\begin{bmatrix}
	{{H}_{J-1}^*}\\{{G}_{J-1}^*}
	\end{bmatrix}_{\downarrow 2}^*&0\\
	0&I_{2^{J-1}}
	\end{bmatrix}
	\cdots
	\begin{bmatrix}
	\begin{bmatrix}
	{{H}_2^*}\\{{G}_2^*}
	\end{bmatrix}_{\downarrow 2}^*&0\\
	0&I_{2^{J}-2^{2}}
	\end{bmatrix}
	\begin{bmatrix}
	\begin{bmatrix}
	{{H}_1^*}\\{{G}_1^*}
	\end{bmatrix}_{\downarrow 2}^*&0\\
	0&I_{2^{J}-2}
	\end{bmatrix}
}_{J\text{ terms}}
\end{equation*}
where subsampling takes precedence over taking the adjoint in order to avoid a multitude of brackets.

\begin{figure}
	\centering
	{\includegraphics[width=.4\columnwidth,trim={0em -1.6em 0em 0em},]{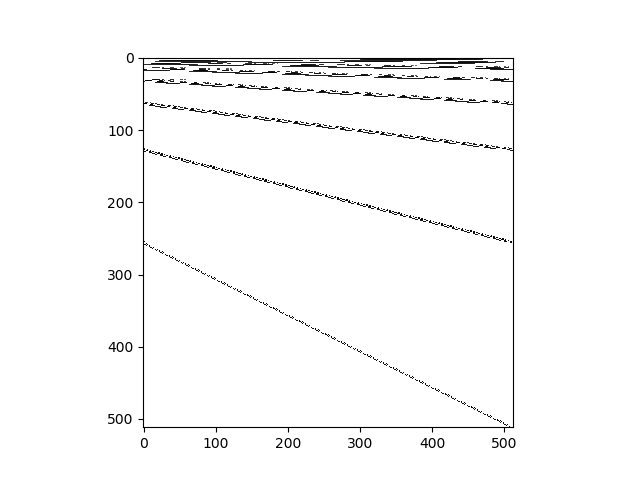}}
	{\includegraphics[width=.4\columnwidth,trim={0em -1.6em 0em 0em},]{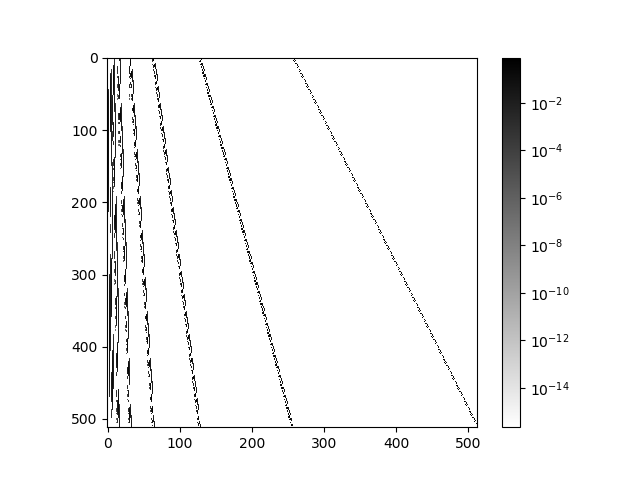}}
	\caption{Structure of DWT (left) and iDWT matrix (right) for \texttt{db2} and $J=9$. \label{fig:dwtstructure}.}
\end{figure}

With the use of cascading filter banks the DWT and iDWT can be implemented in $\mathcal O(N)$ operations \cite[Section~7.4.1: Fast Biorthogonal Wavelet Transform, p. 310]{Mallat2009}. This algorithm is called the fast wavelet transform (FWT) and was introduced in \cite{Mallat1989}.
The same complexity can not be achieved using an ordinary matrix-vector multiply since $W$ and $W^{-1}$ contain $\mathcal O(N\log (N))$ non-zero elements. This is clear by looking at Figure~\ref{fig:dwtstructure} and by the following lemma.
\begin{lemma}\label{lem:dwtstructure}
	Each column of $W_J$ has $\mathcal O(J)$ non-zero elements while each row of its inverse has $\mathcal O(J)$ non-zero elements. Furthermore, both have $\mathcal O(J2^J)$ non-zero elements.
\end{lemma}
\begin{proof}
	Owing to our periodic setting, in the following we say that a matrix $A\in\R^{N\times N}$ is banded if
	\[
	|m-n|\mod N>b\Rightarrow A(m,n)=0.
	\]
	Let $A\in\R^{N\times N}$ have bandwidth $a$ and $B\in\R^{2N\times 2N}$ have bandwidth $b$, then $A[B]_{\downarrow 2}$ has bandwidth $a+b/2$. This is verified by writing
	\[ \left(A[B]_{\downarrow 2}\right)(m,n) = \sum_{i=0}^{N-1}A(m,i)B(2i,n), \]
	which is only non-zero if $|m-i|\mod N\leq a$ and $|2i-n|\mod 2N\leq b$, i.e, if $|m-n|\leq a+b/2$.

	We can rewrite $W$ in $J$ vertical blocks
	\begin{equation*}
		W = \begin{bmatrix}
			B_1\\B_2\\\\\vdots\\\\B_{J-2}\\B_{J-1}\\B_J
			\end{bmatrix} = \begin{bmatrix}
			\begin{bmatrix}\tilde H_{1}^*\\\tilde G_{1}^*\end{bmatrix}_{\downarrow 2}[\tilde H_{2}^*]_{\downarrow 2}\cdots[\tilde H_{J-1}^*]_{\downarrow 2} [\tilde H_J^*]_{\downarrow 2}\\
			[\tilde G_{2}^*]_{\downarrow 2}[\tilde H_{3}^*]_{\downarrow 2}\cdots[\tilde H_{J-1}^*]_{\downarrow 2} [\tilde H_J^*]_{\downarrow 2}\\\\
			\vdots\\\\
			[\tilde G_{J-2}^*]_{\downarrow 2}[\tilde H_{J-1}^*]_{\downarrow 2}[\tilde H_J^*]_{\downarrow 2}\\
			[\tilde G_{J-1}^*]_{\downarrow 2}[\tilde H_J^*]_{\downarrow 2}
			\\
			[\tilde G_J^*]_{\downarrow 2}
		\end{bmatrix}.
	\end{equation*}
The blocks $B_1$ and $B_2$ have size $2\times 2^J$ and block $B_j$ has size $2^{j-1}\times 2^J$ for $j\geq1$. First we show that each block has a bounded number of non-zero elements per column. To that end we denote $K$ and $L$ as the number of non-zero elements in the sequences $\tilde h$ and $\tilde g$ respectively. It is clear that $[\tilde H_j^*]_{\downarrow 2}$ and $[\tilde G_j^*]_{\downarrow 2}$ have $\mathcal O(K/2)$ and $\mathcal O(L/2)$ non-zero elements per column respectively since $H_J$ and $G_J$ are banded with bandwidth $K$ and $L$ respectively. The products also have a bounded number of non-zero elements per column by the first part of the proof. Matrix $[\tilde G_{J-1}^*]_{\downarrow 2}[\tilde H_J^*]_{\downarrow 2}$ has e.g. $L/2+K/4$ non-zero elements per column.

	Since each block contains $\mathcal O((K+L)2^J)$ non-zero elements with a limited number of non-zero per column, $W$ contains $\mathcal O(2^J J)$ non-zero elements and $\mathcal O(J)$ non-zero elements per column. The proof for $W^{-1}$ is entirely analogous; but write $W^{-1}$ in $J$ horizontal blocks.
\end{proof}

If we introduce the dual DWT $\tilde W_J$
\begin{equation*}\label{eq:dualdwt}
\tilde W_J=\underbrace{
	\begin{bmatrix}
	\begin{bmatrix}
	{{H}_1^*}\\{{G}_1^*}
	\end{bmatrix}_{\downarrow 2}&0\\
	0&I_{2^{J}-2}
	\end{bmatrix}
	\begin{bmatrix}
	\begin{bmatrix}
	{{H}_2^*}\\{{G}_2^*}
	\end{bmatrix}_{\downarrow 2}&0\\
	0&I_{2^{J}-2^{2}}
	\end{bmatrix}
	\cdots
	\begin{bmatrix}
	\begin{bmatrix}
	{{H}_{J-1}^*}\\{{G}_{J-1}^*}
	\end{bmatrix}_{\downarrow 2}&0\\
	0&I_{2^{J-1}}
	\end{bmatrix}
	\begin{bmatrix}
	{{H}_J^*}\\{{G}_J^*}
	\end{bmatrix}_{\downarrow 2}
}_{J\text{ terms}}
\end{equation*}
and analogously define the iDWT $\tilde W^{-1}_J$, we can verify that
\begin{equation}\label{eq:dwtrelations}
	W_J^*=\tilde W_J^{-1}\qquad(W_J^{-1})^*=\tilde W_J.
\end{equation}

\section{Periodic wavelets on the interval and discrete duals}\label{wavs:periodicwavelet}

\subsection{Periodization}

For simplicity of the exposition we again assume that $N=2^J$, $J\in\N$. We introduce the periodic and scaled father function with period $1$ as
\begin{equation*}\label{eq:periodizedfatherfunction}
	\phi_{N}(t) = \sum_{k\in\Z} 2^{J/2}\phi\left(2^{J}(t-k)\right).
\end{equation*}
The periodic scaling basis that consists of $N$ translated father functions $\phi_{kN}(t)=\phi_{N}(t-\tfrac kN)$ is
\begin{equation*}
	\Phi_{N} = \{ \phi_{kN} \}_{k=0}^{N-1}.
\end{equation*}
Similarly, the dual scaling function in $L^2(0,1)$ is
\begin{equation*}
\tilde\Phi_{N} = \{ \tilde\phi_{kN} \}_{k=0}^{N-1}.
\end{equation*}
with translates $\tilde\phi_{kN}(t)=\tilde\phi_{N}(t-\tfrac kN)$ of the dual father function
\[
\tilde\phi_{N}(t) = \sum_{k\in\Z} 2^{J/2}\tilde\phi\left(2^{J}(t-k)\right).
\]

Wavelet bases with period $1$ can be devised in the same way as above for the scaling bases, by summing over their translations. Alternatively, we can define them using periodic scaling bases and the iDWT matrices $W^{-1}_J$ and $\tilde W^{-1}_J$:
\begin{align}\label{eq:primalwaveletbasis}
	\Psi_{N} = \{\psi_{kN}\}_{k=0}^{N-1},\qquad \psi_{kN}(t)=\sum_{l=0}^{N-1}\phi_{lN}(t)(W^{-1}_J)(l,k),\\
	\tilde\Psi_{N} = \{\tilde\psi_{kN}\}_{k=0}^{N-1},\qquad \tilde\psi_{kN}(t)=\sum_{l=0}^{N-1}\tilde\phi_{lN}(t)(\tilde W^{-1}_J)(l,k)\label{eq:dualwaveletbasis}.
\end{align}
Note that the index $k$ in $\psi_{kN}$ incorporates both the scale and the translation of the corresponding wavelet in this notation. The $N$ wavelets are defined on all scales $0,1,\ldots,J-1$, as in \eqref{eq:waveletvector}.

\subsection{Compactly supported wavelets and discrete evaluation}

\begin{figure}
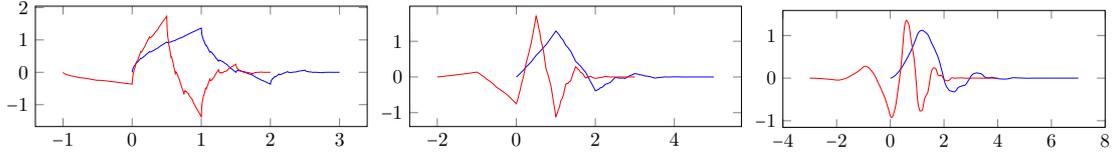

	\centering
	\resizebox{.33\columnwidth}{!}{\input{img/db2.tikz}}%
	\resizebox{.33\columnwidth}{!}{\input{img/db3.tikz}}%
	\resizebox{.33\columnwidth}{!}{\input{img/db4.tikz}}%
	\caption{The father (blue) and mother (red) function of \texttt{db2}, \texttt{db3}, and \texttt{db4} (left to right). \label{fig:daubechieswavelets}}
\end{figure}

We focus on two well-known families of compactly supported wavelet families. The first family are the Daubechies orthogonal wavelet bases \cite{Daubechies1988}. It was shown in \cite{Daubechies1988} that the scaling function for orthogonal wavelets with $p$ vanishing moments (a regularity condition on the wavelets) has a support of length at least $2p-1$. Daubechies wavelets are optimal in the sense that they have a minimum support length for a given number of vanishing moments \cite{Daubechies1988}\cite[Theorem 7.9]{Mallat2009}. In Figure~\ref{fig:daubechieswavelets}, $\phi$ and $\psi$ are shown for \texttt{db2}, \texttt{db3} and \texttt{db4}, i.e., the  father and mother functions of the Daubechies wavelet with 2, 3 and 4 vanishing moments. These standard wavelets are widely used in applications.

The Daubechies wavelet and scaling functions are defined by their compactly supported sequences $h$ and $g$. No closed form formula is known for the functions themselves. Note that the functions associated with \texttt{db2} are continuous, but nowhere differentiable. However, one can evaluate compactly supported scaling functions at dyadic points $k/2^j$, $k\in\Z$ from the sequence $h$ using the following procedure  \cite{Daubechies1992a}.

We create a matrix system by evaluating the two-scale relation~\eqref{eq:twoscalerelation} in integer points. For a case where $\phi(k)=0$ for $k<0$ or $k>5$ we obtain for example
\begin{align*}
	\begin{bmatrix}
		\phi(0)\\
		\phi(1)\\
		\phi(2)\\
		\phi(3)\\
		\phi(4)\\
		\phi(5)
	\end{bmatrix}
	 =
	 \sqrt 2 \begin{bmatrix}
	 h_0&&&&&\\
	 h_2&h_1&h_0&&&\\
	 h_4&h_3&h_2&h_1&h_0&\\
	 &h_5&h_4&h_3&h_2&h_1\\
	 &&&h_5&h_4&h_3\\
	 &&&&&h_5\\
	 \end{bmatrix}
	 \begin{bmatrix}
	 \phi(0)\\
	 \phi(1)\\
	 \phi(2)\\
	 \phi(3)\\
	 \phi(4)\\
	 \phi(5)
	 \end{bmatrix}.
\end{align*}
In this notation, the two-scale relation implies that the matrix shown has an eigenvalue $1$. The corresponding eigenvector represents $\phi$ evaluated at integer points. To evaluate at finer dyadic levels, it suffices to use the two-scale relation repeatedly, since it also states that
\begin{equation*}
	\phi\left(\tfrac k {2^{j+1}}\right) = \sum_{l\in\Z}h_l\phi\left(\tfrac k {2^j}-l\right).
\end{equation*}
The wavelet function may be evaluated at the dyadic points by first evaluating the corresponding scaling function, after which equation~\eqref{eq:orthonormalwaveletbasis} can be applied.

\begin{figure}
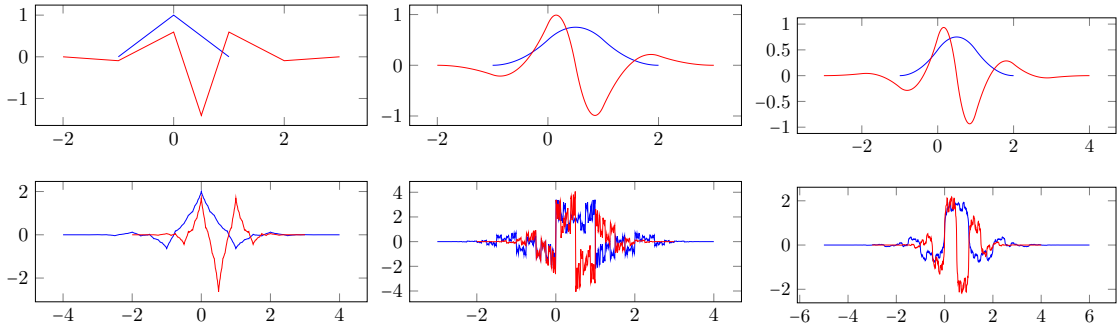

	\centering
	\resizebox{.33\columnwidth}{!}{\input{img/cdf24.tikz}}%
	\resizebox{.33\columnwidth}{!}{\input{img/cdf33.tikz}}%
	\resizebox{.33\columnwidth}{!}{\input{img/cdf35.tikz}}%
	\caption{The father (blue) and mother (red) function of \texttt{cdf24}, \texttt{cdf33}, and \texttt{cdf35} (left to right) both primal (top) and dual (bottom).}\label{fig:cdfwavelet}
\end{figure}

The second wavelet family we consider are the biorthogonal CDF wavelets. These are compactly supported and symmetric \cite{Cohen1992}. More specifically, among the family of CDF wavelets we use those that have the centered B-spline as primal father function \cite[\S 6.A]{Cohen1992}. All sequences $h$, $g$, $\tilde h$ and $\tilde g$ are compactly supported and symmetric as well. Figure~\ref{fig:cdfwavelet} shows the primal and dual father and mother functions for \texttt{cdf24}, \texttt{cdf33} and \texttt{cdf35}. The first digit in the name indicates the number of vanishing moments of the dual wavelet, the second digit refers to those of the primal wavelet. Contrary to the Daubechies wavelets, closed-form formulas exist for some primal scaling functions. The primal scaling functions of \texttt{cdf}$p\tilde p$ shown in Figure~\ref{fig:cdfwavelet} are the centered B-splines of order $p$. However, the dual scaling function can typically only be evaluated in dyadic points using the above procedure.\footnote{The values of $p$ and $\tilde{p}$ do not uniquely determine a biorthogonal multiresolution analysis. We use filters corresponding to B-splines as described in \cite[\S 6.A]{Cohen1992}. Thus, the meaning of \texttt{cdf44} in this paper differs from the widely used CDF filters of primal and dual order $4$ in signal processing, e.g., in the JPEG2000 standard.}

\subsection{Continuous dual bases}

Both the Daubechies and CDF wavelets have at least one compact dual in the continuous sense, i.e., the dual is biorthogonal in $L^2(\R)$:
\begin{equation*}
	\langle\phi_{jk},\tilde\phi_{jl}\rangle_{L^2(\R)}=\int_{\R}\phi_{jk}(t)\tilde\phi_{jl}(t)\d t=\delta_{kl},\qquad \forall j,k,l\in\Z.
\end{equation*}
The dual scaling function is compactly supported as well. So, a compactly supported basis biorthogonal to $\Phi_N$ ($\Psi_N$) in the continuous sense is $\tilde\Phi_N$ ($\tilde\Psi_N$). Note that there might exist multiple dual bases for the same primal basis.

\subsection{Discrete dual scaling bases}

The various dual bases of CDF are well understood. In the context of this paper, for the efficient computation of wavelet approximations on irregular domains we will require a notion of discrete duality. We therefore introduce the following bilinear form:
\begin{equation}\label{eq:discreteinnerproduct}
	\langle f,g\rangle_{\osi} = \sum_{m\in\Z} f\left(\tfrac mq\right)g\left(\tfrac mq\right).
\end{equation}
Here, $\osi \in \mathbb{N}$, $\osi \geq 1$, plays the role of an oversampling factor.

We aim for a periodic dual scaling basis $\tilde\Phi^\osi_N=\{\tilde\phi^\osi_{kN}\}_{k=0}^{N-1}$ defined on $[0,1]$. In order to construct such duals, we focus first on the samples of the father function on the whole real line. Define the sequence $b$ by sampling $\phi$ in the oversampled grid,
\begin{equation*}
	b^\osi_m=\phi\left(\tfrac m\osi\right), \qquad m\in\Z.
\end{equation*}
The sequence $b$ is compactly supported because $\phi$ is compactly supported. Biorthogonality with respect to~\eqref{eq:discreteinnerproduct} between the integer shifts of $\phi$ and those of a discrete dual father function $\tilde \phi^\osi(t)$, with samples $\tilde b$, translates into the conditions
\begin{equation}\label{eq:discretebiorthogonality}
		\langle \phi, \tilde \phi^\osi(\cdot-k) \rangle_\osi = \sum_{m}b_m\tilde b_{m-k\osi} = \delta_{0k}, \qquad \forall k\in\Z.
\end{equation}
Note that shifts by an integer $k$ of a continuous function correspond to shifts of $k\osi$ samples of its sampled sequence in the discrete grid, because $\osi$ is the oversampling factor.

Once a dual sequence $\tilde{b}$ satisfying \eqref{eq:discretebiorthogonality} is found, it does not immediately give rise to a continous representation of the dual functions $\tilde\phi^\osi(t)$. However, we do know its evaluations in the points $\tfrac m\osi$. We can define suitable discrete periodized dual functions as
\begin{equation}\label{eq:discreteperiodicdual}
 \tilde \phi^\osi_{kN}\left(\tfrac mq \right) =  N^{-1/2} \sum_{l \in \Z} \tilde{b}_{m-\osi k-N\osi l}, \qquad m=0,\ldots,Nq-1.
\end{equation}
Note that there are $Nq$ samples in $[0,1)$, and that the summation over $l$ introduces periodization. By construction, these functions satisfy the discrete biorthogonality conditions
\begin{equation*}\label{eq:discrete_dual_biorthogonality}
 \langle \tilde\phi^\osi_{kN},\phi_{lN}\rangle_{N, \osi}=\delta_{kl},
\end{equation*}
where $\langle \cdot, \cdot \rangle_{N,q}$ is a scaled analogue of \eqref{eq:discreteinnerproduct} restricted to $[0,1]$:
\begin{equation}\label{eq:discreteinnerproduct_N}
	\langle f,g\rangle_{N,\osi} = \sum_{m = 0}^{Nq-1} f\left(\tfrac m{Nq}\right)g\left(\tfrac m{Nq}\right).
\end{equation}

%\MAGENTA{Ik heb de periodieke discrete functies hier geprobeerd iets formeler te defini{\"e}ren. De formule is nog niet volledig. De bedoeling is om met $l$ te sommeren over veelvouden van $1$ (periodizatie), en om de juiste samples van $\tilde b$ op de juiste plaats te krijgen ifv $k$ en $N$ en $m$. Is dit zinvol? Er moet ook een schalering zijn ifv $N$ denk ik. Het doel is dan de volgende projectie:}\GREEN{Ik heb de scaleringen en indices ingevuld, na uitschrijven op papier. Zou moeten kloppen. Ik heb  ondertussen enkele fouten gehaald uit eerdere vergelijken. Het is zeker zinvol in het begrijpen wat te betekenis van deze duale is. }

The discrete dual leads to the discrete projections
\begin{equation*}\label{eq:discrete_projection}
 {\mathcal P}^\osi_N f(t) = \sum_{k=0}^{N-1} \langle f,  \tilde \phi^\osi_{kN} \rangle_{N,\osi} \, \phi_{kN}(t) = \sum_{k=0}^{N-1} v_{kN}^\osi\,   \phi_{kN}(t).
\end{equation*}
The discrete duals play the role of the continuous dual $\tilde \phi_{kN}$ in \eqref{eq:dualprojection}. However, the discrete inner product does not actually require the evaluation of integrals. By construction, the projection is exact on the span of the scaling functions:
\[
 {\mathcal P}^\osi_N f = f, \quad \forall f \in \SPAN \Phi_N. 
\]
In the terminology of splines literature, these reproducing projections are examples of a \emph{quasi-interpolation} method, in which the global approximation is constructed using local approximations \cite{deboor1985controlled,chui1984quasi,dahmen1984order}.
%\GREEN{Is de standaard niet q=1? Ken jij referenties van quasi-`interpolatie` die q>1 gebruiken?} 

\begin{figure}
	\centering
	\resizebox{.33\columnwidth}{!}{\input{img/discretedb24.tikz}}%
	\resizebox{.33\columnwidth}{!}{\input{img/discretedb34.tikz}}%
	\resizebox{.33\columnwidth}{!}{\input{img/discretedb44.tikz}}%
	\caption{The sampled father function (top) and a compact discrete dual (bottom) for \texttt{db2}, \texttt{db3}, \texttt{db4} with $\osi=4$.}\label{fig:compactduals_db}
\end{figure}
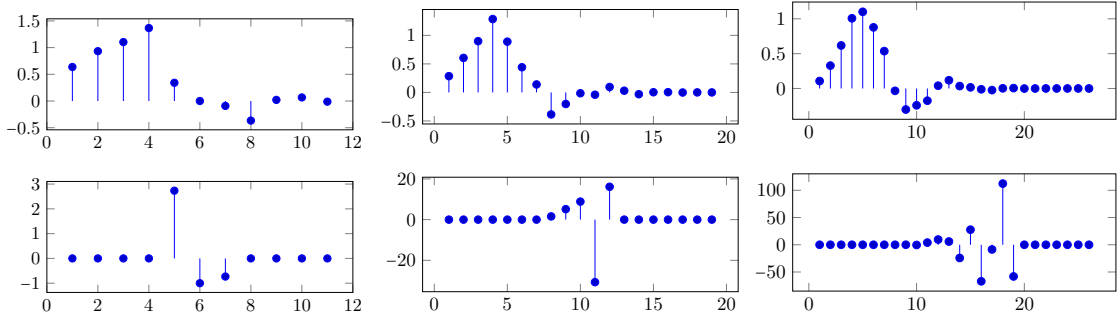
\begin{figure}
	\centering
	\resizebox{.33\columnwidth}{!}{\input{img/discretecdf314.tikz}}%
	\resizebox{.33\columnwidth}{!}{\input{img/discretecdf424.tikz}}%
	\resizebox{.33\columnwidth}{!}{\input{img/discretecdf514.tikz}}%
	\caption{The sampled father function (top) and a compact discrete dual (bottom) for \texttt{cdf31}, \texttt{cdf42}, \texttt{cdf51} with $\osi=4$}\label{fig:compactduals_cdf}
\end{figure}
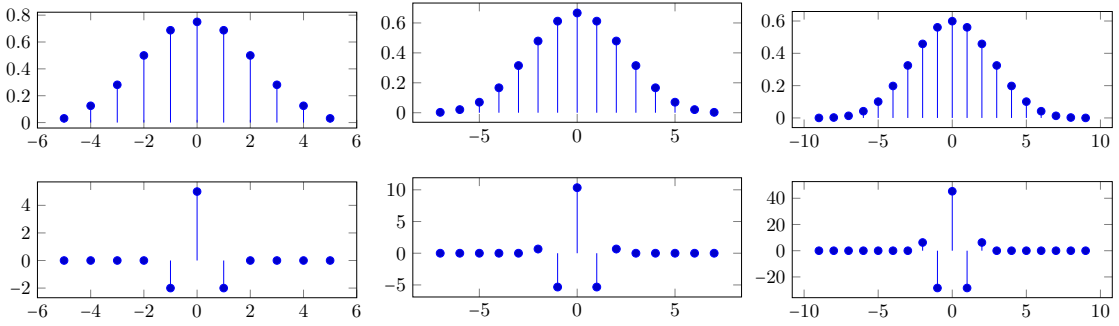

Unlike the Daubechies and CDF scaling functions and their continuous duals, the discrete duals we employ are non-standard. However, it was recently shown by the authors in \cite[Theorem 3.6]{coppe2019splines} that discrete compact dual sequences $\tilde b$ always exist for a basis consisting of translates of B-splines, which is exactly the setting of our choice of CDF primal scaling bases. Experiments show that compact duals of Daubechies scaling functions can also be found. Figures~\ref{fig:compactduals_db} and~\ref{fig:compactduals_cdf} show some of these compact duals, namely those with the smallest possible support. In Table \ref{tab:discretevalues} of the appendix we give the values of the primal and dual sequences used in the experiments below. We observe that the discrete duals for the Daubechies case are somewhat larger than those of the CDF duals. These discrete duals are readily found by solving the linear algebraic conditions~\eqref{eq:discretebiorthogonality}, noting that the system is finite because $b$ and $\tilde{b}$ are compactly supported.

The size of the discrete duals can be reduced by allowing for a larger support and solving the corresponding conditions~\eqref{eq:discretebiorthogonality} in a least squares sense, thereby minimizing the norm of the solution to an underdetermined system. The relevance of the size of the discrete duals is formalized in the following lemma.
\begin{lemma}\label{lem:discrete_error}
 Let $f \in L^\infty([0,1])$ and let the synthesis operator $T_N a = \sum_{k=0}^{N-1} a_k \, \phi_{kN}$ of $\Phi_N$ be a bounded operator from $\C^N$ to $L^2([0,1])$ with bound $B > 0$, i.e.,
 \[
 \Vert T_N \mathbf{a} \Vert_{L^2([0,1])} =  \left\Vert \sum_{k=0}^{N-1} a_k \phi_{kN} \right\Vert_{L^2([0,1])} \leq B \, \Vert \mathbf{a} \Vert_{\ell_2(\C^N)}.
 \]
 Define $\Vert \cdot \Vert_{N,q}^2 = \langle \cdot, \cdot \rangle_{N,q}$. If the support of the sequence $\tilde{b}$ is less than $Nq$, then
 \begin{equation}\label{eq:discrete_error}
 \Vert {\mathcal P}_N f - \tilde{{\mathcal P}}^\osi_N f \Vert \leq B \, \Vert f - {\mathcal P}_N f \Vert_{N,q}  \, \Vert \tilde{b} \Vert_{\ell_2}.
 \end{equation}
\end{lemma}
\begin{proof}
%Using~\eqref{eq:discreteperiodicdual} and~\eqref{eq:discreteinnerproduct_N}, we find for $u \in L^\infty([0,1])$ the following inequality:
%\begin{align*}
%\left| \langle u, \tilde \phi^\osi_{kN} \rangle_{N,\osi} \right| &= \left|\sum_{m=0}^{Nq-1} u\left(\tfrac {m}{\osi N}\right) \, \tilde{\phi}^\osi_{kN} \left(\tfrac {m}{\osi N}\right)\right| \\
%&\leq \left\Vert u\left(\tfrac {\cdot}{\osi N}\right) \right\Vert_{\ell_2(\C^{Nq})} \, \left\Vert \tilde{\phi}^\osi_{kN} \left(\tfrac {\cdot}{\osi N}\right) \right\Vert_{\ell_2(\C^{Nq})} \\
%&= N^{-1/2} \, \Vert u \Vert_{N,q} \, \Vert \tilde{b} \Vert_{\ell_2}.
%\end{align*}
%\GREEN{\RED{is niet duidelijker alles in de $N,q$ `norm' te houden, wat inderdaad equivalent is met de $\ell_2(\C^{Nq})$ norm op de samples van de functies?}
Using~\eqref{eq:discreteinnerproduct_N}, Cauchy-Schwartz and~\eqref{eq:discreteperiodicdual}, we find for $u \in L^\infty([0,1])$ the following inequality:
\begin{align*}
\left| \langle u, \tilde \phi^\osi_{kN} \rangle_{N,\osi} \right| &\leq \| u\|_{N,\osi}\|\tilde \phi^\osi_{kN} \|_{N,\osi}\\
&= N^{-1/2} \, \Vert u \Vert_{N,q} \, \Vert \tilde{b} \Vert_{\ell_2}.
\end{align*}
In the last line, we have used the assumption on the discrete support of $\tilde{b}$ to note that the periodic copies of $\tilde{b}$ in \eqref{eq:discreteperiodicdual} do not actually overlap, hence the equality of norms.

Next, let $g_N = f - {\mathcal P}_N f = f - f_N$, such that $f = f_N + g_N$.  Since both ${\mathcal P}_N$ and $\tilde{{\mathcal P}}_N$ are exact on the span of $\Phi_N$, we have that ${\mathcal P}_N f_N = \tilde{{\mathcal P}}_N f_N = f_N$. In addition, ${\mathcal P}_N g_N = 0$. With norms in $L^2([0,1])$ unless noted otherwise, this means that
\begin{align*}
 \Vert {\mathcal P}_N f - \tilde{{\mathcal P}}^\osi_N f \Vert &= \Vert {\mathcal P}_N g_N - \tilde{{\mathcal P}}^\osi_N g_N \Vert = \left\Vert \tilde{{\mathcal P}}^\osi_N g_N \right\Vert = \left\Vert \sum_{k=0}^{N-1} \langle g_N, \tilde \phi^\osi_{kN} \rangle_{N,\osi} \, \phi_{kN} \right\Vert \\
 &\leq B \, \left\Vert \left\{\langle g_N,  \tilde \phi^\osi_{kN} \rangle_{N,\osi}\right\}_{k=0}^{N-1}\, \,  \right\Vert_{\ell_2(\C^N)} \\
 &\leq B \, \sqrt{N} \left\Vert \left\{\langle g_N,  \tilde \phi^\osi_{kN} \rangle_{N,\osi}\right\}_{k=0}^{N-1}\, \,  \right\Vert_{\ell_\infty(\C^N)} \\
 &\leq B \, \Vert g_N \Vert_{N,q} \, \Vert \tilde{b} \Vert_{\ell_2}.
\end{align*}
\end{proof}
The lemma shows that the discrete projection yields a result that is close to the actual continuous dual projection, as long as the dual sequence $\tilde{b}$ does not grow too large in norm. In particular, the convergence rate with $N$ remains the same, only the constant factor is affected. In practice, the assumption on the support of $\tilde{b}$ is not restrictive unless $N$ is very small.

\subsection{Discrete dual wavelet bases?}

We have replaced the continuous dual basis with a discrete one. In view of the close correspondence between the continuous inner product coefficients $v_{kN}= \langle f, \tilde{\phi}_{Nk} \rangle$ and their discrete analogues $v_{kN}^\osi = \langle f, \tilde{\phi}^\osi_{Nk} \rangle_{N,\osi}$, as quantified by Lemma~\ref{lem:discrete_error}, we forego the construction of a discrete dual wavelet basis. We simply retain the primal and dual scaling functions and wavelets, and we continue to use the wavelet transform given by $W_j$ and $\tilde{W}_j$ respectively. Thus, the construction involving the discrete sequence $\tilde{b}$ of the previous section can be thought of merely as a quadrature scheme to approximate $v_{kN}$ using equispaced samples of $f$. Compared to other quadrature schemes for wavelet and scaling coefficients, such as Sweldens quadrature \cite{sweldens1994quadrature}, our quadrature discretization has an additional discrete orthogonality structure that will be used later on.

The construction of a fully discrete wavelet basis based on $\tilde{\phi}^\osi_{kN}$ would have several disadvantages. First of all, the sequence $\tilde{b}$ has no multiscale structure. Therefore, $\tilde{\phi}^\osi_{kN}$ does not satisfy a two-scale relation and further structure would have to be imposed on $\tilde{b}$. Secondly, and more importantly, changing the dual scaling function from $\tilde{\phi}$ to $\tilde{\phi}^\osi_{kN}$ would affect the primal wavelet. Finally, there is no guarantee that a discrete dual wavelet basis exists for a given primal scaling function. We choose to retain the Daubechies and CDF family of scaling functions and wavelets.

Note that the number of dual vanishing moments apparent as $\tilde p$ in the naming of the CDF wavelets \texttt{cdf}$p\tilde p$ has no influence on the regularity of the discrete duals $\tilde{\phi}^\osi$ defined above. It does change the shape of the primal wavelet basis elements through the alternating-flip relation $g_k=(-1)^k\overline{\tilde{h}}_{1-k}$. The higher $\tilde p$, the larger the support of the wavelets. Another effect of the number of dual vanishing moments is seen in the wavelet transforms. If $\tilde p < p$ then the norm of $\tilde W_J$ may increase significantly. This is illustrated in Table \ref{tab:dwtnorms}. For the purposes of this paper, it seems best to consider $\tilde p \geq p$.

\begin{table}[h]
	\centering
	\input{dwttabular1013.tikz}
	\input{dwttabular1046.tikz}
	\caption{\label{tab:dwtnorms}Norms of CDF DWT matrices $W_J$ and $W^{-1}_J$ with $J=10$.}
\end{table}

\section{The approximation problem}\label{wavs:approximationproblem}

We formally define the approximation problem. To that end, we introduce notation for the multivariate approximation. We largely adopt the same notation as was used for B-splines \cite{coppe2019splines}, but we replace the spline bases with a tensor product of the wavelet bases~\eqref{eq:primalwaveletbasis}. Bold letters such as $\mathbf N$ denote a vector of length $d$, $\mathbf N=(N_1,\dots N_d)$, and $I_{\mathbf N}$ denotes the index set
\[
\{ (i_1,\dots,i_d)\, |\, j=1\dots,d,\,  i_j=0,\dots,N_j-1\}.
\]
Furthermore, the tensor product of wavelet bases \eqref{eq:primalwaveletbasis} is written as
\begin{align*}
	\Psi_{\mathbf N}=\Psi_{N_1}\otimes\cdots\otimes\Psi_{N_d}=\{\psi_{\mathbf k\mathbf N}\}_{\mathbf k\in I_{\mathbf N}}\qquad \psi_{\mathbf k\mathbf N}(\mathbf t) = \Psi_{k_1N_1}(t_1)\times\cdots\times\Psi_{k_dN_d}(t_d).
\end{align*}
The tensor products of scaling bases and dual bases are denoted analogously.

\subsection{Continuous projection}

The approximation problem can be discretized in two ways, namely, using inner products and point evaluation. The former is denoted the \emph{continuous projection} and leads to the system $A_{\mathbf N}x=b_{\mathbf N}$, with
\begin{align}
	A_{\mathbf N}(\mathbf k,\mathbf l) &= \langle \psi_{\mathbf l\mathbf N}, \tilde \psi_{\mathbf k\mathbf N}\rangle_{\Omega}, \qquad \mathbf k, \mathbf l\in I_\mathbf N,\label{eq:contsystem}\\
	b_{\mathbf N}(\mathbf k) &= \langle f,\tilde\psi_{\mathbf k\mathbf N}\rangle_{\Omega}\nonumber.
\end{align}
Note that the inner products are defined over $\Omega$, but the wavelet basis has been defined on $\Xi$. If $\Omega=\Xi$ then $A_{\mathbf N}$ is the identity matrix due to the continuous duality of $\Psi_{\mathbf N}$ and $\tilde\Psi_{\mathbf N}$. In our setting $\Omega \subset \Xi$, hence the matrix entries of $A_{\mathbf N}$ may differ from $0/1$ if one of the functions involved overlaps with the boundary. 
It may be difficult to evaluate the inner products on $\Omega$ numerically, especially in the multivariate setting.

We explicitly define the wavelet and scaling basis elements that overlap with the boundary of a given domain $\Omega$. They are contained in the sets
\begin{equation}
\label{eq:splineboundaryset}
\K_{\mathbf N}(\Omega) = \{ {\mathbf k} \in I_{\mathbf N} \, | \, \, \mysupp \phi_{\mathbf k\mathbf N} \cap \Omega \neq \varnothing\text{ and }\mysupp \phi_{\mathbf k\mathbf N} \cap {\Omega^c} \neq \varnothing \}
\end{equation}
and
\begin{equation}
\label{eq:waveletboundaryset}
\L_{\mathbf N}(\Omega) = \{ {\mathbf k} \in I_{\mathbf N} \, | \, \, \mysupp \psi_{\mathbf k\mathbf N} \cap \Omega \neq \varnothing\text{ and }\mysupp \psi_{\mathbf k\mathbf N} \cap {\Omega^c} \neq \varnothing \}
\end{equation}
respectively, where $\Omega^c=\Xi\setminus\Omega$ is the complement of $\Omega$ in $\Xi$.

\subsection{Discrete projection}

For the fully discrete method, we need to define sampling points. Given a basis~$\Psi_{\mathbf N}$, we sample in a regular (cartesian) grid oversampled by an integer $q_i>1$ in each dimension:
\begin{equation*}
	\mathcal T^{\mathbf\osi}_{\mathbf N} = \left\{\left.\left(\tfrac{k_1}{\osi_1 N_1},\dots,\tfrac{k_d}{\osi_d N_d}\right)\, \right| \, k_i=0,\dots,N_i\osi_i-1,\quad i=1,\dots, d \right\}.
\end{equation*}
Since Daubechies wavelets can only be evaluated in points $k/2^j$, $k\in\Z, j\in\N$, we choose each component of $\mathbf \osi$ to be dyadic, i.e., $q_i=2^j$, $j\in\N_0$, when approximating using a Daubechies wavelet basis.

Only the points in the intersection $\mathcal T^{\mathbf\osi,\Omega}_{\mathbf N}=\mathcal T^{\mathbf\osi}_{\mathbf N}\cap \Omega$ are of interest. The number of points in that set, $M=\#\mathcal T^{\mathbf\osi,\Omega}_{\mathbf N}$, should be larger than $N=\Pi_{j=1}^dN_j$ in order to obtain the oversampled system $A^{\mathbf\osi}_{\mathbf N}x=b^{\mathbf\osi}_{\mathbf N}$. That system is given by
\begin{align}
A^{\mathbf\osi}_{\mathbf N}(m,\mathbf l) &= \psi_{\mathbf l\mathbf N}(\mathbf t_m), \qquad \mathbf l\in I_\mathbf N, m\in I_M\label{eq:discsystem}\\
b^{\mathbf\osi}_{\mathbf N}(m) &= f(\mathbf t_m)\nonumber,
\end{align}
where $\mathbf t_m\in \mathcal T^{\mathbf\osi,\Omega}_{\mathbf N}$.

The discrete nature of the support in the discrete setting is mirrored in the definition
\begin{equation*}
\mysupp_{\mathbf\osi} \phi_{\mathbf l\mathbf N} = \mysupp\phi_{\mathbf l\mathbf N} \cap \mathcal T^{\mathbf\osi}_{\mathbf N}.
\end{equation*}
This notion of discrete support is used to determine the number of basis elements that overlap with the boundary of $\Omega$ as
\begin{equation}\label{eq:splineboundaryset_discrete}
\K_{\mathbf N}^{\mathbf\osi}(\Omega) = \{ \mathbf k \in I_{\mathbf N}, \, | \, \mysupp_{\mathbf\osi} \phi_{\mathbf k\mathbf N} \cap \Omega \neq \varnothing\text{ and }\mysupp_{\mathbf\osi} \phi_{\mathbf k\mathbf N} \cap {\Omega^c} \neq \varnothing \}
\end{equation}
and
\begin{equation}
\label{eq:waveletboundaryset_discrete}
\L^{\mathbf\osi}_{\mathbf N}(\Omega) =\{ \mathbf k \in I_{\mathbf N}, \, | \, \mysupp_{\mathbf\osi} \psi_{\mathbf k\mathbf N} \cap \Omega \neq \varnothing\text{ and }\mysupp_{\mathbf\osi} \psi_{\mathbf k\mathbf N} \cap {\Omega^c} \neq \varnothing \}.
\end{equation}

\begin{assumption}\label{ass:boundary}
	We assume that the dimension of the boundary of $\Omega\subset\Xi$ is exactly one less than the dimension of $\Omega$ itself. This means that we will not consider fractal domains. In other words, the sets $\K_{\mathbf N}(\Omega)$ and $\K^{\mathbf\osi}_{\mathbf N}(\Omega)$ grow in size as $\mathcal O\left(N^{(d-1)/d}\right)$ since we further assume that the oversampling is linear, i.e., $M=\gamma N$, with $\gamma>1$.
\end{assumption}
\begin{lemma}
	Firstly, we have
	\begin{equation}
		\#\L_{\mathbf N}(\Omega)=\mathcal O(J\#\K_{\mathbf N}(\Omega))\quad\mbox{and}\quad \#\L^{\mathbf\osi}_{\mathbf N}(\Omega)=\mathcal O(J\#\K^{\mathbf\osi}_{\mathbf N}(\Omega))
	\end{equation}
	 where $J=\Pi_{j=1}^dJ_i=\log_2(N)$. 
	 Secondly, provided Assumption \ref{ass:boundary} is satisfied
	 \begin{equation}
	 	\#\L_{\mathbf N}(\Omega) = \#\L^{\mathbf\osi}_{\mathbf N}(\Omega)=\mathcal O\left(N^{(d-1)/d}\log(N)\right). 
	 \end{equation}
\end{lemma}
\begin{proof}
The former is a direct consequence of the compact nature of the scaling basis, while the latter is a combination of the former and Lemma~\ref{lem:dwtstructure}.
\end{proof}

\section{The AZ algorithm}\label{wavs:az}

In general, the systems introduced in~\eqref{eq:contsystem} and~\eqref{eq:discsystem} are severely ill-conditioned. This is the result of the inherent redundancy of extension frame approximations. One illuminating interpretation of the redundancy is that an approximation can take any form outside of $\Omega$ while not influencing the behavior on $\Omega$. Extension frames, their ill-conditioning and further implications are studied in detail in~\cite{Adcock2019,adcock2017frames}. There, it is advised to solve the ill-conditioned systems using regularization and oversampling to obtain an numerically stable and accurate approximation. One can, e.g., use a truncated singular value decomposition (SVD) as a solver. Unfortunately, this solver and other regularized solvers generally have cubic complexity in $N$.

\begin{algorithm}[h]
	\caption{The AZ algorithm~\cite{az}}\label{alg:AZ}
	{\bf Input:} $A,Z \in \C^{M\times N}$, $b\in\C^M$ \\
	{\bf Output:} $x\in\C^N$ such that $Ax \approx b$
	\begin{algorithmic}[1]
		\State Solve $(I-AZ^*)Ax_1 = (I-AZ^*)b$ using a randomized low-rank solver
		\State $x_2 \gets Z^*(b-Ax_1)$
		\State $x \gets x_1 + x_2$
	\end{algorithmic}
\end{algorithm}

The AZ algorithm (Algorithm~\ref{alg:AZ}) consists of three simple steps and was introduced in~\cite{az} to reduce the computational complexity. It is a generalization of the algorithms proposed in~\cite{matthysen2016fast,matthysen2018function} for the more specific Fourier extension problem. There, similar ill-conditioned systems have to be solved since Fourier extension frames are similarly redundant as wavelet extension frames. The cost of Fourier extension was reduced from cubic to~$\mathcal O(N\log^2(N))$ in 1-D \cite{matthysen2016fast} and to~$\mathcal O(N^2\log^2(N))$ in 2-D \cite{matthysen2018function}.

The AZ algorithm was also successfully applied in~\cite{coppe2019splines} to reduce the cost of spline extension approximations to~$\mathcal O(N)$ in 1-D, $\mathcal O\left(N^{3/2}\right)$ in 2-D and~$\mathcal O\left(N^{3(d-1)/d}\right)$ in~$d$-D with~$d>1$. It is this latter application of AZ, and the corresponding analysis in \cite{coppe2019splines}, that we set out to extend to the wavelet case.

The AZ algorithm is shown in pseudocode in Algorithm~\ref{alg:AZ}. It solves~$Ax=b$ with a time complexity
\begin{equation}\label{eq:aztimings}
	\mathcal O(r\texttt T_\text{mult}+r^2M),
\end{equation}
where $r$ is the rank of the system in the first step and $T_\text{mult}$ is the time complexity of applying $A$ and $Z^*$ to a vector \cite{az}. The residual of the solution corresponds the approximation error in our setting, and it is equal to the residual of the solution in step one of the algorithm. The matrix $Z$ can in principle be chosen arbitrarily. However, the goal is to choose $Z$ such that the rank $r$ of the system in step 1 is small. Loosely speaking, this corresponds to choosing $Z^*$ as a pseudo-inverse to a large subspace of the range of $A$. This is where the discrete biorthogonality properties of the discrete duals constructed in \S\ref{wavs:periodicwavelet} play a decisive role.

\subsection{The choice of $Z$}\label{s:choosingz}

First, we will simply state our choice of $Z$ for both the continuous and discrete setting, based on analogy to previously studied cases. Later on, we will prove why these choices indeed give rise to a low-rank matrix~$A-AZ^*A$ in step 1 of the AZ algorithm. Briefly, we intend to make sure that $Z^*A$ approximates the identity matrix up to a small perturbation.

In the continuous setting, it suffices to choose $Z_{\mathbf N}$ equal to the identity. Indeed, recall from \eqref{eq:contsystem} that $A_{\mathbf N}$ itself is close to the identity matrix in this setting. The perturbation is related to basis functions that overlap with the boundary.

In order to reuse the results of \cite{coppe2019splines}, we introduce the scaling system matrices $\hat A$ and $\hat Z$, before we consider the wavelet system matrices $A$ and $Z$ of \eqref{eq:contsystem}-\eqref{eq:discsystem} above.  In \cite{coppe2019splines}, the B-spline system matrices do not contain inner products with the dual basis as in \eqref{eq:contsystem}, but with the primal scaling basis itself. Thus, $\hat A_{\mathbf N}$ is defined as the Gram matrix of $\Phi_N$,
\[
\hat A_{\mathbf N}(\mathbf k,\mathbf l)= \langle  \phi_{\mathbf k\mathbf N},\phi_{\mathbf l\mathbf N}\rangle_{\Omega}, \quad \mathbf k,\mathbf l\in I_{\mathbf N}.
\]
This change is also reflected in the structure of $\hat Z_{\mathbf N}$. It is not the identity matrix, but rather the Gram matrix of $\tilde\Phi_{\mathbf N}$:
\begin{equation*}
	\hat Z_{\mathbf N}(\mathbf k,\mathbf l) =  \left\langle \tilde\phi_{\mathbf k\mathbf N}, \tilde\phi_{\mathbf l\mathbf N}\right\rangle_{L^2(0,1)^d}, \qquad \mathbf k,\mathbf l\in I_{\mathbf N}.
\end{equation*}
If we denote by $W_{\mathbf J}$ the Kronecker product of the $d$ DWT matrices $\{W_{J_i}\}_{i=1}^d$, one can verify that the relation between the wavelet system matrix $A_{\mathbf N}$ and the scaling system matrix $\hat A_{\mathbf N}$ is given by
\begin{equation}\label{eq:AZtoscaling}
	A_{\mathbf N}= W_{\mathbf J}\hat Z_{\mathbf N}^*\hat A_{\mathbf N}W_{\mathbf J}^{-1},
\end{equation}
using~\eqref{eq:dualwaveletbasis} and~\eqref{eq:dwtrelations}. Therefore, $Z_{\mathbf N}^*A_{\mathbf N}= W_{\mathbf J}\hat Z_{\mathbf N}^*\hat A_{\mathbf N}W_{\mathbf J}^{-1}$ as well.

In the discrete setting, we can define both $A^{\mathbf\osi}_{\mathbf N}$ and our chosen matrix $Z^{\mathbf\osi}_{\mathbf N}$ in terms of the pointwise evaluations of the discrete dual scaling functions, followed by the continuous dual wavelet transform. Let
\[
 \hat Z^{\mathbf\osi}_{\mathbf N}(m,\mathbf l)=\tilde\phi^{\mathbf\osi}_{\mathbf l\mathbf N}(\mathbf t_m) \quad \mbox{and} \quad \hat A^{\mathbf\osi}_{\mathbf N}(m,\mathbf l)=\phi^{\mathbf\osi}_{\mathbf l\mathbf N}(\mathbf t_m), \qquad m\in I_M, \mathbf l\in I_{\mathbf N}, \mathbf t_m\in \mathcal T^{\mathbf\osi,\Omega}_{\mathbf N}.
\]
Then
\begin{equation}\label{eq:AandZtoscaling}
	Z^{\mathbf\osi}_{\mathbf N} = \hat Z^{\mathbf\osi}_{\mathbf N} W^{-1}_{\mathbf J}\quad\mbox{and}\quad A^{\mathbf\osi}_{\mathbf N} = \hat A^{\mathbf\osi}_{\mathbf N} W^{-1}_{\mathbf J}
\end{equation}
such that $\left(Z^{\mathbf\osi}_{\mathbf N}\right)^*A^{\mathbf\osi}_{\mathbf N}=W_{\mathbf J} \left(\hat Z^{\mathbf\osi}_{\mathbf N}\right)^*\hat A^{\mathbf\osi}_{\mathbf N}W_{\mathbf J}^{-1}$ again using~\eqref{eq:dwtrelations}.

In the wavelet case $T_\text{mult}$ of \eqref{eq:aztimings} is $\mathcal O(N)$, since all matrices are combinations of matrices that can be applied in $\mathcal O(N)$ operations. This is because they either contain $\mathcal O(N)$ non-zero elements or they can be represented by a (i)DWT.

\subsection{The rank and sparsity structure of $A-AZ^*A$}

In order to study the properties of the matrix $A-AZ^*A$ in the wavelet case, we again make use of the results of \cite{coppe2019splines} for the scaling bases. In the current notation, the approximation in the scaling basis leads to the matrix $\hat A-\hat A\hat Z^*\hat A$. We can write $A-AZ^*A$ as a product of this matrix with DWTs:
\begin{align}
	A_{\mathbf N}-A_{\mathbf N}Z_{\mathbf N}^*A_{\mathbf N}
	&= W_{\mathbf J}\hat A_{\mathbf N}W^{-1}_{\mathbf J}-W_{\mathbf J}\hat A_{\mathbf N}W^{-1}_{\mathbf J} W_{\mathbf J}\hat Z_{\mathbf N}^*\hat A_{\mathbf N}W_{\mathbf J}^{-1} \nonumber\\
	&= 	W_{\mathbf J}\left(\hat A_{\mathbf N}-\hat A_{\mathbf N}\hat Z_{\mathbf N}^*\hat A_{\mathbf N}\right)W^{-1}_{\mathbf J}\label{eq:AAZAcontinuousderivation}.
\end{align}
Here, we used~\eqref{eq:dwtrelations} and~\eqref{eq:AZtoscaling}. Similarly, for the discrete projection,
\begin{align}
	A^{\mathbf\osi}_{\mathbf N}-A^{\mathbf\osi}_{\mathbf N}(Z^{\mathbf\osi}_{\mathbf N})^*A^{\mathbf\osi}_{\mathbf N}
	&=\hat A^{\mathbf\osi}_{\mathbf N}W^{-1}_{\mathbf J}-\hat A^{\mathbf\osi}_{\mathbf N}W^{-1}_{\mathbf J}(\hat Z^{\mathbf\osi}_{\mathbf N}\tilde W^{-1}_{\mathbf J})^*\hat A^{\mathbf\osi}_{\mathbf N}W^{-1}_{\mathbf J}\nonumber\\
	&= 	\left(\hat A^{\mathbf\osi}_{\mathbf N}-\hat A^{\mathbf\osi}_{\mathbf N}W^{-1}_{\mathbf J}(\tilde W^{-1}_{\mathbf J})^*(\hat Z^{\mathbf\osi}_{\mathbf N})^*\hat A^{\mathbf\osi}_{\mathbf N}\right)W^{-1}_{\mathbf J}\nonumber\\
	&=\left(\hat A^{\mathbf\osi}_{\mathbf N}-\hat A^{\mathbf\osi}_{\mathbf N}(\hat Z^{\mathbf\osi}_{\mathbf N})^*\hat A^{\mathbf\osi}_{\mathbf N}\right)W^{-1}_{\mathbf J}\label{eq:AAZAdiscretederivation}
\end{align}
using~\eqref{eq:dwtrelations} and~\eqref{eq:AandZtoscaling}.

We restate the results of Theorems 6.1 and 6.2 and Corollary 6.6 in \cite{coppe2019splines} using the notation above. The proofs are algebraically tedious, yet conceptually straightforward. They rely on two central observations: (i) the matrices $A$ and $Z$ are highly sparse due to the compact support of the basis functions and (ii) the effect of $\Omega \subset \Xi$ compared to the case $\Omega = \Xi$ is confined to those basis functions that overlap with the boundary. The main technical difficulty is to accurately describe these basis functions and their corresponding index sets.

\begin{lemma}[{\cite[Theorem 6.1, Theorem 6.2, Corrollary~6.6]{coppe2019splines}}]\label{lem:scalinglowrank}
	For the AZ pairs $(\hat A_{\mathbf N},\hat Z_{\mathbf N})$ and $(\hat A^{\mathbf\osi}_{\mathbf N},\hat Z^{\mathbf\osi}_{\mathbf N})$, the matrix $A-AZ^*A$ has
	\begin{enumerate}
		\item at most $\#\K_{\mathbf N}(\Omega)$ and  $\#\K^{\mathbf\osi}_{\mathbf N}(\Omega)$ non-zero columns,
		\item at most rank $\#\K_{\mathbf N}(\Omega)$ and  $\#\K^{\mathbf\osi}_{\mathbf N}(\Omega)$,
		\item $\mathcal O(\#\K_{\mathbf N}(\Omega))$ and  $\mathcal O(\#\K^{\mathbf\osi}_{\mathbf N}(\Omega))$ non-zero rows,
		\item $\mathcal O(\#\K_{\mathbf N}(\Omega))$ and  $\mathcal O(\#\K^{\mathbf\osi}_{\mathbf N}(\Omega))$ non-zero elements,
	\end{enumerate}
	respectively. The constants in the big $\mathcal O$ notation are independent of $N$.
\end{lemma} 

We will also add a more precise statement on the number of non-zero rows. For this we need to introduce two more index sets:
\begin{eqnarray}
\nonumber\M_{\mathbf N}(\Omega) &=& \{ \mathbf k \in I_{\mathbf N}\, |\, \forall\mathbf l\in\K_{\mathbf N}(\Omega), \forall\mathbf i\in I_{\mathbf N}:\\
&&\qquad\mysupp \tilde \phi_{\mathbf i\mathbf N}\cap \mysupp \phi_{\mathbf l\mathbf N}\neq\emptyset\text{ and }\mysupp \phi_{\mathbf i\mathbf N}\cap\mysupp\phi_{\mathbf k\mathbf N}\neq\emptyset   \}\label{eq:Mcont}
\end{eqnarray}
and 
\begin{eqnarray}
\M^{\mathbf\osi}_{\mathbf N}(\Omega) = \{  m \in I_{M}\, |\, \forall\mathbf l\in \K^{\mathbf\osi}_{\mathbf N}, \forall\mathbf i\in I_{\mathbf N}:
 \mysupp_{\mathbf\osi} \tilde \phi_{\mathbf i\mathbf N}\cap \mysupp_{\mathbf\osi} \phi_{\mathbf l\mathbf N}\neq\emptyset\text{ and }\phi_{\mathbf i\mathbf N}(\mathbf t_m)\neq 0   \}.\label{eq:Mdiscr}
\end{eqnarray}
The set $\M_{\mathbf N}(\Omega)$ corresponds to the indices of all dual basis functions that overlap with any primal basis function that overlaps with the boundary. Similarly, the set $\M^{\mathbf\osi}_{\mathbf N}(\Omega)$ has indices of all points in the support of any dual discrete basis function which overlaps with any primal basis function that overlaps with the boundary.

\begin{theorem}\label{thm:lowrow}
	For the AZ pairs $(\hat A_{\mathbf N},\hat Z_{\mathbf N})$ and $(\hat A^{\mathbf\osi}_{\mathbf N},\hat Z^{\mathbf\osi}_{\mathbf N})$, the matrix $A-AZ^*A$ has non-zero row indices $\M_{\mathbf N}(\Omega)$ and  $\M^{\mathbf\osi}_{\mathbf N}(\Omega)$
	respectively. Furthermore, $\#\M_{\mathbf N}(\Omega)=\mathcal O(\#\K_{\mathbf N}(\Omega))$ and  $\#\M^{\mathbf\osi}_{\mathbf N}(\Omega)=\mathcal O(\#\K^{\mathbf\osi}_{\mathbf N}(\Omega))$.	
\end{theorem}
\begin{proof}
For the AZ pair $(\hat A_{\mathbf N},\hat Z_{\mathbf N})$ we proceed similarly as in the proof of \cite[Theorem 6.1]{coppe2019splines}. First we note that:
\begin{align}\label{eq:proofnonzerorows1}
\left(I - Z^*{A}\right)(\mathbf k,\mathbf l) = \left\{
\begin{array}{ll}
\delta_{\mathbf k\mathbf l} & \text{if }\mysupp\phi_{\mathbf N,\mathbf l}\subset\Omega^c,\\
0 & \text{if }\mysupp\phi_{\mathbf l\mathbf N}\subset\Omega,\\
\delta_{\mathbf k\mathbf l}-(\phi_{\mathbf l\mathbf N},\tilde\phi_{\mathbf k\mathbf N})_{L^2(\Omega)} &\text{otherwise}.
\end{array}\right.
\end{align}
This is due to the compact support of the basis functions and to the continuous biorthogonality. In case $\Omega=\Xi$ all matrix entries would be zero, here they differ if $\phi_{\mathbf{N},\mathbf l}$ is supported outside of $\Omega$ or if it overlaps with the boundary.

Secondly, multiplication by $A$ on the left yields: 
\begin{align}\label{eq:plungestructure2}
\left( A (I - Z^*A) \right)(\mathbf k,\mathbf l) = \left\{
\begin{array}{ll}
0 & \text{if }\mysupp\phi_{\mathbf l\mathbf N}\subset\Omega^c,\\
0 & \text{if }\mysupp\phi_{\mathbf l\mathbf N}\subset\Omega,\\
a_{\mathbf k\mathbf l} &\text{otherwise},
\end{array}\right.
\end{align}
where the values $a_{\mathbf k\mathbf l}$ may or may not be zero.

We note that $A-AZ^*A$ has several columns that are identically zero. The proof of \cite[Theorem 6.1]{coppe2019splines} also shows that the non-zero column indices of $A-AZ^*A$ are given by $\K_{\mathbf N}(\Omega)$, hereafter abbreviated by $\K$.

If we let~$E\in\{0,1\}^{N\times \#\K}$ be the extension matrix that extends $\K$ to $I_{\mathbf N}$,
\begin{equation*}
E(\mathbf k,\mathbf l) = \delta_{\mathbf k\mathbf l}, \qquad \mathbf k\in I_{\mathbf N}, \mathbf l\in \K,
\end{equation*}
then the matrix~${A}(I-{Z}^*{ A})E$ contains all non-zero columns of ${ A}(I-{ Z}^*{ A})$ but has size~$N\times\#\K$ instead of~$N\times N$.

The matrix $(I- Z^* A)E$ contains a small number of non-zero rows since
\begin{align*}
	[(I-Z^*A)E](\mathbf k,\mathbf l) = \left\{
	\begin{array}{ll}
	\delta_{\mathbf k\mathbf l}-(\phi_{\mathbf l \mathbf N},\tilde\phi_{\mathbf k \mathbf N})_{L^2(\Omega)} & \text{if }\mathbf l\in\K\\
	0 & \text{otherwise}
	\end{array}\right.
\end{align*}
using \eqref{eq:splineboundaryset} and~\eqref{eq:proofnonzerorows1}.
The non-zero row indices of this matrix are in the index set 
\begin{equation}\label{eq:I1}
 I_1 = \K \cup \{ \mathbf k \in I_{\mathbf N}\, |\, \forall \mathbf l\in\K: \mysupp \tilde \phi_{\mathbf k \mathbf N}\cap \mysupp \phi_{\mathbf l \mathbf N}\neq \emptyset \}
\end{equation}
so we can rewrite ${ A}(I-{ Z}^*{ A})E$ as 
\[  { A}E_1E_1^*(I-{ Z}^*{ A})E \]
with $E_1\in \{0,1\}^{N\times \# I_{\mathbf N}}$ an extension matrix derived from the index set $I_1$ 
\begin{equation*}
E_1(\mathbf k,\mathbf l) = \delta_{\mathbf k\mathbf l}, \qquad \mathbf k\in I_{\mathbf N}, \mathbf l\in I_1. 
\end{equation*}

The non-zero row index set of ${ A}(I-{ Z}^*{ A})$ is thus 
\begin{equation*}
\{\mathbf k \in I_{\mathbf N}\, |\, 
\forall\mathbf i\in  I_1: \mysupp \phi_{\mathbf i \mathbf N}\cap\mysupp\phi_{\mathbf k \mathbf N}\neq\emptyset   \}
\end{equation*}
which is equivalent to $\M_{\mathbf N}(\Omega)$ in \eqref{eq:Mcont} after substituting the earlier found expression for $ I_1$.

For the AZ pair  $(\hat A^{\mathbf\osi}_{\mathbf N},\hat Z^{\mathbf\osi}_{\mathbf N})$, $I_1$ retains the form of~\eqref{eq:I1} but the meaning of $\K$ is now the one in~\eqref{eq:splineboundaryset_discrete}. That is why the non-zero row index set of ${ A}(I-{ Z}^*{A})$ is~\eqref{eq:Mdiscr}.

Because of Lemma~\ref{lem:scalinglowrank}, $\#\M_{\mathbf N}(\Omega)=\mathcal O(\#\K_{\mathbf N}(\Omega))$ and  $\#\M^{\mathbf\osi}_{\mathbf N}(\Omega)=\mathcal O(\#\K^{\mathbf\osi}_{\mathbf N}(\Omega))$. 
\end{proof}

Lemmas \ref{lem:dwtstructure} and \ref{lem:scalinglowrank}, Theorem~\ref{thm:lowrow} and the relations between the matrices $A$ and $Z$ for the scaling and wavelet bases ~\eqref{eq:AAZAcontinuousderivation}-\eqref{eq:AAZAdiscretederivation} combined lead to the following theorem.
\begin{theorem}\label{thm:waveletlowrank}
	For the AZ pairs $(A_{\mathbf N},Z_{\mathbf N})$ and $(A^{\mathbf\osi}_{\mathbf N}, Z^{\mathbf\osi}_{\mathbf N})$, the matrix $A-AZ^*A$ has
	\begin{enumerate}
		\item at most $\#\L_{\mathbf N}(\Omega)$ and  $\#\L^{\mathbf\osi}_{\mathbf N}(\Omega)$ non-zero columns,
		\item at most rank $\#\K_{\mathbf N}(\Omega)$ and  $\#\K^{\mathbf\osi}_{\mathbf N}(\Omega)$,
		\item  $\mathcal O(\#\L^{\mathbf\osi}_{\mathbf N}(\Omega))$ and at most $\#\M_{\mathbf N}(\Omega)$  non-zero rows,
		\item $\mathcal O(J\#\L_{\mathbf N}(\Omega))$ and  $\mathcal O(J\#\K^{\mathbf\osi}_{\mathbf N}(\Omega))$ non-zero elements,
	\end{enumerate}
	respectively. The constants in the big $\mathcal O$ notation are independent of $N$. The index sets are defined by \eqref{eq:splineboundaryset}-\eqref{eq:waveletboundaryset_discrete}.
\end{theorem}
\begin{proof}
 	\begin{enumerate}
		\item One can follow a similar reasoning as in the proofs of Theorem~6.1, Theorem~6.2 in \cite{coppe2019splines} to obtain this statement. There it is shown that the non-zero column indices of $\hat A-\hat A\hat Z^*\hat A$  are those of the \emph{scaling} basis elements that overlap with the boundary, i.e., $\K_{\mathbf N}(\Omega)$ and  $K^{\mathbf\osi}_{\mathbf N}(\Omega)$. This reasoning applied to $A-AZ^*A$ says that the non-zero column indices are those of the \emph{wavelet} basis elements that overlap with the boundary, i.e., $\L_{\mathbf N}(\Omega)$ and  $\L^{\mathbf\osi}_{\mathbf N}(\Omega)$.

		%Another way to obtain the same result is by noting that $\hat A-\hat A\hat Z^*\hat A$ --- a matrix with $\#\K_{\mathbf N}(\Omega)$ ($\#\K^{\mathbf\osi}_{\mathbf N}(\Omega)$) non-zero columns --- is right multiplied with $W^{-1}_{\mathbf J}$ --- a matrix that has $\mathcal O(J)$ non-zero-elements per row.
		\item The rank of $A-AZ^*A$ is less than or equal to the rank of $\hat A-\hat A\hat Z^*\hat A$ since $A-AZ^*A=W(\hat A-\hat A\hat Z^*\hat A)W^{-1}$ or $A-AZ^*A=(\hat A-\hat A\hat Z^*\hat A)W^{-1}$ and $W$ is of full rank.
		\item In the discrete case $\hat A-\hat A\hat Z^*\hat A$ has at most $\#\M^{\mathbf\osi}_{\mathbf N}(\Omega)$ non-zero rows, so $A-AZ^*A$ which is the former right-multiplied with a discrete wavelet transform has at most  $\#\M^{\mathbf\osi}_{\mathbf N}(\Omega)$ non-zero rows. In the continuous case, $A-AZ^*A=W(\hat A-\hat A\hat Z^*\hat A)W^{-1}$ where $\hat A-\hat A\hat Z^*\hat A$ has $\mathcal O(\#\K_{\mathbf N}(\Omega))$ non-zero rows and $W$ has $\mathcal O(J)$ non-zero elements per column, thus $A-AZ^*A$ has $\mathcal O(J \#\K_{\mathbf N}(\Omega))=\mathcal O(\#\L_{\mathbf N}(\Omega))$ non-zero rows.
		\item  The compact support of the scaling bases ensures that the number of non-zero elements per column and per row of $\hat A-\hat A\hat Z^*\hat A$ is bounded by a constant independent of $N$. Hence  $\hat A-\hat A\hat Z^*\hat A$ has  $\mathcal O(\#\K_{\mathbf N}(\Omega))$ ($\mathcal O(\#\K^{\mathbf\osi}_{\mathbf N}(\Omega))$) non-zero elements (as shown in Lemma~\ref{lem:scalinglowrank}). The bound on the number of non-zero elements per row and right-multiplication of $\hat A-\hat A\hat Z^*\hat A$ with $W^{-1}$ results in a matrix with $\mathcal O(J)$ elements per row. Therefore, the matrix has $\mathcal O(J\#\K^{\mathbf\osi}_{\mathbf N})$ non-zero elements. 
		
		The same holds for the continuous case. Furthermore, since the matrix still has a bounded number of non-zero elements per column after right-multiplication with $W^{-1}$, left-multiplication with $W$ results in a matrix where the number of non-zero elements per row and per column grow like $\mathcal O(J)$ and $A-AZ^*A$ contains in the continuous case $\mathcal O(J\#\L^{\mathbf\osi}_{\mathbf N})$ non-zero elements.
	\end{enumerate}
\end{proof}

\subsection{The vanilla AZ algorithm}\label{ss:vanillaAZ}

\begin{figure}[t]
	\centering
	\resizebox{\textwidth}{!}{\input{img/AZtimings1d-3d.tikz}}
	\caption{\label{fig:AZtimings1d-3d}Timings in seconds of the AZ algorithm (Algorithm~\ref{alg:AZ}) applied to the approximation of (left) $f(x)=e^x$ on $[0,1/2]$ using a 1-D wavelet basis on $[0,1]$, (middle) $f(x,y)=e^{xy}$ on the disk with center $[1/2,1/2]$ and radius $0.35$ using a 2-D wavelet basis on $[0,1]^2$, (right) $f(x,y,z)=e^{xyz}$ on ball with center $[1/2,1/2,1/2]$ and radius $0.4$ using a 3-D wavelet basis on $[0,1]^3$. We approximate using several wavelets of primal orders $2$ to $4$, shown with different markers. The expected asymptotic results of Theorem~\ref{thm:aztimings} are shown by the black dashed line: $\mathcal O(N)$ in 1-D, $\mathcal O(N^2)$ in 2-D and $\mathcal O(N^{7/3})$ in 3-D.}
\end{figure}
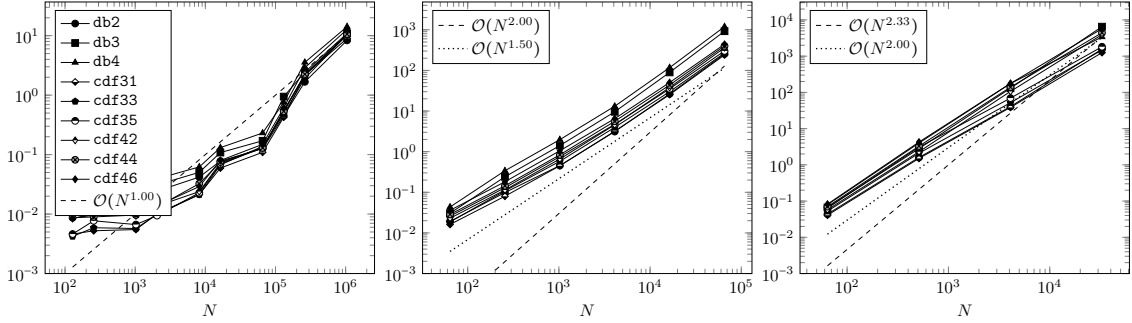

We will refer to Algorithm~\ref{alg:AZ} as the vanilla AZ algorithm. Here, the matrices $A$ and $Z$, as well as the matrix $A-AZ^*A$ have dimension $M \times N$, where $M$ is the total number of sample points and $N$ is the total number of degrees of freedom. We have already established that this algorithm is not optimal, because the matrix $A-AZ^*A$ has a large number of zero-rows and zero-columns. Still, because of its apparent simplicity, we state the expected computational complexity. Also, a surprising feature of the low-rank solver we have used in our implementation in step 1 is that its computational complexity is actually much better.

\begin{theorem}\label{thm:aztimings}
	Provided Assumption \ref{ass:boundary} is satisfied, the AZ algorithm (Algorithm~\ref{alg:AZ}) using the AZ pairs $(A_{\mathbf N},Z_{\mathbf N})$ and $(A^{\mathbf\osi}_{\mathbf N},Z^{\mathbf\osi}_{\mathbf N})$ can be implemented with $\mathcal O(N)$ operations in 1-D, $\mathcal O(N^{2})$ operations in 2-D and $\mathcal O(N^{(3d-2)/d}) $ operations in $d$-D, $d>1$.
\end{theorem}
\begin{proof}
	Recall that the computational of the AZ algorithm, Algorithm~\ref{alg:AZ}, is
	\[ 	\mathcal O(r\texttt T_\text{mult}+r^2M). \]
	Assumption~\ref{ass:boundary} combined with Theorem~\ref{thm:waveletlowrank} shows that $r=\mathcal O(N^{(d-1)/d})$. In \S\ref{s:choosingz} it was shown that $\texttt T_\text{mult}$ is $\mathcal O(N)$. The full AZ-algorithm therefore requires $\mathcal O(N)$ operations in 1-D, and $\mathcal O(N^{(3d-2)/d}) $ operations in $d$-D, $d>1$.
\end{proof}

Theorem~\ref{thm:aztimings} is illustrated in Figure~\ref{fig:AZtimings1d-3d}. In these numerical results we only consider the discrete setting, because the inner product integrals in~\eqref{eq:contsystem} are not trivial to compute. The complete code for this experiment and the following experiments is available online \cite{FrameFunTranslates.jl,FrameFunWavelets.jl}.

As a low-rank solver in step 1 of the algorithm we have used the low-rank QR solver implemented in \cite{LowRankApprox.jl}. It is an algorithm that builds the QR factorization of a matrix $A$ by random sampling. This means that  $A$ is sampled by multiplying $A$ with $r$ random vectors. It is an adaptive algorithm that determines $r$ adaptively in order to obtain an accurate factorization \cite{halko2011finding}. The low-rank QR is in general a bit more efficient than the low-rank SVD implemented in the same package. We note in Figure~\ref{fig:AZtimings1d-3d}, to our surprise, that the AZ algorithm outperforms the expected complexity of Theorem~\ref{thm:aztimings}. The experimental time complexity follows the dotted line more closely than the dashed one in 2-D and 3-D. We did not explicitly indicate our knowledge of the low number of non-zero rows (see Theorem~\ref{thm:waveletlowrank}) to the low-rank solver. However, the solver apparently took notice automatically and used it to its advantage, lowering the time complexity to that of Theorem~\ref{thm:reducedaztimings} in the following section. There, we do assume an algorithm that explicitly uses knowledge of the non-zero rows and columns.

\subsection{The reduced AZ algorithm}\label{ss:reducedAZ}

Following Theorem~\ref{thm:waveletlowrank} the non-zero rows and columns are described by the index sets $\L^{\mathbf\osi}_{\mathbf N}(\Omega)$ (or $\L^{\mathbf\osi}_{\mathbf N}(\Omega)$ in the discrete setting) and $\M^{\mathbf\osi}_{\mathbf N}(\Omega)$ ($\M^{\mathbf\osi}_{\mathbf N}(\Omega)$). To ease notation, we will use $\L$ and $\M$ instead if the meaning is independent of the context or can be deduced out of the context. 

By iterating once over all wavelet basis functions and checking their support one can determine the indices in $\L$ in $\mathcal O(N)$ operations. Indeed, owing to the completely regular structure of the bases involved, the supports of the basis functions are easily computed. From the definitions \eqref{eq:Mcont}-\eqref{eq:Mdiscr} we see that we can compute $\M$ by first computing $\K$ --- this in $\mathcal O(N)$ operations by iterating over all scaling basis elements --- then checking their overlap with the duals, again in $\mathcal O(N)$ operations.

If these index sets are known we can create the (sparsely representable) extension and restriction matrices $E=\{0,1\}^{N\times\#\L}$ and $R=\{0,1\}^{\#\M\times M}$:
\begin{eqnarray*}
	E(\mathbf k, \mathbf l) &=& \delta_{\mathbf k\mathbf l}, \qquad \mathbf k\in I_{\mathbf N}, \mathbf l\in \L \\
	R(k, l) &=& \delta_{kl}, \qquad  k\in\M, l\in I_m.
\end{eqnarray*}

\begin{algorithm}[t]
	\caption{The reduced AZ algorithm \cite[Algorithm~2]{coppe2019splines}}\label{alg:reducedAZ}
	{\bf Input:} $A,Z \in \C^{M\times N}$, $b\in\C^M$, $\epsilon > 0$ \\
	{\bf Output:} $x\in\C^N$ such that $Ax \approx b$
	\begin{algorithmic}[1]
		\State Determine $E$ which extends the index set $\L$ to $I_{\mathbf N}$
		\State Determine $R$ which restricts $I_{M}$ to the indices of non-zero rows of $(I-AZ^*)AE$.
		\State Solve $R(I-AZ^*)AEx_1 = R(I-AZ^*)b$
		\State $x_2 \gets Z^*(b-AEx_1)$
		\State $x \gets Ex_1 + x_2$
	\end{algorithmic}
\end{algorithm}
With these matrices we construct the matrix $R(I-AZ^*)AE$ that holds the same information as $A-AZ^*A$ but without all known zero rows and columns. This matrix is used in the reduced AZ algorithm~(Algorithm \ref{alg:reducedAZ}) \cite{coppe2019splines}. The number of non-zero columns of $R(I-AZ^*)AE$ still grows at a faster rate than its rank, see Theorem~\ref{thm:waveletlowrank}. That is why in the following theorem we distinguish between a full direct solver, such as a pivoted QR, and a low-rank direct solver, such as a randomized low-rank QR. Note also that reducing the size does not change the complexity of applying the matrices. Both $A-AZ^*A$ and $R(I-AZ^*)AE$ result in a matrix-vector multiply that takes $\mathcal O(N)$ operations.

\begin{theorem}\label{thm:reducedaztimings}
	Provided Assumption \ref{ass:boundary} is satisfied, the reduced AZ algorithm (Algorithm~\ref{alg:reducedAZ}) using the AZ pair $(A_{\mathbf N},Z_{\mathbf N})$ can be implemented with 
	\begin{enumerate}
		\item $\mathcal O(N\log(N))$ operations in 1-D, $\mathcal O(\log(N)^3N^{3/2})$ operations in 2-D and $\mathcal O(\log(N)^3N^{3(d-1)/d})$ operations in $d$-D with $d>1$ if a \emph{full} direct solver is used in the 3rd step.
		\item $\mathcal O(N)$ operations in 1-D, $\mathcal O(\log(N)N^{3/2})$ operations in 2-D and $\mathcal O(\log(N)N^{3(d-1)/d})$ operations in $d$-D with $d>1$ if a \emph{low-rank} direct solver is used in the 3rd step.
	\end{enumerate}
	Using the AZ pair  $(A^{\mathbf\osi}_{\mathbf N},Z^{\mathbf\osi}_{\mathbf N})$ the reduced AZ algorithm can be implemented with 
	\begin{enumerate}
		\item $\mathcal O(N\log(N))$ operations in 1-D, $\mathcal O(\log(N)^2N^{3/2})$ operations in 2-D and $\mathcal O(\log(N)^2N^{3(d-1)/d})$ operations in $d$-D with $d>1$ if a \emph{full} direct solver is used in the 3rd step.
		\item $\mathcal O(N)$ operations in 1-D, $\mathcal O(N^{3/2})$ operations in 2-D and $\mathcal O(N^{3(d-1)/d})$ operations in $d$-D with $d>1$ if a \emph{low-rank} direct solver is used in the 3rd step.
	\end{enumerate}
\end{theorem}
\begin{proof}
	The proof largely follows \cite[Theorem 6.8]{coppe2019splines}, with the exception of a factor $\mathcal O(J)$ in some places and using Theorem~\ref{thm:waveletlowrank}.
	
	The extension matrix $E$ can be constructed in $\mathcal O(N)$ operations by iterating once over all wavelet basis functions. The matrix $R$ can be constructed in $\mathcal O(N)$ by creating the set $\M$. For the discrete case, this immediately gives the non-zero row indices. For the continuous case, an additional DWT is needed to find these non-zero row indices.\footnote{For completeness, in our implementation we have performed a DWT on a vector of length $N$ with a \texttt{NaN} instead of a floating point number at the indices $\K$ (non-zero indices of $\hat A_{\mathbf N}-\hat A_{\mathbf N}\hat Z_{\mathbf N}^*\hat A_{\mathbf N}$). Since a mathematical operation between a \texttt{NaN} and a floating point number results in a  \texttt{NaN}, the resulting vector thus contains \texttt{NaN}s at the non-zero indices of $W_{\mathbf J}(\hat A_{\mathbf N}-\hat A_{\mathbf N}\hat Z_{\mathbf N}^*\hat A_{\mathbf N})$.}  The matrix $R(A-A Z^*A)E$ has size $m \times n$, where \begin{enumerate}
		\item $n=\#\L = \mathcal O(JN^{(d-1)/d})$ 
		\item $m=\mathcal O(\# \L^{\mathbf\osi})=\mathcal O(N^{(d-1)/d})$ in the continuous setting and $m=\#\M =\mathcal O(N^{(d-1)/d})$ in the discrete setting.
	\end{enumerate} 
However, by Theorem~\ref{thm:waveletlowrank} and Assumption~\ref{ass:boundary} its rank $r$ is only $\mathcal O(\K) = \mathcal O(N^{(d-1)/d})$.

	Next, creating and solving the $m\times n$ linear system with a direct solver requires $\mathcal O(n\texttt T_\text{mult}+mn^2)$ operations, while it takes $\mathcal O(r\texttt T_\text{mult}+mr^2)$ with a low-rank solver. Therefore, the full time complexity of the algorithm is
	\[
	\mathcal O(N + n\texttt T_\text{mult} +  nM + mn^2)  = \mathcal O(nN + mn^2)
	\]
	with a full direct solver and 
	\[
	\mathcal O(N + r\texttt T_\text{mult} +  nM + mr^2)  = \mathcal O(rN + mr^2)
	\]
	with a low-rank direct solver, where we take linear oversampling $M=\mathcal O(N)$ into account and all matrices can be applied in $\mathcal O(N)$ operations, so  $T_\text{mult}=\mathcal O(N)$. Filling in the results obtained earlier in the proof results in the statement of the theorem.
\end{proof}

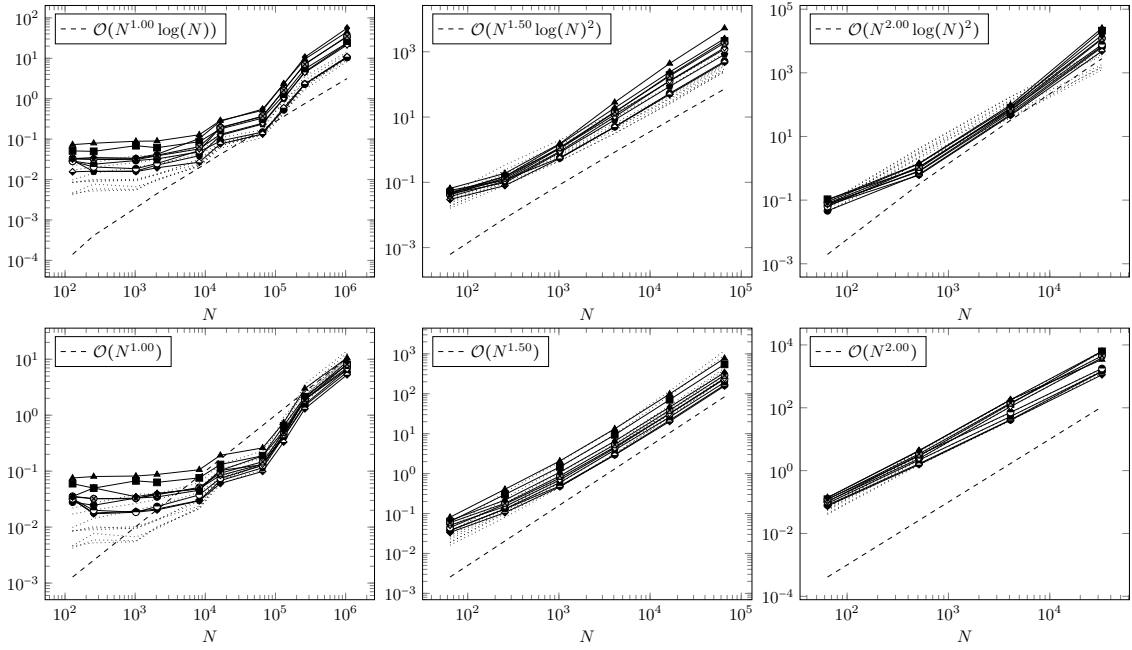
\begin{figure}
	\centering
	\resizebox{\textwidth}{!}{\input{img/AZR1timings1d-3d.tikz}}\\
	\resizebox{\textwidth}{!}{\input{img/AZR2timings1d-3d.tikz}}
	\caption{\label{fig:AZR1timings1d-3d}Timings in seconds of the reduced AZ algorithm (Algorithm~\ref{alg:reducedAZ}) applied to the approximation problems of Figure~\ref{fig:AZtimings1d-3d}. Top row: A full direct solver is used in step 3 of the algorithm. Bottom row: A low-rank solver. The expected asymptotic results of Theorem~\ref{thm:reducedaztimings} are shown by the black dashed line. The timings of Figure~\ref{fig:AZtimings1d-3d} (vanilla AZ) are repeated with black dotted lines.}
\end{figure}
\begin{figure}
	\centering
	\resizebox{\textwidth}{!}{\input{img/AZR2timings2d.tikz}}
	\caption{\label{fig:AZR2timings2d}The bottom middle panel of Figure \ref{fig:AZR1timings1d-3d} in more detail.}
\end{figure}
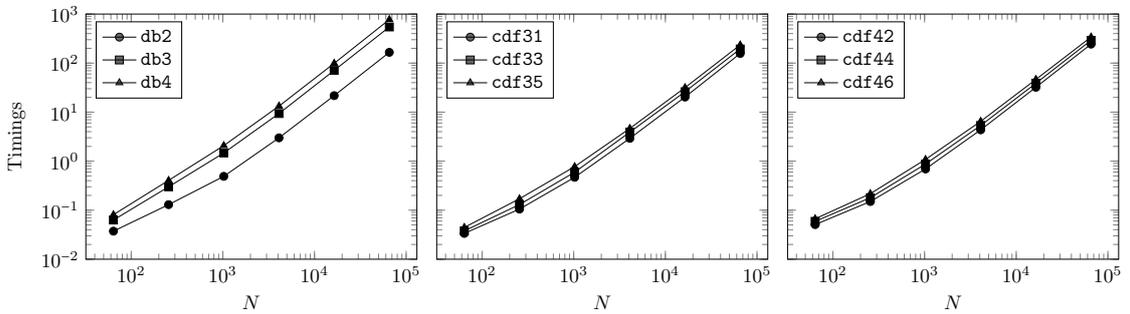

The statements in Theorem~\ref{thm:reducedaztimings} are corroborated in Figure~\ref{fig:AZR1timings1d-3d}. The difference in time complexity between a full and a low-rank direct solver is only a matter of logarithmic factors. However, the experiments indicate that the low-rank solver is the better choice, at least in our implementation, since apart from slightly improved time complexity its cost seems also substantially lower in absolute terms. As noted at the end of \S\ref{ss:vanillaAZ}, the reduced AZ algorithm does not perform much better than vanilla AZ (in the dotted lines) since the randomized solver seems to take advantage of the non-zero rows. There are however some advantages to explicitly removing the zero columns and rows. First, the sampling matrix in the low-rank solver reduces in size from $\mathcal O(N\times N^{(d-1)/d})$ to $O(N^{(d-1)/d}\log(N)\times N^{(d-1)/d})$. Secondly, all matrices and decompositions stored are smaller. Finally, we do not need to select a solver that implicitly takes care of the zero rows.

From Figure \ref{fig:AZR2timings2d} it is clear that Daubechies wavelets are less efficient than CDF wavelets (in our implementation). Furthermore, we see that wavelets with higher $\tilde p$ and hence larger support are less efficient than those with a lower number of dual vanishing moments. Of course, their wider support results in a larger number of basis functions that overlap with the boundary.

The left bottom panel of Figure \ref{fig:AZerrors1d-3d} shows the residual of the experiments in Figure \ref{fig:AZR2timings2d}, but with $\mathbf\osi=(4,4)$ instead of $\mathbf\osi=(2,2)$. The residual directly corresponds to the approximation error in the point samples. 
We see the expected algebraic convergence and that wavelets with a higher order converge faster. Note that the oversampling factors $\mathbf\osi=(4,4)$ were needed here to make a clear distinction between the experimental CDF convergence rates. With $\mathbf\osi=(2,2)$, the choice in the previous experiments including Figure~\ref{fig:AZR2timings2d},  both \texttt{cdf}3$\tilde p$ and \texttt{cdf}4$\tilde p$ appear to converge approximately at the same rate. Increasing the oversampling factor makes the difference in convergence rates more pronounced.

\subsection{The sparse AZ algorithm}\label{ss:sparseAZ}

\begin{algorithm}[t]
	\caption{The sparse AZ algorithm \cite[Algorithm~3]{coppe2019splines}}\label{alg:sparseAZ}
	{\bf Input:} $A,Z \in \C^{M\times N}$, $b\in\C^M$ \\
	{\bf Output:} $x\in\C^N$ such that $Ax \approx b$
	\begin{algorithmic}[1]
		\State Create sparse matrix $(I-AZ^*)A$
		\State Solve $(I-AZ^*)Ax_1 = (I-AZ^*)b$ using sparse QR
		\State $x_2 \gets Z^*(b-Ax_1)$
		\State $x \gets x_1 + x_2$
	\end{algorithmic}
\end{algorithm}

Finally, we consider an algorithm that only exploits the sparsity of the matrices $A$, $Z$ and $A-AZ^*A$. As for B-splines in \cite{coppe2019splines}, here we use the sparse direct rank-revealing QR decomposition of \cite{Davis2009} in the first step of the AZ algorithm instead of a low-rank solver.

The sparse AZ algorithm is formulated in Algorithm~\ref{alg:sparseAZ}. The first step is the creation of the matrix $A-AZ^*A$ in sparse form. It was shown in \cite{coppe2019splines} that a sparse version $\hat A-\hat A\hat Z^*\hat A$ (holding $\mathcal O(N^{(d-1)/d})$ non-zero values and $\mathcal O(1)$ non-zero elements in each column and row) can be created in $\mathcal O(N)$ operations. Given this sparse matrix  $\hat A-\hat A\hat Z^*\hat A$, we can easily create the matrix containing its non-zero columns $\hat A(I-\hat A\hat Z^*)\hat AE$ with 
\begin{equation*}
	E(\mathbf k, \mathbf l) = \delta_{\mathbf k\mathbf l}, \qquad \mathbf k\in I_{\mathbf N}, \mathbf l\in \K.
\end{equation*}
It remains to construct the sparse matrix $E^*W^{-1}$ since
\[
 A-AZ^*A = \hat A(I-\hat A\hat Z^*)\hat AEE^*W^{-1}.
\]
Provided Assumption~\ref{ass:boundary} is satisfied, we deduce from Lemma~\ref{lem:dwtstructure} that $E^*W^{-1}$ has $\mathcal O(N^{(d-1)/d}\log N)$ non-zero elements, i.e., $\mathcal O(N^{(d-1)/d})$ rows with $\mathcal O(\log N)$ elements each. Because of the structure in both matrices the standard sparse matrix-matrix multiply in Julia 1.3 performs the multiplication of $A(I-\hat A\hat Z^*)\hat AE$ with $EW^{-1}$ in $\mathcal O(N^{(d-1)/d}\log(N))$.
%\MAGENTA{Dit is nogal implementatiespecifiek, hier zou een referentie nuttig zijn. Ofwel naar Julia ofwel (beter) naar het "Gustavsen" algoritme?} \GREEN{Dan kies ik toch voor `standard sparse matrix-matrix multiply in Julia 1.3` dit lijkt me de veiligste optie omdat in de julia code letterlijk staat `based on Gustavsen's matrix multiplication algorithm` en niet `implements the Gustavsen's matrix multiplication algorithm`}

Hence, for the best possible time complexity in the first step of the sparse AZ algorithm, we need to construct $E^*W^{-1}$ in $\mathcal O(JN^{(d-1)/d})$ operations. Here, we present an algorithm that in the worst case is $\mathcal O(N\log N)$.  An inspection of \S\ref{ss:dwt} using compactly supported sequences and convolutions instead of matrices yields the insight that each column of $W_{J_i}^{-1}$ (1-D iDWT transform) contains a shifted version of either one of 
\begin{eqnarray}\label{eq:idwtfilters}
	&&g^{J},\\\nonumber
	&&[g^{J-1}]_{\uparrow 2^1}\star h^{J},\\\nonumber
	&&[g^{J-2}]_{\uparrow 2^2}\star [h^{J-1}]_{\uparrow 2^1}\star h^{J}, \\\nonumber
	&&\qquad\qquad\quad\vdots\\\nonumber
	&&[g^{1}]_{\uparrow 2^{J-1}}\star [h^{2}]_{\uparrow 2^{J-1}}\star\cdots\star [h^{J-1}]_{\uparrow 2^1} h^{J}, \\\nonumber
	&&[h^{1}]_{\uparrow 2^{J-1}}\star [h^{2}]_{\uparrow 2^{J-1}}\star\cdots\star [h^{J-1}]_{\uparrow 2^1} h^{J}
\end{eqnarray}
where we left out the dimension subindex to $J_i$ and where $[a]_{\uparrow q}$ denotes upsampling and $a^J$ periodisation:
\begin{equation*}
	([a]_{\uparrow q})_k = a_{kq}, \qquad (a^J)_k = \sum_{l\in\Z} a_{k+l2^J}.
\end{equation*}

A convolution $a\star b$ can be computed in $\mathcal O(n\log n)$ operations with $n$ the sum of the supports of $a$ and $b$ using the Fast Fourier Transform. The top filter in \eqref{eq:idwtfilters}  has constant support, but the support lengths grow steadily up to $J_i2^{J_i}$. Therefore, the $J_i+1$ filters can be computed in $\mathcal O(N_i\log N_i)$ operations.

If we assume that the degrees of freedom are evenly distributed over all dimensions, then we have $N_i=\mathcal O(N^{1/d})$ and the total number of nonzero entries in $E^*W^{-1}$ is $\mathcal O(JN^{(d-1)/d})$. Hence, knowing their locations, the sparse matrix can also be constructed in $\mathcal O(JN^{(d-1)/d})$ operations. In the worst case, if one dimension has all degrees of freedom, then the cost of the full construction algorithm of the sparse $E^*W^{-1}$ may be as large as $\mathcal O (N\log N)$. This is unlikely to be the case in practice.

\begin{figure}
	\centering
	\resizebox{\textwidth}{!}{\input{img/AZStimings1d-3d.tikz}}
	\caption{\label{fig:AZStimings1d-3d}Timings in seconds of the sparse AZ algorithm (Algorithm~\ref{alg:sparseAZ}) applied to the approximation problems of Figure~\ref{fig:AZtimings1d-3d}. In the black dashed line: $\mathcal O(N)$. The timings of Figure~\ref{fig:AZtimings1d-3d} are repeated in the black dotted lines.}
\end{figure}
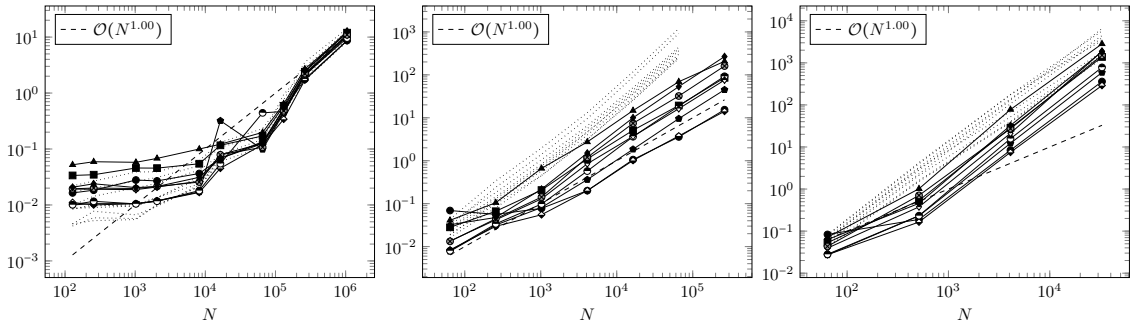
\begin{figure}
	\centering
	\resizebox{\textwidth}{!}{\input{img/AZAStimings1d-3d.tikz}}
	\caption{\label{fig:Stimings1d-3d}Timings in seconds of the direct sparse QR solver applied to the approximation problems of Figure~\ref{fig:AZtimings1d-3d}. In the black dashed line: $\mathcal O(N)$. The timings of Figure~\ref{fig:AZStimings1d-3d} are repeated in the black dotted lines.}
\end{figure}
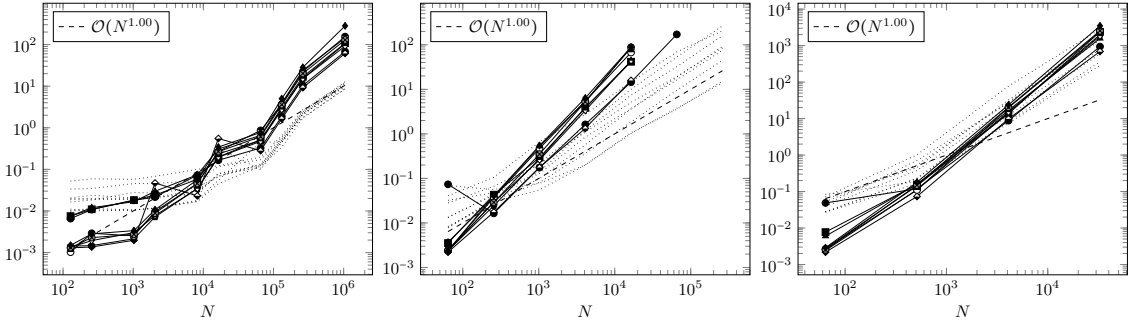

Figure~\ref{fig:AZStimings1d-3d} shows that the sparse AZ algorithm appears to be more efficient than the algorithms above, especially in the lower dimensions. In 3-D, the level of sparsity is not yet high enough to show a possible advantageous effect of using a sparse solver.

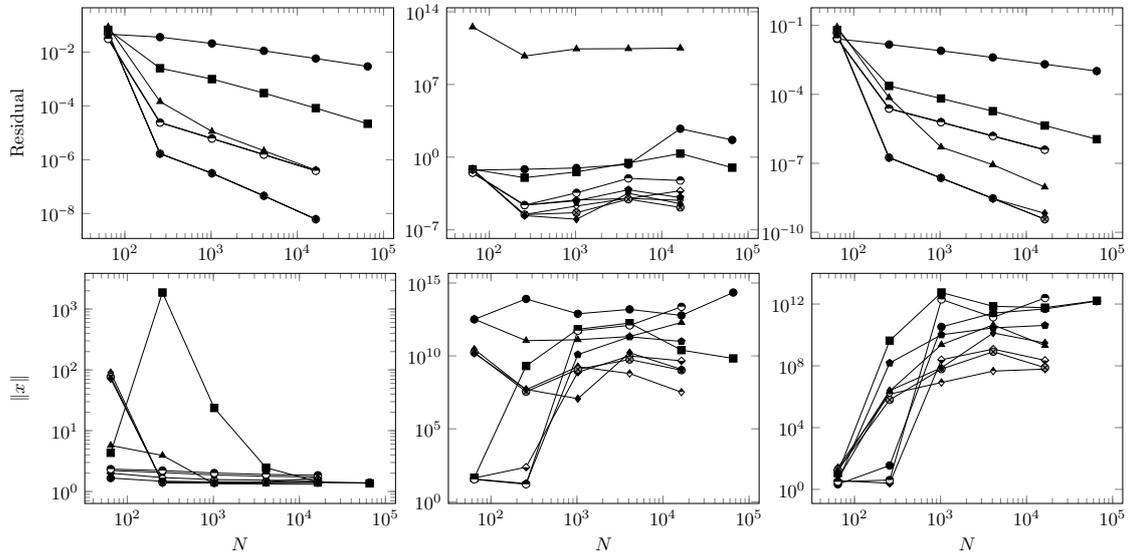
\begin{figure}
	\centering
	\resizebox{\textwidth}{!}{\input{img/smallresidual1d-3d_q4.tikz}}\\
	\resizebox{\textwidth}{!}{\input{img/smallcoefnorm1d-3d_q4.tikz}}
	\caption{\label{fig:AZerrors1d-3d}Residual (top) and coefficient norm (bottom) of 2-D problem in Figure \ref{fig:AZtimings1d-3d}, but with $\mathbf\osi=(4,4)$ instead of $\mathbf\osi=(2,2)$ using reduced AZ, sparse AZ and sparse QR (left to right). }
\end{figure}

%\begin{figure}
%	\centering
%	\resizebox{\textwidth}{!}{\input{img/smallresidual1d-3d.tikz}}\\
%	\resizebox{\textwidth}{!}{\input{img/smallcoefnorm1d-3d.tikz}}
%	\caption{\label{fig:AZerrors1d-3d}Residual (top) and coefficient norm (bottom) of 2-D problem in Figure \ref{fig:AZtimings1d-3d} with reduced AZ, sparse AZ and sparse QR (left to right). \MAGENTA{Vermeld de waarde van $\osi$}}
%\end{figure}

We also compare coefficient norm and residual of the different algorithms. For the coefficient norm, it is known for the vanilla and reduced AZ algorithms that a small-norm coefficient will be returned if it exists, since the solver in step 1 is closely related to a truncated SVD solver \cite{Adcock2019,adcock2017frames,az}. For the direct sparse QR solver of \cite{Davis2009}, no analogous error analysis is known, hence the coefficient norm might be large. While this was not the case for spline extension in \cite{coppe2019splines}, coefficients are indeed larger using wavelet extensions as shown in Figure~\ref{fig:AZerrors1d-3d}. Both the sparse AZ and a direct sparse QR solver are affected. This has a negative impact on the residual, as is also shown in Figure~\ref{fig:AZerrors1d-3d}.

\section{Adaptively smoothed wavelet AZ algorithm}\label{wavs:smooth}

In this final section we compare the use of wavelet-based extensions to the simpler setting of spline-based extensions.

The approximation space spanned by a basis of B-splines is exactly the same as that spanned by spline-based wavelets. Indeed, the wavelet transform is merely a change of basis. Therefore, the best approximations are the same:
\[
 \arg \min_{u \in \SPAN \Phi_N} \Vert f - u \Vert =  \arg \min_{u \in \SPAN \Psi_N} \Vert f - u \Vert.
\]
In typical applications of wavelets their compression properties play a major role. That is not really the case here: the wavelet least squares matrix $A$ is even somewhat less sparse than the corresponding matrix $\hat A$ using B-splines.

Another beneficial property of wavelets is the multiresolution nature of the approximation. In particular, unlike with B-splines, it is possible to associate different weighting factors with different scales. This enables the construction of bases for a range of function spaces with varying smoothness properties. For example, methods for the solution of partial differential equations employ wavelet bases for Sobolev spaces \cite{dahmen1997opeq}. In our setting, weighing different wavelet scales allows one to obtain smoother approximations. The increased smoothness is only visible in the extension $\Xi\setminus\Omega$ of the wavelet frame, since the approximation always resembles the function itself in the interior $\Omega$.

A smooth extension is not guaranteed by the methods described in \S\ref{wavs:az}. A least squares solver aims to minimize the residual of the system with a minimal norm solution. Therefore, the resulting wavelets coefficients do not necessarily decrease with increasing scale, even when approximating smooth functions \cite{adaptivity}. All coefficients have roughly similar size. In contrast, when approximating with a regular basis, the decrease of wavelet coefficient size is guaranteed for smooth functions, depending on the order of the multiresolution analysis at hand.
%A property that is taken advantage of in numerical approximation packages as \texttt{ApproxFun.jl} and \texttt{Chebfun.m}. 

\begin{algorithm}[h]
	\caption{The smoothed AZ algorithm}\label{alg:smoothedAZ}
	{\bf Input:} $A,Z \in \C^{M\times N}$, $b\in\C^M$, $W\in \R_+^{N\times N}$\\
	{\bf Output:} $x\in\C^N$ such that $Ax \approx b$
	\begin{algorithmic}[1]
		\State Solve $(I-AZ^*)AW x_1 = (I-AZ^*)b$ using a randomized low-rank solver
		\State $x_2 \gets Z^*(b-AWx_1)$
		\State $x \gets Wx_1 + x_2$
	\end{algorithmic}
\end{algorithm}

Smoothing can be introduced by switching to a weighted least squares formulation. We weigh the wavelet coefficients using the smoothed AZ algorithm (Algorithm~\ref{alg:smoothedAZ}). We simply add a diagonal weight matrix in the first step of the AZ algorithm, replacing $x_1$ by $W x_1$, and leave the other steps unchanged. The diagonal matrix has weights that depend on the scale of the corresponding entries of $x_1$. Much in the same way, the weighted reduced and sparse AZ algorithms can be formulated and implemented. Note that a weighted least squares problem is only being solved in step 1, not in step 2. The dual does not require modifications in this formulation.

We aim for a coefficient vector in which the coefficients decrease in size with increasing scale. The logic is as follows. Say a function is approximated on a coarse scale $\mathbf J-1$ with approximation error $e_1$. It can be expected that this approximation can be refined by adding wavelet coefficients on finer scales with size on the order of $\mathcal O(e_1)$. We ensure that the coefficients on the finer scale have that size simply by choosing the corresponding diagonal entries of $W$ equal to $e_1$. The function is now approximated on the finer scale, say with an approximation error $e_2$. This error can be used to weight the next scale of coefficients, and so on. Thus, we obtain a diagonal weighting operator parametrized with the weights $e_1,\dots,e_L$ and size $\mathbf N$, $L\leq \log_2(N_i)$, $i=1,\dots,d$:
\begin{equation}
	W(\mathbf k,\mathbf l; [e_1, e_2, \dots, e_L], \mathbf N) = \delta_{\mathbf k,\mathbf l}e_{\min\{1,L-\min_{i=1,\dots,d}\{\log_2(N_i)-L\}\}}\qquad  \mathbf k,\mathbf l \in I_{\mathbf N}.
\end{equation}
We use the size of the right hand side as an initial weight for the first approximation. This way we arrive at the adaptive Algorithm~\ref{alg:adaptivelysmoothedAZ}. 

\begin{algorithm}[h]
	\caption{The adaptively weighted AZ algorithm}\label{alg:adaptivelysmoothedAZ}
	{\bf Input:} $\mathbf N$, $\mathbf\osi$, $f$, $\Omega\subset\Xi$, wavelet types\\
	{\bf Output:} $x$ (the wavelet extension coefficients)
	\begin{algorithmic}[1]
		\State $\mathbf n\gets\mathbf N-\min(N)+1$
		\State $e \gets \|b^{\mathbf\osi}_{\mathbf n}\|$
		\While{$\mathbf n < \mathbf N$}
		\State $W \gets W(\cdot,\cdot; e, \mathbf N)$
		\State $x\gets$ Apply the smoothed AZ algorithm with $A^{\mathbf\osi}_{\mathbf n}$, $Z^{\mathbf\osi}_{\mathbf n}$, $b^{\mathbf\osi}_{\mathbf n}$, $W$
		\State $e \gets [e; \|A^{\mathbf\osi}_{\mathbf n}x-b^{\mathbf\osi}_{\mathbf n}\|]$
		\State $\mathbf n \gets 2\mathbf n$
		\EndWhile
		
	\end{algorithmic}
\end{algorithm}

\begin{figure}
	\begin{center}
		\resizebox{.8\columnwidth}{!}{\input{img/1dadaptivity.tikz}}%
	\end{center}
	\caption{Wavelet extension approximation of $f(x)=e^{x}$ on $[0,0.6]$ using \texttt{cdf33}, $N=256$, $\osi=2$.  Approximation and extension (left) and coefficient size (right) when using the adaptively weighted reduced AZ algorithm (red), a pivoted QR (blue), the reduced AZ algorithm (brown).}\label{fig:1Dcoefs}
\end{figure}
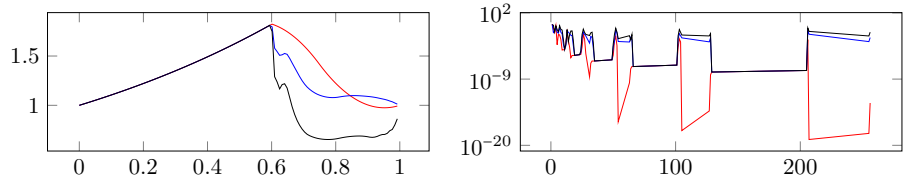

\begin{figure}
	\begin{center}
		\includegraphics[width=.3\linewidth]{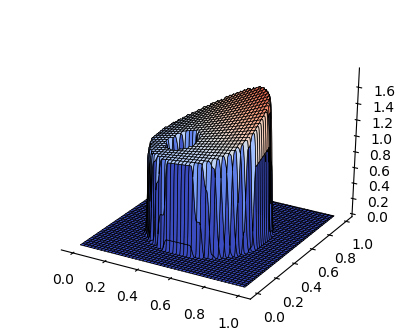}
		\hspace{.3\linewidth}
		\includegraphics[width=.3\linewidth]{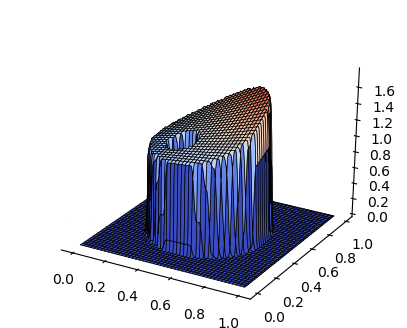}
		\\\vspace{.5em}
		\includegraphics[width=.3\linewidth]{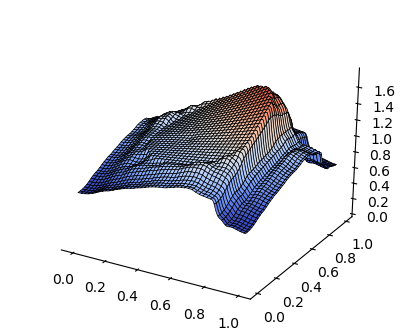}
		\includegraphics[width=.3\linewidth]{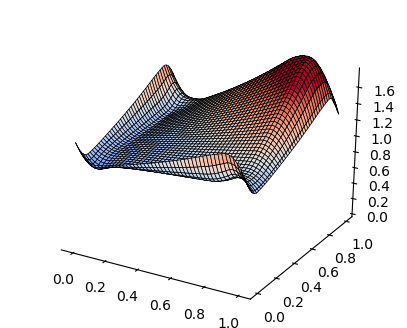}
		\includegraphics[width=.3\linewidth]{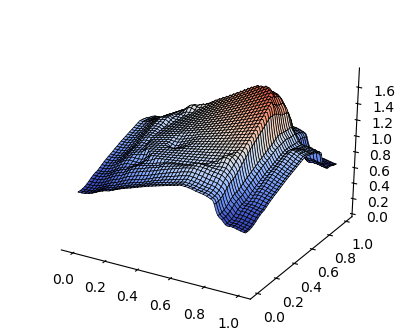}\\\vspace{.5em}
		\includegraphics[width=.3\linewidth]{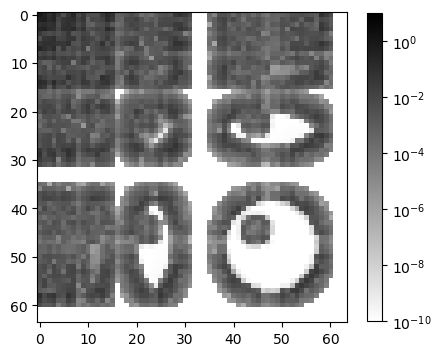}
		\includegraphics[width=.3\linewidth]{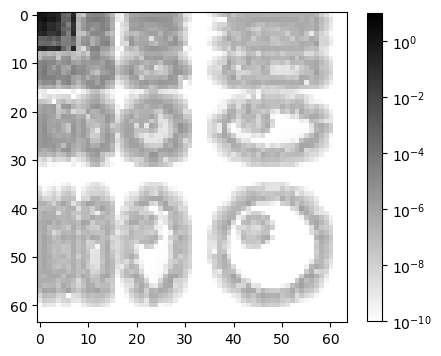}
		\includegraphics[width=.3\linewidth]{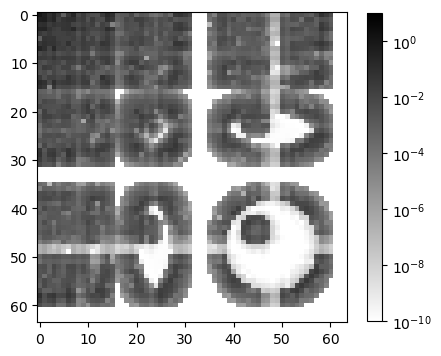}
	\end{center}
	\caption{\label{fig:2Dcoefs}Approximation of $f(x,y)=e^{xy}$ on $[0,1/2]^2$, $\mathbf N=[64,64], \mathbf \osi=[4,4]$. Top row: Cubic spline extension approximation. Mid and bottom row: \texttt{cdf33} wavelet extension approximation. Top and mid row: The approximant on $[0,1]^2$. Bottom row: the coefficient size. Left: The reduced AZ algorithm . Mid: The adaptively weighted reduced AZ algorithm. Right: A pivoted QR.}
\end{figure}
\begin{table}
	\centering
	\input{weightedtabular.tikz}
	\caption{\label{tab:comparison}Wavelet coefficient norms and residual error. }
\end{table}

We illustrate the adaptive Algorithm with a simple example in 1D first. In Figure~\ref{fig:1Dcoefs} we compare to a pivoted QR and the reduced AZ algorithm (both without smoothing). The function $f(x)=e^x$ is approximated on the interval $[0,0.6]$ using a \texttt{cdf33} wavelet extension. The weighted algorithms very clearly leads to the smoothest extension (shown in the left panel). In the right panel, the sizes of the wavelet coefficients are compared. The approximation domain $[0,0.6]$ is visible in all levels of the wavelet coefficients. There, all coefficients have rougly similar size (the three lines overlap) because they approximate the same function. The coefficients differ in the extension $[0.6,1]$. There, it is seen that the smoothed approximation (red line) yields significantly smaller coefficients than the non-smoothed approximations.
%\BLUE{The coefficients of wavelet basis functions in the interior are of the same scale for all three methods. The coefficients in the exterior are however much smaller for the method using weighting than for the other two. This leads to the rugged appearance of the extension. Note that the weighing in the AZ algorithm only effects coefficients of wavelet functions on the boundary.} \MAGENTA{Te overlopen, ik begrijp de figuur niet helemaal.}

In Figure~\ref{fig:2Dcoefs} and Table \ref{tab:comparison} we compare the adaptive Algorithm \ref{alg:adaptivelysmoothedAZ} using both sparse (only in table) and reduced AZ with the original non-smoothed reduced and sparse AZ algorithm. Also we compare with a simple pivoted QR. Finally, we compare wavelet extension with spline extension approximation as in \cite{coppe2019splines} where no adaptive weighting is possible. 

The smoothed approximation (the middle column of Figure~\ref{fig:2Dcoefs} and second row of the table) leads to a visually smoother extension. We also see that the non-smoothed methods tend to result in larger wavelet coefficients at the boundary. For the smoothed result, one can also see that the wavelet coefficients at a coarse scale are larger then those at a finer scale.
The approximant itself also takes another shape. The spline extension approximant will drop down to zero outside of $\Omega$. Choosing a wavelet extension instead of a spline extension has no effect on the residual, nor does weighting. As expected, the pivoted QR leads to the smallest coefficient norm. 

\section{Concluding remarks}\label{wavs:conclusion}

We have shown that wavelet approximation on general domains is possible and efficient using regular wavelets defined on a bounding box. Compared to existing wavelet literature, we have had to compute a discrete dual scaling sequence. Apart from this construction, with values listed exhaustively in the appendix, well-known scaling functions and wavelets could be used.

The ill-conditioning of the extension problem necessitates a least squares formulation with regularization. The proposed algorithms have a complexity of $\mathcal O(N)$ operations in 1-D, $\mathcal O(N^{3/2})$ in 2-D and $\mathcal O(N^{3(d-1)/d})$ in $d$-D, $d>1$. This should be compared to the cubic cost of standard direct solvers. The algorithms required a detailed study of the sparsity structure of all the matrices involved.

The use of a general sparse QR function does not require such detailed study of the structure. As such, it is simpler to implement. The results experimentally appear to be more efficient. However, it also seems less stable.
%
%\RED{\begin{enumerate}
%	\item Both the vanilla AZ and the reduced AZ (using a low-rank solver in the first step)) attain $\mathcal O(N)$ operations in 1-D, $\mathcal O(N^{3/2})$ in 2-D and $\mathcal O(N^{3(d-1)/d})$ in $d$-D, $d>1$. The first one because of a lucky accident. Instead of the cubic complexity of generic least squares solver.
%	\item Sparse AZ appears experimentally more efficient, but	is less stable.
%	\item We had to introduce a discrete dual,
%	\item Wavelet function approximation on generic domains is possible using wavelet extension frames.
%\end{enumerate}}

\appendix
\section{Appendix}

The values of the discrete duals that were used for the experiments in this paper are tabulated in Table~\ref{tab:discretevalues}.

%\begin{table}
%	\input{dualtabular.tikz}
%	\caption{\label{tab:discretevalues} Values of sampled father function and compactly supported dual with minimal support and $\osi=2$, i.e., $b^2_k$ (odd rows) and $\tilde b^2_k$ (even rows). The column \emph{offset} indicates the smallest $k$ for which  $b^2_k$ ($\tilde b^2_k$) is non-zero. The column \emph{values} holds the non-zero values of $b^2_k$ ($\tilde b^2_k$). Compare them with the sequences shown in Figures~\ref{fig:compactduals_db} and~\ref{fig:compactduals_cdf} that show them for $\osi=4$.}
%\end{table}

\begin{table}
	\begin{tabular}{p{0.5\textwidth}@{\hskip -10pt}p{0.5\textwidth}}
		\vspace{0pt}  \resizebox{.47\textwidth}{!}{\input{verticaldualtabulardb2.tikz}}&
		\vspace{0pt} \resizebox{.47\textwidth}{!}{\input{verticaldualtabularcdf31.tikz}}\\
		\vspace{0pt}  \resizebox{.47\textwidth}{!}{\input{verticaldualtabulardb3.tikz}}&
		\vspace{0pt} \resizebox{.47\textwidth}{!}{\input{verticaldualtabularcdf42.tikz}}\\
		\vspace{0pt} 
		\resizebox{.47\textwidth}{!}{\input{verticaldualtabulardb4.tikz}}&
	\end{tabular}
	\caption{Values of sampled father function and compactly supported dual with minimal support and $\osi=2$, i.e., $b^2_k$ (second column) and $\tilde b^2_k$ (third column). The second and third column hold the non-zero values of $b^2_k$ $\tilde b^2_k$, respectively. Compare them with the sequences shown in Figures~\ref{fig:compactduals_db} and~\ref{fig:compactduals_cdf} that show them for $\osi=4$.}\label{tab:discretevalues}
\end{table}

\bibliographystyle{siam}
\bibliography{waveletextension}
\end{document}

%% file: img/db2.tikz
\begin{tikzpicture}
\begin{axis}[width={.5\textwidth}, height={.25\textwidth}]
    \addplot+[mark={none}]
        table[row sep={\\}]
        {
            \\
            0.0  -0.0  \\
            0.0078125  0.0947240455194615  \\
            0.015625  0.13868562803742576  \\
            0.0234375  0.16406685953817435  \\
            0.03125  0.20304985200598893  \\
            0.0390625  0.2148293225396716  \\
            0.046875  0.2402105540404202  \\
            0.0546875  0.27057037559131847  \\
            0.0625  0.2972856162754504  \\
            0.0703125  0.2991079067088336  \\
            0.078125  0.3145319581092828  \\
            0.0859375  0.3349345995598816  \\
            0.09375  0.3516926601437141  \\
            0.1015625  0.37573988246107937  \\
            0.109375  0.39614252391167826  \\
            0.1171875  0.41521115617889387  \\
            0.125  0.4352563509461097  \\
            0.1328125  0.42978947964596004  \\
            0.140625  0.4379243693128761  \\
            0.1484375  0.45103784902994204  \\
            0.15625  0.46050674788024154  \\
            0.1640625  0.47726480846407393  \\
            0.171875  0.4903782881811397  \\
            0.1796875  0.5021577587148225  \\
            0.1875  0.5149137917485052  \\
            0.1953125  0.534339870699104  \\
            0.203125  0.5501213687829364  \\
            0.2109375  0.5645688576833854  \\
            0.21875  0.5799929090838347  \\
            0.2265625  0.5934638354842838  \\
            0.234375  0.6079113243847329  \\
            0.2421875  0.6227162599685653  \\
            0.25  0.6372595264191646  \\
            0.2578125  0.6264566183854817  \\
            0.265625  0.6292554713188652  \\
            0.2734375  0.637032914302398  \\
            0.28125  0.6411657764191645  \\
            0.2890625  0.6525878002694638  \\
            0.296875  0.6603652432529968  \\
            0.3046875  0.6668086770531465  \\
            0.3125  0.6742286733532964  \\
            0.3203125  0.6883187155703623  \\
            0.328125  0.6987641769206617  \\
            0.3359375  0.7078756290875776  \\
            0.34375  0.717963643754494  \\
            0.3515625  0.7260985334214101  \\
            0.359375  0.7352099855883264  \\
            0.3671875  0.7446788844386255  \\
            0.375  0.7538861141556918  \\
            0.3828125  0.7699292813727578  \\
            0.390625  0.7823278677230572  \\
            0.3984375  0.7933924448899733  \\
            0.40625  0.8054335845568895  \\
            0.4140625  0.8155215992238058  \\
            0.421875  0.8265861763907217  \\
            0.4296875  0.8380082002410213  \\
            0.4375  0.8491685549580873  \\
            0.4453125  0.858541676258237  \\
            0.453125  0.8688913600583866  \\
            0.4609375  0.8795984905419193  \\
            0.46875  0.890043951892219  \\
            0.4765625  0.9010127515089853  \\
            0.484375  0.9117198819925183  \\
            0.4921875  0.9223312349259012  \\
            0.5  0.9330127018922194  \\
            0.5078125  0.9183035438585364  \\
            0.515625  0.9171961467919199  \\
            0.5234375  0.9210673397754527  \\
            0.53125  0.9212939518922194  \\
            0.5390625  0.9288097257425187  \\
            0.546875  0.9326809187260516  \\
            0.5546875  0.9352181025262014  \\
            0.5625  0.9387318488263512  \\
            0.5703125  0.948915641043417  \\
            0.578125  0.9554548523937165  \\
            0.5859375  0.9606600545606325  \\
            0.59375  0.9668418192275489  \\
            0.6015625  0.9710704588944647  \\
            0.609375  0.9762756610613812  \\
            0.6171875  0.9818383099116805  \\
            0.625  0.9871392896287469  \\
            0.6328125  0.9992762068458129  \\
            0.640625  1.007768543196112  \\
            0.6484375  1.014926870363028  \\
            0.65625  1.0230617600299445  \\
            0.6640625  1.0292435246968608  \\
            0.671875  1.0364018518637765  \\
            0.6796875  1.043917625714076  \\
            0.6875  1.051171730431142  \\
            0.6953125  1.0566386017312919  \\
            0.703125  1.0630820355314414  \\
            0.7109375  1.069882916014974  \\
            0.71875  1.0764221273652739  \\
            0.7265625  1.08348467698204  \\
            0.734375  1.090285557465573  \\
            0.7421875  1.096990660398956  \\
            0.75  1.103765877365274  \\
            0.7578125  1.1173325813158728  \\
            0.765625  1.1272547043997054  \\
            0.7734375  1.1358428183001543  \\
            0.78125  1.1454074947006037  \\
            0.7890625  1.1530190461010528  \\
            0.796875  1.1616071600015019  \\
            0.8046875  1.1705527205853343  \\
            0.8125  1.1792366120359332  \\
            0.8203125  1.1861332700696159  \\
            0.828125  1.1940064906032988  \\
            0.8359375  1.202237157820364  \\
            0.84375  1.2102061559041966  \\
            0.8515625  1.2186984922544961  \\
            0.859375  1.2269291594715621  \\
            0.8671875  1.235064049138478  \\
            0.875  1.2432690528383288  \\
            0.8828125  1.2496423726055443  \\
            0.890625  1.25699225487276  \\
            0.8984375  1.2646995838233588  \\
            0.90625  1.2721452436407241  \\
            0.9140625  1.2801142417245568  \\
            0.921875  1.287821570675155  \\
            0.9296875  1.2954331220756046  \\
            0.9375  1.3031147875089877  \\
            0.9453125  1.3112753406931195  \\
            0.953125  1.3191742247440177  \\
            0.9609375  1.326977331244766  \\
            0.96875  1.334850551778449  \\
            0.9765625  1.3425835442462632  \\
            0.984375  1.3503866507470117  \\
            0.9921875  1.3582154207649757  \\
            1.0  1.3660254037844388  \\
            1.0078125  1.1843898127455155  \\
            1.015625  1.104279147709587  \\
            1.0234375  1.0613291847080897  \\
            1.03125  0.9911756997724608  \\
            1.0390625  0.9754292587050954  \\
            1.046875  0.9324792957035982  \\
            1.0546875  0.8795721526018015  \\
            1.0625  0.833954171233538  \\
            1.0703125  0.8381220903667713  \\
            1.078125  0.8150864875658731  \\
            1.0859375  0.7820937046646751  \\
            1.09375  0.7563900834970104  \\
            1.1015625  0.7161081388622795  \\
            1.109375  0.683115355961082  \\
            1.1171875  0.6527905914266506  \\
            1.125  0.6205127018922194  \\
            1.1328125  0.6392589444925191  \\
            1.140625  0.6308016651586864  \\
            1.1484375  0.6123872057245546  \\
            1.15625  0.6012619080239557  \\
            1.1640625  0.575558286856291  \\
            1.171875  0.5571438274221588  \\
            1.1796875  0.5413973863547936  \\
            1.1875  0.5236978202874283  \\
            1.1953125  0.49265816238623045  \\
            1.203125  0.46890766621856556  \\
            1.2109375  0.44782518841766716  \\
            1.21875  0.42478958561676905  \\
            1.2265625  0.40566023281587055  \\
            1.234375  0.3845777550149723  \\
            1.2421875  0.36278038384730743  \\
            1.25  0.34150635094610976  \\
            1.2578125  0.37092466701347515  \\
            1.265625  0.3731394611467088  \\
            1.2734375  0.3653970751796427  \\
            1.28125  0.36494385094610976  \\
            1.2890625  0.34991230324551087  \\
            1.296875  0.342169917278445  \\
            1.3046875  0.3370955496781455  \\
            1.3125  0.3300680570778461  \\
            1.3203125  0.3097004726437142  \\
            1.328125  0.29662204994311536  \\
            1.3359375  0.28621164560928275  \\
            1.34375  0.2738481162754504  \\
            1.3515625  0.265390836941618  \\
            1.359375  0.25498043260778563  \\
            1.3671875  0.24385513490718663  \\
            1.375  0.23325317547305485  \\
            1.3828125  0.20897934103892313  \\
            1.390625  0.19199466833832418  \\
            1.3984375  0.17767801400449176  \\
            1.40625  0.16140823467065932  \\
            1.4140625  0.14904470533682698  \\
            1.421875  0.1347280510029944  \\
            1.4296875  0.11969650330239567  \\
            1.4375  0.10518829386826381  \\
            1.4453125  0.09425455126796435  \\
            1.453125  0.08136768366766489  \\
            1.4609375  0.06776592270059892  \\
            1.46875  0.05468750000000013  \\
            1.4765625  0.04056240076646711  \\
            1.484375  0.026960639799401262  \\
            1.4921875  0.01355043393263462  \\
            1.5  0.0  \\
            1.5078125  0.03723081606736548  \\
            1.515625  0.04725811020059895  \\
            1.5234375  0.04732822423353307  \\
            1.53125  0.05468750000000014  \\
            1.5390625  0.047468452299401186  \\
            1.546875  0.04753856633233525  \\
            1.5546875  0.05027669873203583  \\
            1.5625  0.05106170613173636  \\
            1.5703125  0.03850662169760449  \\
            1.578125  0.033240698997005595  \\
            1.5859375  0.030642794663173187  \\
            1.59375  0.026091765329340825  \\
            1.6015625  0.025446985995508372  \\
            1.609375  0.022849081661676002  \\
            1.6171875  0.019536283961077077  \\
            1.625  0.016746824526945207  \\
            1.6328125  0.0002854900928133846  \\
            1.640625  -0.008886682607785511  \\
            1.6484375  -0.015390836941617912  \\
            1.65625  -0.023848116275450307  \\
            1.6640625  -0.028399145609282733  \\
            1.671875  -0.034903299943115154  \\
            1.6796875  -0.04212234764371403  \\
            1.6875  -0.048818057077845894  \\
            1.6953125  -0.05193929967814532  \\
            1.703125  -0.057013667278444764  \\
            1.7109375  -0.06280292824551072  \\
            1.71875  -0.0680688509461096  \\
            1.7265625  -0.07438145017964257  \\
            1.734375  -0.08017071114670848  \\
            1.7421875  -0.08576841701347494  \\
            1.75  -0.09150635094610964  \\
            1.7578125  -0.11082725884730737  \\
            1.765625  -0.12285900501497221  \\
            1.7734375  -0.13222273281587052  \\
            1.78125  -0.14353958561676886  \\
            1.7890625  -0.1509501884176672  \\
            1.796875  -0.16031391621856553  \\
            1.8046875  -0.1703925373862303  \\
            1.8125  -0.17994782028742817  \\
            1.8203125  -0.1859286363547935  \\
            1.828125  -0.19386257742215887  \\
            1.8359375  -0.20251141185629065  \\
            1.84375  -0.21063690802395554  \\
            1.8515625  -0.21980908072455438  \\
            1.859375  -0.22845791515868624  \\
            1.8671875  -0.23691519449251858  \\
            1.875  -0.24551270189221924  \\
            1.8828125  -0.2504468414266505  \\
            1.890625  -0.2573341059610818  \\
            1.8984375  -0.2649362638622796  \\
            1.90625  -0.27201508349701037  \\
            1.9140625  -0.2801405796646752  \\
            1.921875  -0.2877427375658729  \\
            1.9296875  -0.29515334036677127  \\
            1.9375  -0.3027041712335377  \\
            1.9453125  -0.31121277760180144  \\
            1.953125  -0.3191980457035981  \\
            1.9609375  -0.3269917587050952  \\
            1.96875  -0.3349256997724606  \\
            1.9765625  -0.3425791847080897  \\
            1.984375  -0.3503728977095869  \\
            1.9921875  -0.35821793774551536  \\
            2.0  -0.36602540378443876  \\
            2.0078125  -0.2791138582649771  \\
            2.015625  -0.24296477574701283  \\
            2.0234375  -0.22539604424626422  \\
            2.03125  -0.19422555177844975  \\
            2.0390625  -0.19025858124476705  \\
            2.046875  -0.1726898497440184  \\
            2.0546875  -0.1501425281931201  \\
            2.0625  -0.13123978750898824  \\
            2.0703125  -0.137229997075605  \\
            2.078125  -0.12961844567515585  \\
            2.0859375  -0.11702830422455691  \\
            2.09375  -0.10808274364072454  \\
            2.1015625  -0.09184802132335915  \\
            2.109375  -0.0792578798727603  \\
            2.1171875  -0.06800174760554464  \\
            2.125  -0.05576905283832898  \\
            2.1328125  -0.06904842413847874  \\
            2.140625  -0.0687260344715625  \\
            2.1484375  -0.06342505475449658  \\
            2.15625  -0.061768655904197146  \\
            2.1640625  -0.05282309532036475  \\
            2.171875  -0.04752211560329879  \\
            2.1796875  -0.04355514506961611  \\
            2.1875  -0.03861161203593342  \\
            2.1953125  -0.026998033085334546  \\
            2.203125  -0.019029035001502128  \\
            2.2109375  -0.012394046101052948  \\
            2.21875  -0.0047824947006038085  \\
            2.2265625  0.0008759316998453607  \\
            2.234375  0.007510920600294511  \\
            2.2421875  0.014503356184126914  \\
            2.25  0.021234122634725853  \\
            2.2578125  0.0026187146010431496  \\
            2.265625  -0.002394932465573625  \\
            2.2734375  -0.002429989482040648  \\
            2.28125  -0.006109627365274166  \\
            2.2890625  -0.0025001035149747157  \\
            2.296875  -0.002535160531441762  \\
            2.3046875  -0.0039042267312920376  \\
            2.3125  -0.004296730431142306  \\
            2.3203125  0.001980811785923615  \\
            2.328125  0.004613773136223055  \\
            2.3359375  0.005912725303139264  \\
            2.34375  0.008188239970055468  \\
            2.3515625  0.008510629636971669  \\
            2.359375  0.009809581803887876  \\
            2.3671875  0.011465980654187321  \\
            2.375  0.012860710371253253  \\
            2.3828125  0.021091377588319176  \\
            2.390625  0.025677463938618614  \\
            2.3984375  0.028929541105534817  \\
            2.40625  0.033158180772451036  \\
            2.4140625  0.035433695439367224  \\
            2.421875  0.03868577260628343  \\
            2.4296875  0.042295296456582865  \\
            2.4375  0.045643151173648804  \\
            2.4453125  0.04720377247379851  \\
            2.453125  0.049740956273948235  \\
            2.4609375  0.05263558675748119  \\
            2.46875  0.055268548107780625  \\
            2.4765625  0.05842484772454712  \\
            2.484375  0.06131947820808007  \\
            2.4921875  0.0641183311414633  \\
            2.5  0.06698729810778069  \\
            2.5078125  0.04446564007409797  \\
            2.515625  0.03554574300748122  \\
            2.5234375  0.03160443599101418  \\
            2.53125  0.024018548107780667  \\
            2.5390625  0.023721821958080116  \\
            2.546875  0.019780514941613078  \\
            2.5546875  0.014505198741762806  \\
            2.5625  0.010206445041912527  \\
            2.5703125  0.012577737258978454  \\
            2.578125  0.011304448609277902  \\
            2.5859375  0.008697150776194104  \\
            2.59375  0.007066415443110308  \\
            2.6015625  0.0034825551100265153  \\
            2.609375  0.0008752572769427233  \\
            2.6171875  -0.001374593872757828  \\
            2.625  -0.0038861141556919086  \\
            2.6328125  0.0004383030613740235  \\
            2.640625  0.0011181394116734661  \\
            2.6484375  0.00046396657858967135  \\
            2.65625  0.0007863562455058737  \\
            2.6640625  -0.0008443790875779197  \\
            2.671875  -0.001498551920661716  \\
            2.6796875  -0.0017952780703622705  \\
            2.6875  -0.002353673353296344  \\
            2.6953125  -0.004699302053146619  \\
            2.703125  -0.006068368252996899  \\
            2.7109375  -0.007079987769463933  \\
            2.71875  -0.008353276419164488  \\
            2.7265625  -0.009103226802398003  \\
            2.734375  -0.010114846318865035  \\
            2.7421875  -0.011222243385481796  \\
            2.75  -0.012259526419164493  \\
            2.7578125  -0.006505322468565599  \\
            2.765625  -0.004395699384733191  \\
            2.7734375  -0.0036200854842840234  \\
            2.78125  -0.0018679090838348569  \\
            2.7890625  -0.0020688576833856895  \\
            2.796875  -0.0012932437829365215  \\
            2.8046875  -0.0001601831991041142  \\
            2.8125  0.0007112082514947767  \\
            2.8203125  -0.0002046337148225375  \\
            2.828125  -0.00014391318113985126  \\
            2.8359375  0.00027425403592607525  \\
            2.84375  0.00043075211975848396  \\
            2.8515625  0.001110588470057929  \\
            2.859375  0.0015287556871238553  \\
            2.8671875  0.0018511453540400583  \\
            2.875  0.002243649053890338  \\
            2.8828125  0.0008044688211059872  \\
            2.890625  0.0003418510883216371  \\
            2.8984375  0.00023668003892052754  \\
            2.90625  -0.0001301601437141001  \\
            2.9140625  2.6337940118308555e-5  \\
            2.921875  -7.883310928280097e-5  \\
            2.9296875  -0.000279781708833633  \\
            2.9375  -0.00041061627545039225  \\
            2.9453125  -6.25630913185385e-5  \\
            2.953125  2.382095957979702e-5  \\
            2.9609375  1.4427460328409997e-5  \\
            2.96875  7.514799401109603e-5  \\
            2.9765625  -4.359538174363997e-6  \\
            2.984375  -1.3753037425750997e-5  \\
            2.9921875  2.516980538511502e-6  \\
            3.0  -0.0  \\
        }
        ;
    \addplot+[mark={none}]
        table[row sep={\\}]
        {
            \\
            -1.0  0.0  \\
            -0.9921875  -0.025381231500748606  \\
            -0.984375  -0.037160702034431295  \\
            -0.9765625  -0.04396158251796425  \\
            -0.96875  -0.0544070438682637  \\
            -0.9609375  -0.05756334348503018  \\
            -0.953125  -0.06436422396856313  \\
            -0.9453125  -0.07249911363547933  \\
            -0.9375  -0.07965744080239555  \\
            -0.9296875  -0.08014572205239553  \\
            -0.921875  -0.08427858416916204  \\
            -0.9140625  -0.08974545546931172  \\
            -0.90625  -0.09423576426946145  \\
            -0.8984375  -0.10067919806961116  \\
            -0.890625  -0.10614606936976088  \\
            -0.8828125  -0.11125549398652734  \\
            -0.875  -0.11662658773652741  \\
            -0.8671875  -0.11516174398652741  \\
            -0.859375  -0.11734148110329387  \\
            -0.8515625  -0.12085522740344358  \\
            -0.84375  -0.12339241120359332  \\
            -0.8359375  -0.12788272000374304  \\
            -0.828125  -0.13139646630389273  \\
            -0.8203125  -0.13455276592065923  \\
            -0.8125  -0.13797073467065923  \\
            -0.8046875  -0.14317593683757543  \\
            -0.796875  -0.14740457650449162  \\
            -0.7890625  -0.15127576948802454  \\
            -0.78125  -0.15540863160479107  \\
            -0.7734375  -0.1590181554550905  \\
            -0.765625  -0.16288934843862343  \\
            -0.7578125  -0.1668563189723061  \\
            -0.75  -0.1707531754730548  \\
            -0.7421875  -0.16785854498952182  \\
            -0.734375  -0.16860849537275535  \\
            -0.7265625  -0.1706924549393721  \\
            -0.71875  -0.17179985200598885  \\
            -0.7109375  -0.17486037407260563  \\
            -0.703125  -0.17694433363922235  \\
            -0.6953125  -0.17867084652245588  \\
            -0.6875  -0.18065902853892296  \\
            -0.6796875  -0.18443444397230616  \\
            -0.671875  -0.1872332969056894  \\
            -0.6640625  -0.18967470315568932  \\
            -0.65625  -0.1923777785389229  \\
            -0.6484375  -0.19455751565568935  \\
            -0.640625  -0.19699892190568935  \\
            -0.6328125  -0.19953610570583902  \\
            -0.625  -0.20200317547305477  \\
            -0.6171875  -0.20630192917290507  \\
            -0.609375  -0.2096241203727553  \\
            -0.6015625  -0.21258886488922232  \\
            -0.59375  -0.2158152785389229  \\
            -0.5859375  -0.21851835392215646  \\
            -0.578125  -0.22148309843862335  \\
            -0.5703125  -0.22454362050524018  \\
            -0.5625  -0.2275340285389229  \\
            -0.5546875  -0.23004554882185693  \\
            -0.546875  -0.23281873823802451  \\
            -0.5390625  -0.23568770520434174  \\
            -0.53125  -0.23848655813772504  \\
            -0.5234375  -0.2414256391369764  \\
            -0.515625  -0.2442946061032937  \\
            -0.5078125  -0.24713790955239529  \\
            -0.5  -0.24999999999999997  \\
            -0.49218749999999994  -0.2460586929835329  \\
            -0.484375  -0.2457619668338324  \\
            -0.47656250000000006  -0.24679924986751506  \\
            -0.46874999999999994  -0.24685997040119778  \\
            -0.4609375  -0.24887381593488045  \\
            -0.45312500000000006  -0.24991109896856312  \\
            -0.44531249999999994  -0.25059093531886256  \\
            -0.4375  -0.2515324408023956  \\
            -0.42968750000000006  -0.25426117970284473  \\
            -0.42187499999999994  -0.2560133561032939  \\
            -0.4140625  -0.25740808582035973  \\
            -0.40625000000000006  -0.2590644846706593  \\
            -0.39843749999999994  -0.2601975452544916  \\
            -0.390625  -0.26159227497155757  \\
            -0.38281250000000006  -0.2630827822387732  \\
            -0.37499999999999994  -0.2645031754730549  \\
            -0.3671875  -0.26775525263997113  \\
            -0.35937500000000006  -0.2700307673068873  \\
            -0.35156249999999994  -0.27194883529042024  \\
            -0.34375  -0.27412857240718674  \\
            -0.33593750000000006  -0.27578497125748613  \\
            -0.32812499999999994  -0.277703039241019  \\
            -0.3203125  -0.2797168847747018  \\
            -0.31250000000000006  -0.28166061627545036  \\
            -0.30468749999999994  -0.2831254600254504  \\
            -0.296875  -0.28485197290868386  \\
            -0.28906250000000006  -0.286674263342067  \\
            -0.28124999999999994  -0.2884264397425163  \\
            -0.2734375  -0.2903188442088336  \\
            -0.26562500000000006  -0.29214113464221686  \\
            -0.25781249999999994  -0.29393776155838436  \\
            -0.25  -0.29575317547305485  \\
            -0.2421874999999999  -0.29938836284057  \\
            -0.2343750000000001  -0.30204698770808514  \\
            -0.2265625  -0.3043481658922169  \\
            -0.2187499999999999  -0.30691101320958225  \\
            -0.2109375000000001  -0.3089505222604806  \\
            -0.203125  -0.3112517004446124  \\
            -0.1953124999999999  -0.31364865617889404  \\
            -0.1875000000000001  -0.3159754978802416  \\
            -0.1796875  -0.31782345183084043  \\
            -0.1718749999999999  -0.3199330749146728  \\
            -0.1640625000000001  -0.3221384755486549  \\
            -0.15625  -0.324273762149703  \\
            -0.1484374999999999  -0.3265492768166191  \\
            -0.1406250000000001  -0.3287546774506013  \\
            -0.1328125  -0.33093441456736766  \\
            -0.12499999999999989  -0.3331329386826371  \\
            -0.11718750000000011  -0.33484066456736783  \\
            -0.109375  -0.33681005958533206  \\
            -0.10156249999999989  -0.33887523215344606  \\
            -0.09375000000000011  -0.34087029068862595  \\
            -0.0859375  -0.34300557728967407  \\
            -0.07812499999999989  -0.3450707498577879  \\
            -0.07031250000000011  -0.3471102589086863  \\
            -0.0625  -0.3491685549580875  \\
            -0.05468749999999989  -0.35135516859356675  \\
            -0.04687500000000011  -0.35347166819611203  \\
            -0.0390625  -0.3555625042814416  \\
            -0.03124999999999989  -0.35767212736527404  \\
            -0.02343750000000011  -0.3597441764521009  \\
            -0.015625  -0.3618350125374305  \\
            -0.007812499999999889  -0.3639327251414729  \\
            -1.1102230246251565e-16  -0.3660254037844387  \\
            0.0078125  -0.14150996377695252  \\
            0.01562500000000011  -0.038433809973059785  \\
            0.02343749999999989  0.020192480062868723  \\
            0.03125  0.11135832876646204  \\
            0.03906250000000011  0.1374450601347257  \\
            0.04687499999999989  0.1960713501706543  \\
            0.0546875  0.26660794530688214  \\
            0.06250000000000011  0.3284255919760444  \\
            0.07031249999999989  0.33069171314370904  \\
            0.078125  0.3654973929790387  \\
            0.08593750000000011  0.4122133779146677  \\
            0.09374999999999989  0.45021041438323095  \\
            0.1015625  0.5056453477859258  \\
            0.10937500000000011  0.5523613327215549  \\
            0.11718749999999989  0.5958859610239502  \\
            0.125  0.6417468245269452  \\
            0.1328125000000001  0.6265750487604781  \\
            0.1406249999999999  0.6439428316616758  \\
            0.1484375  0.6732209196631729  \\
            0.1562500000000001  0.6937800591976043  \\
            0.1640624999999999  0.7317770956661673  \\
            0.171875  0.7610551836676646  \\
            0.1796875000000001  0.7871419150359282  \\
            0.1874999999999999  0.8155648816047909  \\
            0.1953125  0.859944631339821  \\
            0.2031250000000001  0.8956054326077852  \\
            0.2109374999999999  0.9280748772425158  \\
            0.21875  0.9628805570778456  \\
            0.2265625000000001  0.9930137665119773  \\
            0.2343749999999999  1.025483211146708  \\
            0.2421875  1.0588077772140732  \\
            0.2500000000000001  1.0915063509461096  \\
            0.2578124999999999  1.0635691486467085  \\
            0.265625  1.068171505014972  \\
            0.2734375000000001  1.0846841664835352  \\
            0.2812499999999999  1.0924778794850325  \\
            0.2890625  1.1177094894206616  \\
            0.2968750000000001  1.1342221508892245  \\
            0.3046874999999999  1.147543455724554  \\
            0.3125  1.163200995760483  \\
            0.3203125000000001  1.194815318962579  \\
            0.3281249999999999  1.217710693697609  \\
            0.3359375  1.2374147117994052  \\
            0.3437500000000001  1.2594549651018012  \\
            0.3515624999999999  1.2768227480029988  \\
            0.359375  1.2965267661047954  \\
            0.3671875000000001  1.3170859056392261  \\
            0.3749999999999999  1.3370190528383288  \\
            0.3828125  1.3733058464416226  \\
            0.3906250000000001  1.4008736915778504  \\
            0.3984374999999999  1.4252501800808448  \\
            0.40625  1.4519629037844384  \\
            0.4140625000000001  1.474003157086834  \\
            0.4218749999999999  1.4983796455898277  \\
            0.4296875  1.523611255525457  \\
            0.4375000000000001  1.5482168731257568  \\
            0.4453124999999999  1.5685468835628826  \\
            0.453125  1.5912131292006078  \\
            0.4609375000000001  1.6147344962709673  \\
            0.4687499999999999  1.6376298710059978  \\
            0.4765625  1.6617772304116871  \\
            0.4843750000000001  1.685298597482047  \\
            0.4921874999999999  1.7085908354551012  \\
            0.5  1.7320508075688772  \\
            0.5078124999999999  1.2886689604551023  \\
            0.5156250000000001  1.0954548474820478  \\
            0.5234375  0.9938084804116879  \\
            0.5312499999999999  0.8251298710059987  \\
            0.5390625000000001  0.7905157462709684  \\
            0.546875  0.6888693792006084  \\
            0.5546874999999999  0.5626875085628835  \\
            0.5625000000000001  0.4544668731257571  \\
            0.5703125  0.4689237555254577  \\
            0.5781249999999999  0.4163483955898287  \\
            0.5859375000000001  0.33923753208683416  \\
            0.59375  0.28008790378443865  \\
            0.6015624999999999  0.18501580508084536  \\
            0.6093750000000001  0.10790494157785113  \\
            0.6171875  0.03736834644162307  \\
            0.6249999999999999  -0.03798094716167105  \\
            0.6328125000000001  0.01239840563922745  \\
            0.640625  -0.004254483895203887  \\
            0.6484374999999999  -0.0454428769970005  \\
            0.6562500000000001  -0.06867003489819834  \\
            0.6640625  -0.12781966320059385  \\
            0.6718749999999999  -0.16900805630239052  \\
            0.6796875000000001  -0.2036221810374208  \\
            0.6875  -0.24304900423951692  \\
            0.6953124999999999  -0.3153471692754453  \\
            0.7031250000000001  -0.369684099110775  \\
            0.7109375  -0.41744676057933816  \\
            0.7187499999999999  -0.4700221205149671  \\
            0.7265625000000001  -0.5129720835164643  \\
            0.734375  -0.5607347449850273  \\
            0.7421874999999999  -0.610258976353291  \\
            0.7500000000000001  -0.6584936490538904  \\
            0.7578125  -0.5818172227859261  \\
            0.7656249999999999  -0.5721730388532913  \\
            0.7734375000000001  -0.5870643584880221  \\
            0.78125  -0.583994442922154  \\
            0.7890624999999999  -0.6168469977574836  \\
            0.7968750000000001  -0.6317383173922142  \\
            0.8046875  -0.6400553686601784  \\
            0.8124999999999999  -0.653185118395209  \\
            0.8203125000000001  -0.6991862099640713  \\
            0.828125  -0.7272260663323351  \\
            0.8359374999999999  -0.748691654333832  \\
            0.8437500000000001  -0.7749699408023955  \\
            0.8515625  -0.7916228303368265  \\
            0.8593749999999999  -0.8130884183383238  \\
            0.8671875000000001  -0.8363155762395214  \\
            0.875  -0.8582531754730547  \\
            0.8828124999999999  -0.9138796639760493  \\
            0.8906250000000001  -0.9515449172784447  \\
            0.8984375  -0.9826359022140736  \\
            0.9062499999999999  -1.018539585616769  \\
            0.9140625000000001  -1.0448178720853318  \\
            0.921875  -1.0759088570209605  \\
            0.9296874999999999  -1.1087614118562903  \\
            0.9375000000000001  -1.1403244080239554  \\
            0.9453125  -1.1630795546931172  \\
            0.9531249999999999  -1.1906473998293452  \\
            0.9609375000000001  -1.2199768148652734  \\
            0.96875  -1.2480166712335372  \\
            0.9765624999999999  -1.2786356449371306  \\
            0.9843750000000001  -1.307965059973059  \\
            0.9921875  -1.3368224637769512  \\
            0.9999999999999999  -1.3660254037844388  \\
            1.0078125  -1.041667100141473  \\
            1.015625  -0.9067568875374306  \\
            1.0234375  -0.841189488952101  \\
            1.03125  -0.7248596273652741  \\
            1.0390625  -0.7100546917814418  \\
            1.046875  -0.644487293196112  \\
            1.0546875  -0.5603395435935669  \\
            1.0625  -0.4897935549580875  \\
            1.0703125  -0.5121493214086864  \\
            1.078125  -0.48374262485778813  \\
            1.0859375  -0.43675557728967396  \\
            1.09375  -0.403370290688626  \\
            1.1015625  -0.34278148215344606  \\
            1.109375  -0.29579443458533217  \\
            1.1171875  -0.2537859770673679  \\
            1.125  -0.20813293868263716  \\
            1.1328125  -0.257692227067368  \\
            1.140625  -0.25648905245060144  \\
            1.1484375  -0.23670552681661933  \\
            1.15625  -0.2305237621497031  \\
            1.1640625  -0.19713847554865507  \\
            1.171875  -0.17735494991467282  \\
            1.1796875  -0.16255001433084046  \\
            1.1875  -0.14410049788024154  \\
            1.1953125  -0.10075803117889412  \\
            1.203125  -0.07101732544461248  \\
            1.2109375  -0.046255209760480566  \\
            1.21875  -0.017848513209582348  \\
            1.2265625  0.003269021607783029  \\
            1.234375  0.02803113729191481  \\
            1.2421875  0.054127262159429895  \\
            1.25  0.07924682452694516  \\
            1.2578125  0.009773175941615468  \\
            1.265625  -0.008938009642216992  \\
            1.2734375  -0.00906884420883369  \\
            1.28125  -0.02280143974251637  \\
            1.2890625  -0.009330513342067187  \\
            1.296875  -0.00946134790868398  \\
            1.3046875  -0.014570772525450448  \\
            1.3125  -0.016035616275450415  \\
            1.3203125  0.007392490225298179  \\
            1.328125  0.017218835758980845  \\
            1.3359375  0.022066591242513835  \\
            1.34375  0.03055892759281328  \\
            1.3515625  0.03176210220957975  \\
            1.359375  0.03660985769311273  \\
            1.3671875  0.04279162236002894  \\
            1.375  0.04799682452694514  \\
            1.3828125  0.0787140927612267  \\
            1.390625  0.09582960002844236  \\
            1.3984375  0.10796651724550826  \\
            1.40625  0.12374801532934072  \\
            1.4140625  0.13224035167964013  \\
            1.421875  0.14437726889670607  \\
            1.4296875  0.15784819529715519  \\
            1.4375  0.17034255919760438  \\
            1.4453125  0.17616687718113733  \\
            1.453125  0.1856357760314368  \\
            1.4609375  0.1964386840651194  \\
            1.46875  0.2062650295988021  \\
            1.4765625  0.21804450013248483  \\
            1.484375  0.22884740816616744  \\
            1.4921875  0.23929286951646686  \\
            1.5  0.25000000000000006  \\
            1.5078125  0.1659480279476044  \\
            1.515625  0.1326585188967061  \\
            1.5234375  0.11794936086302339  \\
            1.53125  0.08963844186227479  \\
            1.5390625  0.08853104479565804  \\
            1.546875  0.07382188676197535  \\
            1.5546875  0.05413413867814296  \\
            1.5625  0.038090971461077025  \\
            1.5703125  0.04694075449475971  \\
            1.578125  0.04218877656137648  \\
            1.5859375  0.032458208577843504  \\
            1.59375  0.026372221461077018  \\
            1.6015625  0.012997072610777581  \\
            1.609375  0.0032665046272446277  \\
            1.6171875  -0.005130054172905084  \\
            1.625  -0.014503175473054835  \\
            1.6328125  0.0016357692941608361  \\
            1.640625  0.00417295309431055  \\
            1.6484375  0.0017315468443105525  \\
            1.65625  0.002934721461077027  \\
            1.6640625  -0.003151265655689447  \\
            1.671875  -0.00559267190568945  \\
            1.6796875  -0.0067000689723062095  \\
            1.6875  -0.008784028538922972  \\
            1.6953125  -0.017538034022455927  \\
            1.703125  -0.022647458639222417  \\
            1.7109375  -0.02642287407260565  \\
            1.71875  -0.031174852005988893  \\
            1.7265625  -0.03397370493937212  \\
            1.734375  -0.037749120372755345  \\
            1.7421875  -0.04188198248952183  \\
            1.75  -0.04575317547305484  \\
            1.7578125  -0.024278193972306208  \\
            1.765625  -0.016404973438623523  \\
            1.7734375  -0.013510342955090562  \\
            1.78125  -0.006971131604791121  \\
            1.7890625  -0.007721081988024642  \\
            1.796875  -0.0048264515044916766  \\
            1.8046875  -0.0005978118375754758  \\
            1.8125  0.002654265329340731  \\
            1.8203125  -0.0007637034206592706  \\
            1.828125  -0.0005370913038927881  \\
            1.8359375  0.0010235299962569333  \\
            1.84375  0.0016075887964066562  \\
            1.8515625  0.004144772596556379  \\
            1.859375  0.0057053938967061  \\
            1.8671875  0.006908568513472576  \\
            1.875  0.008373412263472584  \\
            1.8828125  0.003002318513472583  \\
            1.890625  0.0012758056302390654  \\
            1.8984375  0.0008833019303887883  \\
            1.90625  -0.0004857642694614884  \\
            1.9140625  9.829453068823418e-5  \\
            1.921875  -0.00029420916916204296  \\
            1.9296875  -0.0010441595523955607  \\
            1.9375  -0.0015324408023955614  \\
            1.9453125  -0.00023348863547935708  \\
            1.953125  8.890103143684708e-5  \\
            1.9609375  5.384401496981049e-5  \\
            1.96875  0.0002804561317362921  \\
            1.9765625  -1.6270017964262505e-5  \\
            1.984375  -5.132703443129902e-5  \\
            1.9921875  9.393499251387002e-6  \\
            2.0  0.0  \\
        }
        ;
\end{axis}
\end{tikzpicture}

%% file: img/discretedb24.tikz
\begin{tikzpicture}
\begin{groupplot}[width={.5\textwidth}, height={.25\textwidth}, group style={group size={1 by 2}}]
    \nextgroupplot[]
    \addplot+[samples at={1,2,3,4,5,6,7,8,9,10,11}, ycomb, mark={*}]
        table[row sep={\\}]
        {
            \\
            1.0  0.6372595264191646  \\
            2.0  0.9330127018922194  \\
            3.0  1.103765877365274  \\
            4.0  1.3660254037844388  \\
            5.0  0.34150635094610976  \\
            6.0  0.0  \\
            7.0  -0.09150635094610964  \\
            8.0  -0.36602540378443876  \\
            9.0  0.021234122634725853  \\
            10.0  0.06698729810778069  \\
            11.0  -0.012259526419164493  \\
        }
        ;
    \nextgroupplot[]
    \addplot+[samples at={1,2,3,4,5,6,7,8,9,10,11}, ycomb, mark={*}]
        table[row sep={\\}]
        {
            \\
            1.0  0.0  \\
            2.0  0.0  \\
            3.0  0.0  \\
            4.0  0.0  \\
            5.0  2.732050807568877  \\
            6.0  -0.9999999999999996  \\
            7.0  -0.7320508075688783  \\
            8.0  0.0  \\
            9.0  0.0  \\
            10.0  0.0  \\
            11.0  0.0  \\
        }
        ;
\end{groupplot}
\end{tikzpicture}

%% file: img/discretedb34.tikz
\begin{tikzpicture}
\begin{groupplot}[width={.5\textwidth}, height={.25\textwidth}, group style={group size={1 by 2}}]
    \nextgroupplot[]
    \addplot+[samples at={1,2,3,4,5,6,7,8,9,10,11,12,13,14,15,16,17,18,19}, ycomb, mark={*}]
        table[row sep={\\}]
        {
            \\
            1.0  0.2847166242333907  \\
            2.0  0.605178468387557  \\
            3.0  0.8981130495184216  \\
            4.0  1.2863350694256988  \\
            5.0  0.8899160481330288  \\
            6.0  0.4411224814623552  \\
            7.0  0.1394188819005187  \\
            8.0  -0.38583696104587756  \\
            9.0  -0.2029799345209342  \\
            10.0  -0.014970591386629456  \\
            11.0  -0.040567571439028935  \\
            12.0  0.09526754600378076  \\
            13.0  0.02994405989833984  \\
            14.0  -0.031541302974913414  \\
            15.0  0.0030251312919360785  \\
            16.0  0.0042343456163981144  \\
            17.0  -0.0015967977438250268  \\
            18.0  0.00021094451163078802  \\
            19.0  1.0508728152663858e-5  \\
        }
        ;
    \nextgroupplot[]
    \addplot+[samples at={1,2,3,4,5,6,7,8,9,10,11,12,13,14,15,16,17,18,19}, ycomb, mark={*}]
        table[row sep={\\}]
        {
            \\
            1.0  0.0  \\
            2.0  0.0  \\
            3.0  0.0  \\
            4.0  0.0  \\
            5.0  0.0  \\
            6.0  0.0  \\
            7.0  0.0  \\
            8.0  1.574884046305118  \\
            9.0  5.144659499172568  \\
            10.0  8.86142085846009  \\
            11.0  -30.72596107090182  \\
            12.0  16.144996666964037  \\
            13.0  0.0  \\
            14.0  0.0  \\
            15.0  0.0  \\
            16.0  0.0  \\
            17.0  0.0  \\
            18.0  0.0  \\
            19.0  0.0  \\
        }
        ;
\end{groupplot}
\end{tikzpicture}

%% file: img/discretedb44.tikz
\begin{tikzpicture}
\begin{groupplot}[width={.5\textwidth}, height={.25\textwidth}, group style={group size={1 by 2}}]
    \nextgroupplot[]
    \addplot+[samples at={1,2,3,4,5,6,7,8,9,10,11,12,13,14,15,16,17,18,19,20,21,22,23,24,25,26}, ycomb, mark={*}]
        table[row sep={\\}]
        {
            \\
            1.0  0.1069089516202416  \\
            2.0  0.3281394313733513  \\
            3.0  0.6175569012001709  \\
            4.0  1.0071699777256031  \\
            5.0  1.100816874468014  \\
            6.0  0.8772950941928345  \\
            7.0  0.5362551559984561  \\
            8.0  -0.03383695405283656  \\
            9.0  -0.3020570317036663  \\
            10.0  -0.2420075454221948  \\
            11.0  -0.1753436683264306  \\
            12.0  0.03961046271590372  \\
            13.0  0.11822215797083187  \\
            14.0  0.034299729174901344  \\
            15.0  0.016274438999383706  \\
            16.0  -0.011764358205726692  \\
            17.0  -0.023480523154080942  \\
            18.0  0.002152148639222631  \\
            19.0  0.005283787771867855  \\
            20.0  -0.0011979575961769996  \\
            21.0  -0.00040859629443929516  \\
            22.0  0.00012142423810055572  \\
            23.0  -2.6619872719321462e-5  \\
            24.0  1.8829413233544343e-5  \\
            25.0  -1.8329069010770086e-6  \\
            26.0  -2.821962152280963e-7  \\
        }
        ;
    \nextgroupplot[]
    \addplot+[samples at={1,2,3,4,5,6,7,8,9,10,11,12,13,14,15,16,17,18,19,20,21,22,23,24,25,26}, ycomb, mark={*}]
        table[row sep={\\}]
        {
            \\
            1.0  0.0  \\
            2.0  0.0  \\
            3.0  0.0  \\
            4.0  0.0  \\
            5.0  0.0  \\
            6.0  0.0  \\
            7.0  0.0  \\
            8.0  0.0  \\
            9.0  0.05892863871614476  \\
            10.0  -0.3827503798636192  \\
            11.0  3.992495182812906  \\
            12.0  9.604016140213815  \\
            13.0  5.921362042026277  \\
            14.0  -24.267607600461194  \\
            15.0  27.693689557927655  \\
            16.0  -67.11225684737529  \\
            17.0  -8.629287491157518  \\
            18.0  112.2962708060729  \\
            19.0  -58.17486016093983  \\
            20.0  0.0  \\
            21.0  0.0  \\
            22.0  0.0  \\
            23.0  0.0  \\
            24.0  0.0  \\
            25.0  0.0  \\
            26.0  0.0  \\
        }
        ;
\end{groupplot}
\end{tikzpicture}

%% file: img/discretecdf314.tikz
\begin{tikzpicture}
\begin{groupplot}[width={.5\textwidth}, height={.25\textwidth}, group style={group size={1 by 2}}]
    \nextgroupplot[]
    \addplot+[samples at={-5,-4,-3,-2,-1,0,1,2,3,4,5}, ycomb, mark={*}]
        table[row sep={\\}]
        {
            \\
            -5.0  0.03125  \\
            -4.0  0.125  \\
            -3.0  0.28125  \\
            -2.0  0.5  \\
            -1.0  0.6875  \\
            0.0  0.75  \\
            1.0  0.6875  \\
            2.0  0.5  \\
            3.0  0.28125  \\
            4.0  0.125  \\
            5.0  0.03125  \\
        }
        ;
    \nextgroupplot[]
    \addplot+[samples at={-5,-4,-3,-2,-1,0,1,2,3,4,5}, ycomb, mark={*}]
        table[row sep={\\}]
        {
            \\
            -5.0  0.0  \\
            -4.0  0.0  \\
            -3.0  0.0  \\
            -2.0  0.0  \\
            -1.0  -2.000000000000001  \\
            0.0  5.000000000000002  \\
            1.0  -2.000000000000002  \\
            2.0  0.0  \\
            3.0  0.0  \\
            4.0  0.0  \\
            5.0  0.0  \\
        }
        ;
\end{groupplot}
\end{tikzpicture}

%% file: img/discretecdf424.tikz
\begin{tikzpicture}
\begin{groupplot}[width={.5\textwidth}, height={.25\textwidth}, group style={group size={1 by 2}}]
    \nextgroupplot[]
    \addplot+[samples at={-7,-6,-5,-4,-3,-2,-1,0,1,2,3,4,5,6,7}, ycomb, mark={*}]
        table[row sep={\\}]
        {
            \\
            -7.0  0.0026041666666666665  \\
            -6.0  0.020833333333333332  \\
            -5.0  0.0703125  \\
            -4.0  0.16666666666666663  \\
            -3.0  0.31510416666666663  \\
            -2.0  0.47916666666666663  \\
            -1.0  0.6119791666666666  \\
            0.0  0.666666666666667  \\
            1.0  0.611979166666667  \\
            2.0  0.47916666666666696  \\
            3.0  0.31510416666666696  \\
            4.0  0.16666666666666607  \\
            5.0  0.07031249999999845  \\
            6.0  0.02083333333333326  \\
            7.0  0.0026041666666660745  \\
        }
        ;
    \nextgroupplot[]
    \addplot+[samples at={-7,-6,-5,-4,-3,-2,-1,0,1,2,3,4,5,6,7}, ycomb, mark={*}]
        table[row sep={\\}]
        {
            \\
            -7.0  0.0  \\
            -6.0  0.0  \\
            -5.0  0.0  \\
            -4.0  0.0  \\
            -3.0  0.0  \\
            -2.0  0.6666666666664526  \\
            -1.0  -5.333333333332785  \\
            0.0  10.333333333332845  \\
            1.0  -5.333333333333189  \\
            2.0  0.6666666666666703  \\
            3.0  0.0  \\
            4.0  0.0  \\
            5.0  0.0  \\
            6.0  0.0  \\
            7.0  0.0  \\
        }
        ;
\end{groupplot}
\end{tikzpicture}

%% file: img/discretecdf514.tikz
\begin{tikzpicture}
\begin{groupplot}[width={.5\textwidth}, height={.25\textwidth}, group style={group size={1 by 2}}]
    \nextgroupplot[]
    \addplot+[samples at={-9,-8,-7,-6,-5,-4,-3,-2,-1,0,1,2,3,4,5,6,7,8,9}, ycomb, mark={*}]
        table[row sep={\\}]
        {
            \\
            -9.0  0.00016276041666666666  \\
            -8.0  0.0026041666666666665  \\
            -7.0  0.01318359375  \\
            -6.0  0.04166666666666677  \\
            -5.0  0.10091145833333337  \\
            -4.0  0.19791666666666688  \\
            -3.0  0.32486979166666713  \\
            -2.0  0.45833333333333304  \\
            -1.0  0.560872395833333  \\
            0.0  0.598958333333333  \\
            1.0  0.560872395833333  \\
            2.0  0.45833333333333215  \\
            3.0  0.32486979166666163  \\
            4.0  0.19791666666667274  \\
            5.0  0.10091145833333215  \\
            6.0  0.04166666666666785  \\
            7.0  0.013183593750006217  \\
            8.0  0.00260416666667318  \\
            9.0  0.0001627604166678509  \\
        }
        ;
    \nextgroupplot[]
    \addplot+[samples at={-9,-8,-7,-6,-5,-4,-3,-2,-1,0,1,2,3,4,5,6,7,8,9}, ycomb, mark={*}]
        table[row sep={\\}]
        {
            \\
            -9.0  0.0  \\
            -8.0  0.0  \\
            -7.0  0.0  \\
            -6.0  0.0  \\
            -5.0  0.0  \\
            -4.0  0.0  \\
            -3.0  0.0  \\
            -2.0  6.277777777774357  \\
            -1.0  -28.444444444440517  \\
            0.0  45.33333333333291  \\
            1.0  -28.444444444445338  \\
            2.0  6.277777777778058  \\
            3.0  0.0  \\
            4.0  0.0  \\
            5.0  0.0  \\
            6.0  0.0  \\
            7.0  0.0  \\
            8.0  0.0  \\
            9.0  0.0  \\
        }
        ;
\end{groupplot}
\end{tikzpicture}

%% file: dwttabular1013.tikz
\begin{tabular}{|l|l|l|l|} 
 \hline
$p$&$\tilde p$&$ \|W_{10}\|$&$ \|W^{-1}_{10}\|$ \\\hline
1&1&1.00e+00&1.00e+00\\\hline 
&3&1.30e+00&1.37e+00\\\hline 
&5&1.40e+00&1.54e+00\\\hline 
2&2&2.33e+00&1.41e+00\\\hline 
&4&1.76e+00&1.41e+00\\\hline 
&6&1.74e+00&1.41e+00\\\hline 
3&1&5.43e+01&2.00e+00\\\hline 
&3&5.05e+00&2.00e+00\\\hline 
&5&2.97e+00&2.00e+00\\\hline 
\end{tabular} 

%% file: dwttabular1046.tikz
\begin{tabular}{|l|l|l|l|} 
 \hline
$p$&$\tilde p$&$ \|W_{10}\|$&$ \|W^{-1}_{10}\|$ \\\hline
4&2&2.94e+02&2.83e+00\\\hline 
&4&1.77e+01&2.83e+00\\\hline 
&6&5.51e+00&2.83e+00\\\hline 
5&1&9.59e+04&4.27e+00\\\hline 
&3&2.85e+03&4.00e+00\\\hline 
&5&1.46e+02&4.00e+00\\\hline 
6&2&1.01e+06&5.66e+00\\\hline 
&4&3.91e+04&5.66e+00\\\hline 
&6&2.00e+03&5.66e+00\\\hline 
\end{tabular} 

%% file: img/AZtimings1d-3d.tikz
\begin{tikzpicture}
\begin{groupplot}[legend style={fill opacity={1.0}, text opacity={1}}, legend pos={north west}, xlabel={$N$}, xmode={log}, ymode={log}, legend cell align={left}, cycle list name={mark list*}, ymin={0.001}, group style={group size={3 by 1}}]
    \nextgroupplot[]
    \addplot
        table[row sep={\\}]
        {
            \\
            128.0  0.009899876  \\
            256.0  0.014340951  \\
            1024.0  0.020245233  \\
            2048.0  0.024358307  \\
            8192.0  0.041900199  \\
            16384.0  0.077871171  \\
            65536.0  0.126022159  \\
            131072.0  0.449578517  \\
            262144.0  1.673088098  \\
            1.048576e6  8.333767936  \\
        }
        ;
    \addplot
        table[row sep={\\}]
        {
            \\
            128.0  0.017142486  \\
            256.0  0.020231018  \\
            1024.0  0.026454036  \\
            2048.0  0.03096518  \\
            8192.0  0.050358526  \\
            16384.0  0.108172684  \\
            65536.0  0.170933755  \\
            131072.0  0.945229296  \\
            262144.0  2.430539623  \\
            1.048576e6  10.973877053  \\
        }
        ;
    \addplot
        table[row sep={\\}]
        {
            \\
            128.0  0.021897963  \\
            256.0  0.026700621  \\
            1024.0  0.038818468  \\
            2048.0  0.039193275  \\
            8192.0  0.062070302  \\
            16384.0  0.131094383  \\
            65536.0  0.23088579  \\
            131072.0  0.75621866  \\
            262144.0  3.561206865  \\
            1.048576e6  14.173685418  \\
        }
        ;
    \addplot
        table[row sep={\\}]
        {
            \\
            128.0  0.004602924  \\
            256.0  0.005238237  \\
            1024.0  0.005492982  \\
            2048.0  0.010109468  \\
            8192.0  0.021395398  \\
            16384.0  0.05945554  \\
            65536.0  0.109249105  \\
            131072.0  0.4262257  \\
            262144.0  2.187695273  \\
            1.048576e6  8.989004114  \\
        }
        ;
    \addplot
        table[row sep={\\}]
        {
            \\
            128.0  0.004211221  \\
            256.0  0.005893489  \\
            1024.0  0.005718999  \\
            2048.0  0.0099338  \\
            8192.0  0.021205718  \\
            16384.0  0.067279071  \\
            65536.0  0.137144721  \\
            131072.0  0.481626918  \\
            262144.0  1.913128743  \\
            1.048576e6  10.437178284  \\
        }
        ;
    \addplot
        table[row sep={\\}]
        {
            \\
            128.0  0.004613329  \\
            256.0  0.007711363  \\
            1024.0  0.00666258  \\
            2048.0  0.009419637  \\
            8192.0  0.02282761  \\
            16384.0  0.06986202  \\
            65536.0  0.129459413  \\
            131072.0  0.550871963  \\
            262144.0  2.18945405  \\
            1.048576e6  10.284367146  \\
        }
        ;
    \addplot
        table[row sep={\\}]
        {
            \\
            128.0  0.008390612  \\
            256.0  0.010172008  \\
            1024.0  0.009335014  \\
            2048.0  0.013630088  \\
            8192.0  0.023891777  \\
            16384.0  0.072475184  \\
            65536.0  0.131459053  \\
            131072.0  0.52763786  \\
            262144.0  2.183443629  \\
            1.048576e6  10.204141149  \\
        }
        ;
    \addplot
        table[row sep={\\}]
        {
            \\
            128.0  0.008686331  \\
            256.0  0.008958512  \\
            1024.0  0.009639375  \\
            2048.0  0.013123783  \\
            8192.0  0.032531769  \\
            16384.0  0.075033653  \\
            65536.0  0.129769062  \\
            131072.0  0.533436891  \\
            262144.0  2.243546931  \\
            1.048576e6  10.575803695  \\
        }
        ;
    \addplot
        table[row sep={\\}]
        {
            \\
            128.0  0.0085393  \\
            256.0  0.009585086  \\
            1024.0  0.010164779  \\
            2048.0  0.013699522  \\
            8192.0  0.029041157  \\
            16384.0  0.079779141  \\
            65536.0  0.148838902  \\
            131072.0  0.611762917  \\
            262144.0  2.888624073  \\
            1.048576e6  12.759567981  \\
        }
        ;
    \addplot[style={black,dashed}]
        table[row sep={\\}]
        {
            \\
            128.0  0.00128  \\
            256.0  0.00256  \\
            1024.0  0.01024  \\
            2048.0  0.02048  \\
            8192.0  0.08192  \\
            16384.0  0.16384  \\
            65536.0  0.65536  \\
            131072.0  1.31072  \\
            262144.0  2.62144  \\
            1.048576e6  10.48576  \\
        }
        ;
    \legend{{$\texttt{db2}$},{$\texttt{db3}$},{$\texttt{db4}$},{$\texttt{cdf31}$},{$\texttt{cdf33}$},{$\texttt{cdf35}$},{$\texttt{cdf42}$},{$\texttt{cdf44}$},{$\texttt{cdf46}$},{$\mathcal O(N^{1.00})$}}
    \nextgroupplot[]
    \addplot
        table[row sep={\\}]
        {
            \\
            64.0  0.019138978  \\
            256.0  0.098626629  \\
            1024.0  0.454422121  \\
            4096.0  3.123187186  \\
            16384.0  26.484014227  \\
            65536.0  255.193790657  \\
            262144.0  0.0  \\
            1.048576e6  0.0  \\
        }
        ;
    \addplot
        table[row sep={\\}]
        {
            \\
            64.0  0.030564859  \\
            256.0  0.240590857  \\
            1024.0  1.411568874  \\
            4096.0  9.581889099  \\
            16384.0  89.429581742  \\
            65536.0  923.462167195  \\
            262144.0  0.0  \\
            1.048576e6  0.0  \\
        }
        ;
    \addplot
        table[row sep={\\}]
        {
            \\
            64.0  0.043453888  \\
            256.0  0.33926656  \\
            1024.0  1.978699542  \\
            4096.0  13.137903145  \\
            16384.0  117.279186228  \\
            65536.0  1194.663278371  \\
            262144.0  0.0  \\
            1.048576e6  0.0  \\
        }
        ;
    \addplot
        table[row sep={\\}]
        {
            \\
            64.0  0.016160767  \\
            256.0  0.079943752  \\
            1024.0  0.449783366  \\
            4096.0  3.192330432  \\
            16384.0  25.410880038  \\
            65536.0  240.618940611  \\
            262144.0  0.0  \\
            1.048576e6  0.0  \\
        }
        ;
    \addplot
        table[row sep={\\}]
        {
            \\
            64.0  0.018713235  \\
            256.0  0.108382887  \\
            1024.0  0.580150583  \\
            4096.0  3.962469342  \\
            16384.0  29.919168844  \\
            65536.0  268.079083205  \\
            262144.0  0.0  \\
            1.048576e6  0.0  \\
        }
        ;
    \addplot
        table[row sep={\\}]
        {
            \\
            64.0  0.025871177  \\
            256.0  0.129458084  \\
            1024.0  0.755191669  \\
            4096.0  4.716746636  \\
            16384.0  35.295167234  \\
            65536.0  302.152513666  \\
            262144.0  0.0  \\
            1.048576e6  0.0  \\
        }
        ;
    \addplot
        table[row sep={\\}]
        {
            \\
            64.0  0.022967965  \\
            256.0  0.113595252  \\
            1024.0  0.661453073  \\
            4096.0  4.587808727  \\
            16384.0  37.995362847  \\
            65536.0  362.682190331  \\
            262144.0  0.0  \\
            1.048576e6  0.0  \\
        }
        ;
    \addplot
        table[row sep={\\}]
        {
            \\
            64.0  0.028516398  \\
            256.0  0.144397503  \\
            1024.0  0.836187631  \\
            4096.0  5.594721182  \\
            16384.0  43.725694077  \\
            65536.0  400.497684764  \\
            262144.0  0.0  \\
            1.048576e6  0.0  \\
        }
        ;
    \addplot
        table[row sep={\\}]
        {
            \\
            64.0  0.035963478  \\
            256.0  0.181783422  \\
            1024.0  1.054088342  \\
            4096.0  6.505950088  \\
            16384.0  50.38282855  \\
            65536.0  440.102833138  \\
            262144.0  0.0  \\
            1.048576e6  0.0  \\
        }
        ;
    \addplot[style={black,dashed}]
        table[row sep={\\}]
        {
            \\
            64.0  0.00012288  \\
            256.0  0.00196608  \\
            1024.0  0.03145728  \\
            4096.0  0.50331648  \\
            16384.0  8.05306368  \\
            65536.0  128.84901888  \\
        }
        ;
    \addplot[style={black,dotted,thick}]
        table[row sep={\\}]
        {
            \\
            64.0  0.003501322839482538  \\
            256.0  0.028010582715860304  \\
            1024.0  0.22408466172688243  \\
            4096.0  1.7926772938150595  \\
            16384.0  14.341418350520476  \\
            65536.0  114.7313468041638  \\
        }
        ;
    \legend{{},{},{},{},{},{},{},{},{},{$\mathcal O(N^{2.00})$},{$\mathcal O(N^{1.50})$}}
    \nextgroupplot[]
    \addplot
        table[row sep={\\}]
        {
            \\
            64.0  0.043804572  \\
            512.0  1.603994372  \\
            4096.0  40.521657436  \\
            32768.0  1567.190559478  \\
        }
        ;
    \addplot
        table[row sep={\\}]
        {
            \\
            64.0  0.0664589  \\
            512.0  2.809119305  \\
            4096.0  130.509693414  \\
            32768.0  6518.404201037  \\
        }
        ;
    \addplot
        table[row sep={\\}]
        {
            \\
            64.0  0.079444104  \\
            512.0  3.999944357  \\
            4096.0  175.417624062  \\
            32768.0  3393.167662544  \\
        }
        ;
    \addplot
        table[row sep={\\}]
        {
            \\
            64.0  0.040745216  \\
            512.0  1.509244403  \\
            4096.0  38.549607027  \\
            32768.0  1252.696788967  \\
        }
        ;
    \addplot
        table[row sep={\\}]
        {
            \\
            64.0  0.052778126  \\
            512.0  2.104878631  \\
            4096.0  51.968447306  \\
            32768.0  1480.905336022  \\
        }
        ;
    \addplot
        table[row sep={\\}]
        {
            \\
            64.0  0.064147449  \\
            512.0  3.084559341  \\
            4096.0  69.474925273  \\
            32768.0  1817.673329234  \\
        }
        ;
    \addplot
        table[row sep={\\}]
        {
            \\
            64.0  0.057773266  \\
            512.0  2.296781883  \\
            4096.0  106.453014747  \\
            32768.0  3931.600476457  \\
        }
        ;
    \addplot
        table[row sep={\\}]
        {
            \\
            64.0  0.070114044  \\
            512.0  3.27999242  \\
            4096.0  139.123403616  \\
            32768.0  4436.922176358  \\
        }
        ;
    \addplot
        table[row sep={\\}]
        {
            \\
            64.0  0.082277131  \\
            512.0  4.261419345  \\
            4096.0  176.5255296  \\
            32768.0  5725.115170135  \\
        }
        ;
    \addplot[style={black,dashed}]
        table[row sep={\\}]
        {
            \\
            64.0  0.001638400000000001  \\
            512.0  0.20971520000000018  \\
            4096.0  26.843545600000034  \\
            32768.0  3435.9738368000053  \\
        }
        ;
    \addplot[style={black,dotted,thick}]
        table[row sep={\\}]
        {
            \\
            64.0  0.012191732034540988  \\
            512.0  0.7849271135374057  \\
            4096.0  50.53511443826282  \\
            32768.0  3253.547682636673  \\
        }
        ;
    \legend{{},{},{},{},{},{},{},{},{},{$\mathcal O(N^{2.33})$},{$\mathcal O(N^{2.00})$}}
\end{groupplot}
\end{tikzpicture}

%% file: img/AZR1timings1d-3d.tikz
\begin{tikzpicture}
\begin{groupplot}[legend style={fill opacity={0.6}, text opacity={1}}, legend pos={north west}, xlabel={$N$}, xmode={log}, ymode={log}, legend cell align={left}, cycle list name={mark list*}, group style={group size={3 by 1}}]
    \nextgroupplot[]
    \addplot
        table[row sep={\\}]
        {
            \\
            128.0  0.028785607  \\
            256.0  0.023972726  \\
            1024.0  0.029454382  \\
            2048.0  0.032505482  \\
            8192.0  0.048160813  \\
            16384.0  0.090316811  \\
            65536.0  0.149430517  \\
            131072.0  0.530575911  \\
            262144.0  2.293440111  \\
            1.048576e6  10.421557485  \\
        }
        ;
    \addplot
        table[row sep={\\}]
        {
            \\
            128.0  0.051097208  \\
            256.0  0.049864286  \\
            1024.0  0.068847095  \\
            2048.0  0.06134338  \\
            8192.0  0.086393517  \\
            16384.0  0.178710346  \\
            65536.0  0.331446875  \\
            131072.0  1.225222704  \\
            262144.0  5.456210368  \\
            1.048576e6  23.983957382  \\
        }
        ;
    \addplot
        table[row sep={\\}]
        {
            \\
            128.0  0.073080108  \\
            256.0  0.079240175  \\
            1024.0  0.088666025  \\
            2048.0  0.089902834  \\
            8192.0  0.128256426  \\
            16384.0  0.288614  \\
            65536.0  0.51489644  \\
            131072.0  2.33687392  \\
            262144.0  9.539227126  \\
            1.048576e6  44.961696642  \\
        }
        ;
    \addplot
        table[row sep={\\}]
        {
            \\
            128.0  0.01555325  \\
            256.0  0.016112541  \\
            1024.0  0.015734112  \\
            2048.0  0.019537592  \\
            8192.0  0.026858759  \\
            16384.0  0.076360505  \\
            65536.0  0.133971463  \\
            131072.0  0.598336059  \\
            262144.0  2.321888673  \\
            1.048576e6  10.935618125  \\
        }
        ;
    \addplot
        table[row sep={\\}]
        {
            \\
            128.0  0.030097173  \\
            256.0  0.015722057  \\
            1024.0  0.016280436  \\
            2048.0  0.022758175  \\
            8192.0  0.03880451  \\
            16384.0  0.123903554  \\
            65536.0  0.236157609  \\
            131072.0  1.007250673  \\
            262144.0  4.642003714  \\
            1.048576e6  23.190003313  \\
        }
        ;
    \addplot
        table[row sep={\\}]
        {
            \\
            128.0  0.028333771  \\
            256.0  0.02047153  \\
            1024.0  0.018561629  \\
            2048.0  0.025029723  \\
            8192.0  0.052016553  \\
            16384.0  0.195450459  \\
            65536.0  0.360603878  \\
            131072.0  1.601531893  \\
            262144.0  7.262342778  \\
            1.048576e6  33.985946868  \\
        }
        ;
    \addplot
        table[row sep={\\}]
        {
            \\
            128.0  0.032104267  \\
            256.0  0.031543168  \\
            1024.0  0.031033282  \\
            2048.0  0.038629741  \\
            8192.0  0.057363631  \\
            16384.0  0.13108437  \\
            65536.0  0.247412076  \\
            131072.0  1.039419232  \\
            262144.0  4.466148948  \\
            1.048576e6  21.349092142  \\
        }
        ;
    \addplot
        table[row sep={\\}]
        {
            \\
            128.0  0.033630143  \\
            256.0  0.032148546  \\
            1024.0  0.032823547  \\
            2048.0  0.039985706  \\
            8192.0  0.065251952  \\
            16384.0  0.189983638  \\
            65536.0  0.371287496  \\
            131072.0  1.618163885  \\
            262144.0  7.094397074  \\
            1.048576e6  35.130374276  \\
        }
        ;
    \addplot
        table[row sep={\\}]
        {
            \\
            128.0  0.032901649  \\
            256.0  0.035075157  \\
            1024.0  0.033755717  \\
            2048.0  0.042235605  \\
            8192.0  0.103336303  \\
            16384.0  0.277744364  \\
            65536.0  0.552240934  \\
            131072.0  2.313051877  \\
            262144.0  10.641079911  \\
            1.048576e6  55.769864103  \\
        }
        ;
    \addplot[style={black,dashed}]
        table[row sep={\\}]
        {
            \\
            128.0  0.00014062237294951805  \\
            256.0  0.00041201610558312966  \\
            1024.0  0.0019926119926326406  \\
            2048.0  0.004499915934384578  \\
            8192.0  0.019643558074764252  \\
            16384.0  0.04202413027264982  \\
            65536.0  0.17295326115366558  \\
            131072.0  0.36340874900141257  \\
            262144.0  0.7576936141019321  \\
            1.048576e6  3.1412529644741496  \\
        }
        ;
    \legend{{},{},{},{},{},{},{},{},{},{$\mathcal O(N^{1.00}\log(N))$}}
    \addplot[style={black,dotted}]
        table[row sep={\\}]
        {
            \\
            128.0  0.009899876  \\
            256.0  0.014340951  \\
            1024.0  0.020245233  \\
            2048.0  0.024358307  \\
            8192.0  0.041900199  \\
            16384.0  0.077871171  \\
            65536.0  0.126022159  \\
            131072.0  0.449578517  \\
            262144.0  1.673088098  \\
            1.048576e6  8.333767936  \\
        }
        ;
    \addplot[style={black,dotted}]
        table[row sep={\\}]
        {
            \\
            128.0  0.017142486  \\
            256.0  0.020231018  \\
            1024.0  0.026454036  \\
            2048.0  0.03096518  \\
            8192.0  0.050358526  \\
            16384.0  0.108172684  \\
            65536.0  0.170933755  \\
            131072.0  0.945229296  \\
            262144.0  2.430539623  \\
            1.048576e6  10.973877053  \\
        }
        ;
    \addplot[style={black,dotted}]
        table[row sep={\\}]
        {
            \\
            128.0  0.021897963  \\
            256.0  0.026700621  \\
            1024.0  0.038818468  \\
            2048.0  0.039193275  \\
            8192.0  0.062070302  \\
            16384.0  0.131094383  \\
            65536.0  0.23088579  \\
            131072.0  0.75621866  \\
            262144.0  3.561206865  \\
            1.048576e6  14.173685418  \\
        }
        ;
    \addplot[style={black,dotted}]
        table[row sep={\\}]
        {
            \\
            128.0  0.004602924  \\
            256.0  0.005238237  \\
            1024.0  0.005492982  \\
            2048.0  0.010109468  \\
            8192.0  0.021395398  \\
            16384.0  0.05945554  \\
            65536.0  0.109249105  \\
            131072.0  0.4262257  \\
            262144.0  2.187695273  \\
            1.048576e6  8.989004114  \\
        }
        ;
    \addplot[style={black,dotted}]
        table[row sep={\\}]
        {
            \\
            128.0  0.004211221  \\
            256.0  0.005893489  \\
            1024.0  0.005718999  \\
            2048.0  0.0099338  \\
            8192.0  0.021205718  \\
            16384.0  0.067279071  \\
            65536.0  0.137144721  \\
            131072.0  0.481626918  \\
            262144.0  1.913128743  \\
            1.048576e6  10.437178284  \\
        }
        ;
    \addplot[style={black,dotted}]
        table[row sep={\\}]
        {
            \\
            128.0  0.004613329  \\
            256.0  0.007711363  \\
            1024.0  0.00666258  \\
            2048.0  0.009419637  \\
            8192.0  0.02282761  \\
            16384.0  0.06986202  \\
            65536.0  0.129459413  \\
            131072.0  0.550871963  \\
            262144.0  2.18945405  \\
            1.048576e6  10.284367146  \\
        }
        ;
    \addplot[style={black,dotted}]
        table[row sep={\\}]
        {
            \\
            128.0  0.008390612  \\
            256.0  0.010172008  \\
            1024.0  0.009335014  \\
            2048.0  0.013630088  \\
            8192.0  0.023891777  \\
            16384.0  0.072475184  \\
            65536.0  0.131459053  \\
            131072.0  0.52763786  \\
            262144.0  2.183443629  \\
            1.048576e6  10.204141149  \\
        }
        ;
    \addplot[style={black,dotted}]
        table[row sep={\\}]
        {
            \\
            128.0  0.008686331  \\
            256.0  0.008958512  \\
            1024.0  0.009639375  \\
            2048.0  0.013123783  \\
            8192.0  0.032531769  \\
            16384.0  0.075033653  \\
            65536.0  0.129769062  \\
            131072.0  0.533436891  \\
            262144.0  2.243546931  \\
            1.048576e6  10.575803695  \\
        }
        ;
    \addplot[style={black,dotted}]
        table[row sep={\\}]
        {
            \\
            128.0  0.0085393  \\
            256.0  0.009585086  \\
            1024.0  0.010164779  \\
            2048.0  0.013699522  \\
            8192.0  0.029041157  \\
            16384.0  0.079779141  \\
            65536.0  0.148838902  \\
            131072.0  0.611762917  \\
            262144.0  2.888624073  \\
            1.048576e6  12.759567981  \\
        }
        ;
    \nextgroupplot[]
    \addplot
        table[row sep={\\}]
        {
            \\
            64.0  0.056576514  \\
            256.0  0.118264754  \\
            1024.0  0.562381487  \\
            4096.0  5.048140974  \\
            16384.0  53.106818328  \\
            65536.0  512.187825592  \\
            262144.0  0.0  \\
            1.048576e6  0.0  \\
        }
        ;
    \addplot
        table[row sep={\\}]
        {
            \\
            64.0  0.049194495  \\
            256.0  0.152783645  \\
            1024.0  1.088946741  \\
            4096.0  14.044710851  \\
            16384.0  192.72162135  \\
            65536.0  2147.305844798  \\
            262144.0  0.0  \\
            1.048576e6  0.0  \\
        }
        ;
    \addplot
        table[row sep={\\}]
        {
            \\
            64.0  0.065231312  \\
            256.0  0.188725191  \\
            1024.0  1.530841872  \\
            4096.0  28.258512301  \\
            16384.0  437.82847934  \\
            65536.0  5315.291758423  \\
            262144.0  0.0  \\
            1.048576e6  0.0  \\
        }
        ;
    \addplot
        table[row sep={\\}]
        {
            \\
            64.0  0.029915049  \\
            256.0  0.079484915  \\
            1024.0  0.539204651  \\
            4096.0  4.959143755  \\
            16384.0  49.451855151  \\
            65536.0  477.210915073  \\
            262144.0  0.0  \\
            1.048576e6  0.0  \\
        }
        ;
    \addplot
        table[row sep={\\}]
        {
            \\
            64.0  0.044524478  \\
            256.0  0.097913618  \\
            1024.0  0.822628081  \\
            4096.0  8.205749447  \\
            16384.0  86.502309052  \\
            65536.0  830.877618906  \\
            262144.0  0.0  \\
            1.048576e6  0.0  \\
        }
        ;
    \addplot
        table[row sep={\\}]
        {
            \\
            64.0  0.047894964  \\
            256.0  0.118429054  \\
            1024.0  1.090181744  \\
            4096.0  11.895675304  \\
            16384.0  125.923547877  \\
            65536.0  1212.52554859  \\
            262144.0  0.0  \\
            1.048576e6  0.0  \\
        }
        ;
    \addplot
        table[row sep={\\}]
        {
            \\
            64.0  0.036526379  \\
            256.0  0.108763292  \\
            1024.0  0.885114279  \\
            4096.0  10.172692738  \\
            16384.0  129.75150655  \\
            65536.0  1332.459703915  \\
            262144.0  0.0  \\
            1.048576e6  0.0  \\
        }
        ;
    \addplot
        table[row sep={\\}]
        {
            \\
            64.0  0.041222356  \\
            256.0  0.134769004  \\
            1024.0  1.168860154  \\
            4096.0  14.4037286  \\
            16384.0  180.329564135  \\
            65536.0  1867.760056217  \\
            262144.0  0.0  \\
            1.048576e6  0.0  \\
        }
        ;
    \addplot
        table[row sep={\\}]
        {
            \\
            64.0  0.044817987  \\
            256.0  0.150062061  \\
            1024.0  1.458198166  \\
            4096.0  18.968794882  \\
            16384.0  239.451459821  \\
            65536.0  2518.242587119  \\
            262144.0  0.0  \\
            1.048576e6  0.0  \\
        }
        ;
    \addplot[style={black,dashed}]
        table[row sep={\\}]
        {
            \\
            64.0  0.000617957867936042  \\
            256.0  0.007871742180035812  \\
            1024.0  0.08487863562994433  \\
            4096.0  0.8415876207262829  \\
            16384.0  7.9410051063668465  \\
            65536.0  72.54597593121002  \\
        }
        ;
    \legend{{},{},{},{},{},{},{},{},{},{$\mathcal O(N^{1.50}\log(N)^2)$}}
    \addplot[style={black,dotted}]
        table[row sep={\\}]
        {
            \\
            64.0  0.019138978  \\
            256.0  0.098626629  \\
            1024.0  0.454422121  \\
            4096.0  3.123187186  \\
            16384.0  26.484014227  \\
            65536.0  255.193790657  \\
            262144.0  0.0  \\
            1.048576e6  0.0  \\
        }
        ;
    \addplot[style={black,dotted}]
        table[row sep={\\}]
        {
            \\
            64.0  0.030564859  \\
            256.0  0.240590857  \\
            1024.0  1.411568874  \\
            4096.0  9.581889099  \\
            16384.0  89.429581742  \\
            65536.0  923.462167195  \\
            262144.0  0.0  \\
            1.048576e6  0.0  \\
        }
        ;
    \addplot[style={black,dotted}]
        table[row sep={\\}]
        {
            \\
            64.0  0.043453888  \\
            256.0  0.33926656  \\
            1024.0  1.978699542  \\
            4096.0  13.137903145  \\
            16384.0  117.279186228  \\
            65536.0  1194.663278371  \\
            262144.0  0.0  \\
            1.048576e6  0.0  \\
        }
        ;
    \addplot[style={black,dotted}]
        table[row sep={\\}]
        {
            \\
            64.0  0.016160767  \\
            256.0  0.079943752  \\
            1024.0  0.449783366  \\
            4096.0  3.192330432  \\
            16384.0  25.410880038  \\
            65536.0  240.618940611  \\
            262144.0  0.0  \\
            1.048576e6  0.0  \\
        }
        ;
    \addplot[style={black,dotted}]
        table[row sep={\\}]
        {
            \\
            64.0  0.018713235  \\
            256.0  0.108382887  \\
            1024.0  0.580150583  \\
            4096.0  3.962469342  \\
            16384.0  29.919168844  \\
            65536.0  268.079083205  \\
            262144.0  0.0  \\
            1.048576e6  0.0  \\
        }
        ;
    \addplot[style={black,dotted}]
        table[row sep={\\}]
        {
            \\
            64.0  0.025871177  \\
            256.0  0.129458084  \\
            1024.0  0.755191669  \\
            4096.0  4.716746636  \\
            16384.0  35.295167234  \\
            65536.0  302.152513666  \\
            262144.0  0.0  \\
            1.048576e6  0.0  \\
        }
        ;
    \addplot[style={black,dotted}]
        table[row sep={\\}]
        {
            \\
            64.0  0.022967965  \\
            256.0  0.113595252  \\
            1024.0  0.661453073  \\
            4096.0  4.587808727  \\
            16384.0  37.995362847  \\
            65536.0  362.682190331  \\
            262144.0  0.0  \\
            1.048576e6  0.0  \\
        }
        ;
    \addplot[style={black,dotted}]
        table[row sep={\\}]
        {
            \\
            64.0  0.028516398  \\
            256.0  0.144397503  \\
            1024.0  0.836187631  \\
            4096.0  5.594721182  \\
            16384.0  43.725694077  \\
            65536.0  400.497684764  \\
            262144.0  0.0  \\
            1.048576e6  0.0  \\
        }
        ;
    \addplot[style={black,dotted}]
        table[row sep={\\}]
        {
            \\
            64.0  0.035963478  \\
            256.0  0.181783422  \\
            1024.0  1.054088342  \\
            4096.0  6.505950088  \\
            16384.0  50.38282855  \\
            65536.0  440.102833138  \\
            262144.0  0.0  \\
            1.048576e6  0.0  \\
        }
        ;
    \nextgroupplot[]
    \addplot
        table[row sep={\\}]
        {
            \\
            64.0  0.046137501  \\
            512.0  0.640945697  \\
            4096.0  49.017103161  \\
            32768.0  6403.128738893  \\
        }
        ;
    \addplot
        table[row sep={\\}]
        {
            \\
            64.0  0.10637981  \\
            512.0  0.979274778  \\
            4096.0  78.50739144  \\
            32768.0  21153.021808371  \\
        }
        ;
    \addplot
        table[row sep={\\}]
        {
            \\
            64.0  0.100404271  \\
            512.0  1.418515723  \\
            4096.0  97.827736564  \\
            32768.0  25434.090457654  \\
        }
        ;
    \addplot
        table[row sep={\\}]
        {
            \\
            64.0  0.069427404  \\
            512.0  0.606718756  \\
            4096.0  45.696162735  \\
            32768.0  4819.834943504  \\
        }
        ;
    \addplot
        table[row sep={\\}]
        {
            \\
            64.0  0.07615282  \\
            512.0  0.789936568  \\
            4096.0  56.358015279  \\
            32768.0  6940.225594922  \\
        }
        ;
    \addplot
        table[row sep={\\}]
        {
            \\
            64.0  0.062915036  \\
            512.0  1.087903191  \\
            4096.0  67.968636619  \\
            32768.0  7751.61698233  \\
        }
        ;
    \addplot
        table[row sep={\\}]
        {
            \\
            64.0  0.072931585  \\
            512.0  0.82371748  \\
            4096.0  64.254159512  \\
            32768.0  11511.443801543  \\
        }
        ;
    \addplot
        table[row sep={\\}]
        {
            \\
            64.0  0.079351081  \\
            512.0  1.087513084  \\
            4096.0  74.36877903  \\
            32768.0  12176.09044375  \\
        }
        ;
    \addplot
        table[row sep={\\}]
        {
            \\
            64.0  0.088285393  \\
            512.0  1.383435539  \\
            4096.0  87.436652056  \\
            32768.0  15362.945259054  \\
        }
        ;
    \addplot[style={black,dashed}]
        table[row sep={\\}]
        {
            \\
            64.0  0.001967935545008953  \\
            512.0  0.3163944283832535  \\
            4096.0  32.242655969426686  \\
            32768.0  2781.3031323220157  \\
        }
        ;
    \legend{{},{},{},{},{},{},{},{},{},{$\mathcal O(N^{2.00}\log(N)^2)$}}
    \addplot[style={black,dotted}]
        table[row sep={\\}]
        {
            \\
            64.0  0.043804572  \\
            512.0  1.603994372  \\
            4096.0  40.521657436  \\
            32768.0  1567.190559478  \\
        }
        ;
    \addplot[style={black,dotted}]
        table[row sep={\\}]
        {
            \\
            64.0  0.0664589  \\
            512.0  2.809119305  \\
            4096.0  130.509693414  \\
            32768.0  6518.404201037  \\
        }
        ;
    \addplot[style={black,dotted}]
        table[row sep={\\}]
        {
            \\
            64.0  0.079444104  \\
            512.0  3.999944357  \\
            4096.0  175.417624062  \\
            32768.0  3393.167662544  \\
        }
        ;
    \addplot[style={black,dotted}]
        table[row sep={\\}]
        {
            \\
            64.0  0.040745216  \\
            512.0  1.509244403  \\
            4096.0  38.549607027  \\
            32768.0  1252.696788967  \\
        }
        ;
    \addplot[style={black,dotted}]
        table[row sep={\\}]
        {
            \\
            64.0  0.052778126  \\
            512.0  2.104878631  \\
            4096.0  51.968447306  \\
            32768.0  1480.905336022  \\
        }
        ;
    \addplot[style={black,dotted}]
        table[row sep={\\}]
        {
            \\
            64.0  0.064147449  \\
            512.0  3.084559341  \\
            4096.0  69.474925273  \\
            32768.0  1817.673329234  \\
        }
        ;
    \addplot[style={black,dotted}]
        table[row sep={\\}]
        {
            \\
            64.0  0.057773266  \\
            512.0  2.296781883  \\
            4096.0  106.453014747  \\
            32768.0  3931.600476457  \\
        }
        ;
    \addplot[style={black,dotted}]
        table[row sep={\\}]
        {
            \\
            64.0  0.070114044  \\
            512.0  3.27999242  \\
            4096.0  139.123403616  \\
            32768.0  4436.922176358  \\
        }
        ;
    \addplot[style={black,dotted}]
        table[row sep={\\}]
        {
            \\
            64.0  0.082277131  \\
            512.0  4.261419345  \\
            4096.0  176.5255296  \\
            32768.0  5725.115170135  \\
        }
        ;
\end{groupplot}
\end{tikzpicture}

%% file: img/AZR2timings1d-3d.tikz
\begin{tikzpicture}
\begin{groupplot}[legend style={fill opacity={0.6}, text opacity={1}}, legend pos={north west}, xlabel={$N$}, xmode={log}, ymode={log}, legend cell align={left}, cycle list name={mark list*}, group style={group size={3 by 1}}]
    \nextgroupplot[]
    \addplot
        table[row sep={\\}]
        {
            \\
            128.0  0.027939146  \\
            256.0  0.024583273  \\
            1024.0  0.033061306  \\
            2048.0  0.034609956  \\
            8192.0  0.045600383  \\
            16384.0  0.10301961  \\
            65536.0  0.13037735  \\
            131072.0  0.412444579  \\
            262144.0  1.486334213  \\
            1.048576e6  5.867627439  \\
        }
        ;
    \addplot
        table[row sep={\\}]
        {
            \\
            128.0  0.059801954  \\
            256.0  0.049743542  \\
            1024.0  0.06705925  \\
            2048.0  0.063068876  \\
            8192.0  0.075896842  \\
            16384.0  0.131411982  \\
            65536.0  0.188549522  \\
            131072.0  0.646795738  \\
            262144.0  2.134892086  \\
            1.048576e6  8.065798021  \\
        }
        ;
    \addplot
        table[row sep={\\}]
        {
            \\
            128.0  0.075512984  \\
            256.0  0.080186719  \\
            1024.0  0.081803642  \\
            2048.0  0.087189459  \\
            8192.0  0.106083456  \\
            16384.0  0.192242484  \\
            65536.0  0.260205951  \\
            131072.0  0.736124628  \\
            262144.0  3.004933602  \\
            1.048576e6  10.597766467  \\
        }
        ;
    \addplot
        table[row sep={\\}]
        {
            \\
            128.0  0.031193202  \\
            256.0  0.017303574  \\
            1024.0  0.018437838  \\
            2048.0  0.020270957  \\
            8192.0  0.029484037  \\
            16384.0  0.060319814  \\
            65536.0  0.099031967  \\
            131072.0  0.331781945  \\
            262144.0  1.302819808  \\
            1.048576e6  5.283219912  \\
        }
        ;
    \addplot
        table[row sep={\\}]
        {
            \\
            128.0  0.031555487  \\
            256.0  0.017882006  \\
            1024.0  0.019160163  \\
            2048.0  0.021050463  \\
            8192.0  0.029972146  \\
            16384.0  0.067647076  \\
            65536.0  0.115066201  \\
            131072.0  0.411669403  \\
            262144.0  1.56491553  \\
            1.048576e6  6.543870072  \\
        }
        ;
    \addplot
        table[row sep={\\}]
        {
            \\
            128.0  0.03226801  \\
            256.0  0.019988518  \\
            1024.0  0.018468182  \\
            2048.0  0.023370993  \\
            8192.0  0.037649917  \\
            16384.0  0.074734348  \\
            65536.0  0.126581224  \\
            131072.0  0.445793848  \\
            262144.0  1.908377325  \\
            1.048576e6  7.102689085  \\
        }
        ;
    \addplot
        table[row sep={\\}]
        {
            \\
            128.0  0.034812919  \\
            256.0  0.032184817  \\
            1024.0  0.03208214  \\
            2048.0  0.039657445  \\
            8192.0  0.04867808  \\
            16384.0  0.084346021  \\
            65536.0  0.132815614  \\
            131072.0  0.441039582  \\
            262144.0  1.903669783  \\
            1.048576e6  6.660420314  \\
        }
        ;
    \addplot
        table[row sep={\\}]
        {
            \\
            128.0  0.035541312  \\
            256.0  0.032230766  \\
            1024.0  0.032846505  \\
            2048.0  0.037903043  \\
            8192.0  0.051293737  \\
            16384.0  0.092194214  \\
            65536.0  0.147016481  \\
            131072.0  0.49381357  \\
            262144.0  1.869707639  \\
            1.048576e6  8.135582009  \\
        }
        ;
    \addplot
        table[row sep={\\}]
        {
            \\
            128.0  0.035267879  \\
            256.0  0.050091475  \\
            1024.0  0.034929296  \\
            2048.0  0.039857534  \\
            8192.0  0.048585919  \\
            16384.0  0.102922547  \\
            65536.0  0.187437679  \\
            131072.0  0.521997933  \\
            262144.0  2.008862049  \\
            1.048576e6  9.081065894  \\
        }
        ;
    \addplot[style={black,dashed}]
        table[row sep={\\}]
        {
            \\
            128.0  0.00128  \\
            256.0  0.00256  \\
            1024.0  0.01024  \\
            2048.0  0.02048  \\
            8192.0  0.08192  \\
            16384.0  0.16384  \\
            65536.0  0.65536  \\
            131072.0  1.31072  \\
            262144.0  2.62144  \\
            1.048576e6  10.48576  \\
        }
        ;
    \legend{{},{},{},{},{},{},{},{},{},{$\mathcal O(N^{1.00})$}}
    \addplot[style={black,dotted}]
        table[row sep={\\}]
        {
            \\
            128.0  0.009899876  \\
            256.0  0.014340951  \\
            1024.0  0.020245233  \\
            2048.0  0.024358307  \\
            8192.0  0.041900199  \\
            16384.0  0.077871171  \\
            65536.0  0.126022159  \\
            131072.0  0.449578517  \\
            262144.0  1.673088098  \\
            1.048576e6  8.333767936  \\
        }
        ;
    \addplot[style={black,dotted}]
        table[row sep={\\}]
        {
            \\
            128.0  0.017142486  \\
            256.0  0.020231018  \\
            1024.0  0.026454036  \\
            2048.0  0.03096518  \\
            8192.0  0.050358526  \\
            16384.0  0.108172684  \\
            65536.0  0.170933755  \\
            131072.0  0.945229296  \\
            262144.0  2.430539623  \\
            1.048576e6  10.973877053  \\
        }
        ;
    \addplot[style={black,dotted}]
        table[row sep={\\}]
        {
            \\
            128.0  0.021897963  \\
            256.0  0.026700621  \\
            1024.0  0.038818468  \\
            2048.0  0.039193275  \\
            8192.0  0.062070302  \\
            16384.0  0.131094383  \\
            65536.0  0.23088579  \\
            131072.0  0.75621866  \\
            262144.0  3.561206865  \\
            1.048576e6  14.173685418  \\
        }
        ;
    \addplot[style={black,dotted}]
        table[row sep={\\}]
        {
            \\
            128.0  0.004602924  \\
            256.0  0.005238237  \\
            1024.0  0.005492982  \\
            2048.0  0.010109468  \\
            8192.0  0.021395398  \\
            16384.0  0.05945554  \\
            65536.0  0.109249105  \\
            131072.0  0.4262257  \\
            262144.0  2.187695273  \\
            1.048576e6  8.989004114  \\
        }
        ;
    \addplot[style={black,dotted}]
        table[row sep={\\}]
        {
            \\
            128.0  0.004211221  \\
            256.0  0.005893489  \\
            1024.0  0.005718999  \\
            2048.0  0.0099338  \\
            8192.0  0.021205718  \\
            16384.0  0.067279071  \\
            65536.0  0.137144721  \\
            131072.0  0.481626918  \\
            262144.0  1.913128743  \\
            1.048576e6  10.437178284  \\
        }
        ;
    \addplot[style={black,dotted}]
        table[row sep={\\}]
        {
            \\
            128.0  0.004613329  \\
            256.0  0.007711363  \\
            1024.0  0.00666258  \\
            2048.0  0.009419637  \\
            8192.0  0.02282761  \\
            16384.0  0.06986202  \\
            65536.0  0.129459413  \\
            131072.0  0.550871963  \\
            262144.0  2.18945405  \\
            1.048576e6  10.284367146  \\
        }
        ;
    \addplot[style={black,dotted}]
        table[row sep={\\}]
        {
            \\
            128.0  0.008390612  \\
            256.0  0.010172008  \\
            1024.0  0.009335014  \\
            2048.0  0.013630088  \\
            8192.0  0.023891777  \\
            16384.0  0.072475184  \\
            65536.0  0.131459053  \\
            131072.0  0.52763786  \\
            262144.0  2.183443629  \\
            1.048576e6  10.204141149  \\
        }
        ;
    \addplot[style={black,dotted}]
        table[row sep={\\}]
        {
            \\
            128.0  0.008686331  \\
            256.0  0.008958512  \\
            1024.0  0.009639375  \\
            2048.0  0.013123783  \\
            8192.0  0.032531769  \\
            16384.0  0.075033653  \\
            65536.0  0.129769062  \\
            131072.0  0.533436891  \\
            262144.0  2.243546931  \\
            1.048576e6  10.575803695  \\
        }
        ;
    \addplot[style={black,dotted}]
        table[row sep={\\}]
        {
            \\
            128.0  0.0085393  \\
            256.0  0.009585086  \\
            1024.0  0.010164779  \\
            2048.0  0.013699522  \\
            8192.0  0.029041157  \\
            16384.0  0.079779141  \\
            65536.0  0.148838902  \\
            131072.0  0.611762917  \\
            262144.0  2.888624073  \\
            1.048576e6  12.759567981  \\
        }
        ;
    \nextgroupplot[]
    \addplot
        table[row sep={\\}]
        {
            \\
            64.0  0.037363809  \\
            256.0  0.129675316  \\
            1024.0  0.490512019  \\
            4096.0  2.988598512  \\
            16384.0  21.692907044  \\
            65536.0  165.694992535  \\
            262144.0  0.0  \\
            1.048576e6  0.0  \\
        }
        ;
    \addplot
        table[row sep={\\}]
        {
            \\
            64.0  0.06298549  \\
            256.0  0.299770763  \\
            1024.0  1.469789632  \\
            4096.0  9.359693171  \\
            16384.0  70.839099651  \\
            65536.0  545.766034318  \\
            262144.0  0.0  \\
            1.048576e6  0.0  \\
        }
        ;
    \addplot
        table[row sep={\\}]
        {
            \\
            64.0  0.080940838  \\
            256.0  0.405282879  \\
            1024.0  2.074103325  \\
            4096.0  13.261050742  \\
            16384.0  98.888064533  \\
            65536.0  773.008723492  \\
            262144.0  0.0  \\
            1.048576e6  0.0  \\
        }
        ;
    \addplot
        table[row sep={\\}]
        {
            \\
            64.0  0.033609684  \\
            256.0  0.105222211  \\
            1024.0  0.468625687  \\
            4096.0  2.918111079  \\
            16384.0  20.284184485  \\
            65536.0  157.044999864  \\
            262144.0  0.0  \\
            1.048576e6  0.0  \\
        }
        ;
    \addplot
        table[row sep={\\}]
        {
            \\
            64.0  0.037789995  \\
            256.0  0.130306863  \\
            1024.0  0.595921652  \\
            4096.0  3.90306844  \\
            16384.0  26.039977853  \\
            65536.0  194.022848639  \\
            262144.0  0.0  \\
            1.048576e6  0.0  \\
        }
        ;
    \addplot
        table[row sep={\\}]
        {
            \\
            64.0  0.044468445  \\
            256.0  0.170732079  \\
            1024.0  0.774501183  \\
            4096.0  4.646368331  \\
            16384.0  31.594873981  \\
            65536.0  234.462088463  \\
            262144.0  0.0  \\
            1.048576e6  0.0  \\
        }
        ;
    \addplot
        table[row sep={\\}]
        {
            \\
            64.0  0.050805153  \\
            256.0  0.150231661  \\
            1024.0  0.692563389  \\
            4096.0  4.351103516  \\
            16384.0  32.080777101  \\
            65536.0  243.692554659  \\
            262144.0  0.0  \\
            1.048576e6  0.0  \\
        }
        ;
    \addplot
        table[row sep={\\}]
        {
            \\
            64.0  0.059211904  \\
            256.0  0.179694843  \\
            1024.0  0.873479647  \\
            4096.0  5.340812488  \\
            16384.0  38.734456389  \\
            65536.0  291.096097486  \\
            262144.0  0.0  \\
            1.048576e6  0.0  \\
        }
        ;
    \addplot
        table[row sep={\\}]
        {
            \\
            64.0  0.066712142  \\
            256.0  0.217208491  \\
            1024.0  1.084570803  \\
            4096.0  6.483454941  \\
            16384.0  46.387085803  \\
            65536.0  339.935983784  \\
            262144.0  0.0  \\
            1.048576e6  0.0  \\
        }
        ;
    \addplot[style={black,dashed}]
        table[row sep={\\}]
        {
            \\
            64.0  0.00256  \\
            256.0  0.02048  \\
            1024.0  0.16384  \\
            4096.0  1.31072  \\
            16384.0  10.48576  \\
            65536.0  83.88608  \\
        }
        ;
    \legend{{},{},{},{},{},{},{},{},{},{$\mathcal O(N^{1.50})$}}
    \addplot[style={black,dotted}]
        table[row sep={\\}]
        {
            \\
            64.0  0.019138978  \\
            256.0  0.098626629  \\
            1024.0  0.454422121  \\
            4096.0  3.123187186  \\
            16384.0  26.484014227  \\
            65536.0  255.193790657  \\
            262144.0  0.0  \\
            1.048576e6  0.0  \\
        }
        ;
    \addplot[style={black,dotted}]
        table[row sep={\\}]
        {
            \\
            64.0  0.030564859  \\
            256.0  0.240590857  \\
            1024.0  1.411568874  \\
            4096.0  9.581889099  \\
            16384.0  89.429581742  \\
            65536.0  923.462167195  \\
            262144.0  0.0  \\
            1.048576e6  0.0  \\
        }
        ;
    \addplot[style={black,dotted}]
        table[row sep={\\}]
        {
            \\
            64.0  0.043453888  \\
            256.0  0.33926656  \\
            1024.0  1.978699542  \\
            4096.0  13.137903145  \\
            16384.0  117.279186228  \\
            65536.0  1194.663278371  \\
            262144.0  0.0  \\
            1.048576e6  0.0  \\
        }
        ;
    \addplot[style={black,dotted}]
        table[row sep={\\}]
        {
            \\
            64.0  0.016160767  \\
            256.0  0.079943752  \\
            1024.0  0.449783366  \\
            4096.0  3.192330432  \\
            16384.0  25.410880038  \\
            65536.0  240.618940611  \\
            262144.0  0.0  \\
            1.048576e6  0.0  \\
        }
        ;
    \addplot[style={black,dotted}]
        table[row sep={\\}]
        {
            \\
            64.0  0.018713235  \\
            256.0  0.108382887  \\
            1024.0  0.580150583  \\
            4096.0  3.962469342  \\
            16384.0  29.919168844  \\
            65536.0  268.079083205  \\
            262144.0  0.0  \\
            1.048576e6  0.0  \\
        }
        ;
    \addplot[style={black,dotted}]
        table[row sep={\\}]
        {
            \\
            64.0  0.025871177  \\
            256.0  0.129458084  \\
            1024.0  0.755191669  \\
            4096.0  4.716746636  \\
            16384.0  35.295167234  \\
            65536.0  302.152513666  \\
            262144.0  0.0  \\
            1.048576e6  0.0  \\
        }
        ;
    \addplot[style={black,dotted}]
        table[row sep={\\}]
        {
            \\
            64.0  0.022967965  \\
            256.0  0.113595252  \\
            1024.0  0.661453073  \\
            4096.0  4.587808727  \\
            16384.0  37.995362847  \\
            65536.0  362.682190331  \\
            262144.0  0.0  \\
            1.048576e6  0.0  \\
        }
        ;
    \addplot[style={black,dotted}]
        table[row sep={\\}]
        {
            \\
            64.0  0.028516398  \\
            256.0  0.144397503  \\
            1024.0  0.836187631  \\
            4096.0  5.594721182  \\
            16384.0  43.725694077  \\
            65536.0  400.497684764  \\
            262144.0  0.0  \\
            1.048576e6  0.0  \\
        }
        ;
    \addplot[style={black,dotted}]
        table[row sep={\\}]
        {
            \\
            64.0  0.035963478  \\
            256.0  0.181783422  \\
            1024.0  1.054088342  \\
            4096.0  6.505950088  \\
            16384.0  50.38282855  \\
            65536.0  440.102833138  \\
            262144.0  0.0  \\
            1.048576e6  0.0  \\
        }
        ;
    \nextgroupplot[]
    \addplot
        table[row sep={\\}]
        {
            \\
            64.0  0.079672346  \\
            512.0  1.654888631  \\
            4096.0  41.519971792  \\
            32768.0  1435.674066914  \\
        }
        ;
    \addplot
        table[row sep={\\}]
        {
            \\
            64.0  0.123150415  \\
            512.0  2.901869838  \\
            4096.0  132.461603419  \\
            32768.0  6244.613953765  \\
        }
        ;
    \addplot
        table[row sep={\\}]
        {
            \\
            64.0  0.139056705  \\
            512.0  4.22519171  \\
            4096.0  179.942765089  \\
            32768.0  3435.862462515  \\
        }
        ;
    \addplot
        table[row sep={\\}]
        {
            \\
            64.0  0.074595379  \\
            512.0  1.548492023  \\
            4096.0  40.293632129  \\
            32768.0  1120.432404464  \\
        }
        ;
    \addplot
        table[row sep={\\}]
        {
            \\
            64.0  0.090928457  \\
            512.0  2.165700197  \\
            4096.0  52.666816863  \\
            32768.0  1438.012365318  \\
        }
        ;
    \addplot
        table[row sep={\\}]
        {
            \\
            64.0  0.102510194  \\
            512.0  3.125498993  \\
            4096.0  71.286416044  \\
            32768.0  1779.262631705  \\
        }
        ;
    \addplot
        table[row sep={\\}]
        {
            \\
            64.0  0.106275277  \\
            512.0  2.403090091  \\
            4096.0  108.79914676  \\
            32768.0  3908.389439779  \\
        }
        ;
    \addplot
        table[row sep={\\}]
        {
            \\
            64.0  0.121295429  \\
            512.0  3.398403529  \\
            4096.0  141.412618188  \\
            32768.0  4502.208914758  \\
        }
        ;
    \addplot
        table[row sep={\\}]
        {
            \\
            64.0  0.136964715  \\
            512.0  4.357564809  \\
            4096.0  179.589013309  \\
            32768.0  6261.553602863  \\
        }
        ;
    \addplot[style={black,dashed}]
        table[row sep={\\}]
        {
            \\
            64.0  0.0004096  \\
            512.0  0.0262144  \\
            4096.0  1.6777216  \\
            32768.0  107.3741824  \\
        }
        ;
    \legend{{},{},{},{},{},{},{},{},{},{$\mathcal O(N^{2.00})$}}
    \addplot[style={black,dotted}]
        table[row sep={\\}]
        {
            \\
            64.0  0.043804572  \\
            512.0  1.603994372  \\
            4096.0  40.521657436  \\
            32768.0  1567.190559478  \\
        }
        ;
    \addplot[style={black,dotted}]
        table[row sep={\\}]
        {
            \\
            64.0  0.0664589  \\
            512.0  2.809119305  \\
            4096.0  130.509693414  \\
            32768.0  6518.404201037  \\
        }
        ;
    \addplot[style={black,dotted}]
        table[row sep={\\}]
        {
            \\
            64.0  0.079444104  \\
            512.0  3.999944357  \\
            4096.0  175.417624062  \\
            32768.0  3393.167662544  \\
        }
        ;
    \addplot[style={black,dotted}]
        table[row sep={\\}]
        {
            \\
            64.0  0.040745216  \\
            512.0  1.509244403  \\
            4096.0  38.549607027  \\
            32768.0  1252.696788967  \\
        }
        ;
    \addplot[style={black,dotted}]
        table[row sep={\\}]
        {
            \\
            64.0  0.052778126  \\
            512.0  2.104878631  \\
            4096.0  51.968447306  \\
            32768.0  1480.905336022  \\
        }
        ;
    \addplot[style={black,dotted}]
        table[row sep={\\}]
        {
            \\
            64.0  0.064147449  \\
            512.0  3.084559341  \\
            4096.0  69.474925273  \\
            32768.0  1817.673329234  \\
        }
        ;
    \addplot[style={black,dotted}]
        table[row sep={\\}]
        {
            \\
            64.0  0.057773266  \\
            512.0  2.296781883  \\
            4096.0  106.453014747  \\
            32768.0  3931.600476457  \\
        }
        ;
    \addplot[style={black,dotted}]
        table[row sep={\\}]
        {
            \\
            64.0  0.070114044  \\
            512.0  3.27999242  \\
            4096.0  139.123403616  \\
            32768.0  4436.922176358  \\
        }
        ;
    \addplot[style={black,dotted}]
        table[row sep={\\}]
        {
            \\
            64.0  0.082277131  \\
            512.0  4.261419345  \\
            4096.0  176.5255296  \\
            32768.0  5725.115170135  \\
        }
        ;
\end{groupplot}
\end{tikzpicture}

%% file: img/AZR2timings2d.tikz
\begin{tikzpicture}
\begin{groupplot}[xmode={log}, ymode={log}, legend cell align={left}, legend style={fill opacity={0.6}, text opacity={1}}, legend pos={north west}, width={0.5\textwidth}, height={0.4\textwidth}, xlabel={$N$}, ymin={0.01}, ymax={1000.0}, ylabel={Timings}, cycle list name={mark list*}, group style={x descriptions at={edge bottom}, ylabels at={edge left}, y descriptions at={edge left}, horizontal sep={1em}, group size={3 by 1}}]
    \nextgroupplot[]
    \addplot
        table[row sep={\\}]
        {
            \\
            64.0  0.037363809  \\
            256.0  0.129675316  \\
            1024.0  0.490512019  \\
            4096.0  2.988598512  \\
            16384.0  21.692907044  \\
            65536.0  165.694992535  \\
            262144.0  0.0  \\
            1.048576e6  0.0  \\
        }
        ;
    \addlegendentry {$\texttt{db2}$}
    \addplot
        table[row sep={\\}]
        {
            \\
            64.0  0.06298549  \\
            256.0  0.299770763  \\
            1024.0  1.469789632  \\
            4096.0  9.359693171  \\
            16384.0  70.839099651  \\
            65536.0  545.766034318  \\
            262144.0  0.0  \\
            1.048576e6  0.0  \\
        }
        ;
    \addlegendentry {$\texttt{db3}$}
    \addplot
        table[row sep={\\}]
        {
            \\
            64.0  0.080940838  \\
            256.0  0.405282879  \\
            1024.0  2.074103325  \\
            4096.0  13.261050742  \\
            16384.0  98.888064533  \\
            65536.0  773.008723492  \\
            262144.0  0.0  \\
            1.048576e6  0.0  \\
        }
        ;
    \addlegendentry {$\texttt{db4}$}
    \nextgroupplot[]
    \addplot
        table[row sep={\\}]
        {
            \\
            64.0  0.033609684  \\
            256.0  0.105222211  \\
            1024.0  0.468625687  \\
            4096.0  2.918111079  \\
            16384.0  20.284184485  \\
            65536.0  157.044999864  \\
            262144.0  0.0  \\
            1.048576e6  0.0  \\
        }
        ;
    \addlegendentry {$\texttt{cdf31}$}
    \addplot
        table[row sep={\\}]
        {
            \\
            64.0  0.037789995  \\
            256.0  0.130306863  \\
            1024.0  0.595921652  \\
            4096.0  3.90306844  \\
            16384.0  26.039977853  \\
            65536.0  194.022848639  \\
            262144.0  0.0  \\
            1.048576e6  0.0  \\
        }
        ;
    \addlegendentry {$\texttt{cdf33}$}
    \addplot
        table[row sep={\\}]
        {
            \\
            64.0  0.044468445  \\
            256.0  0.170732079  \\
            1024.0  0.774501183  \\
            4096.0  4.646368331  \\
            16384.0  31.594873981  \\
            65536.0  234.462088463  \\
            262144.0  0.0  \\
            1.048576e6  0.0  \\
        }
        ;
    \addlegendentry {$\texttt{cdf35}$}
    \nextgroupplot[]
    \addplot
        table[row sep={\\}]
        {
            \\
            64.0  0.050805153  \\
            256.0  0.150231661  \\
            1024.0  0.692563389  \\
            4096.0  4.351103516  \\
            16384.0  32.080777101  \\
            65536.0  243.692554659  \\
            262144.0  0.0  \\
            1.048576e6  0.0  \\
        }
        ;
    \addlegendentry {$\texttt{cdf42}$}
    \addplot
        table[row sep={\\}]
        {
            \\
            64.0  0.059211904  \\
            256.0  0.179694843  \\
            1024.0  0.873479647  \\
            4096.0  5.340812488  \\
            16384.0  38.734456389  \\
            65536.0  291.096097486  \\
            262144.0  0.0  \\
            1.048576e6  0.0  \\
        }
        ;
    \addlegendentry {$\texttt{cdf44}$}
    \addplot
        table[row sep={\\}]
        {
            \\
            64.0  0.066712142  \\
            256.0  0.217208491  \\
            1024.0  1.084570803  \\
            4096.0  6.483454941  \\
            16384.0  46.387085803  \\
            65536.0  339.935983784  \\
            262144.0  0.0  \\
            1.048576e6  0.0  \\
        }
        ;
    \addlegendentry {$\texttt{cdf46}$}
\end{groupplot}
\end{tikzpicture}

%% file: img/AZStimings1d-3d.tikz
\begin{tikzpicture}
\begin{groupplot}[legend style={fill opacity={0.6}, text opacity={1}}, legend pos={north west}, xlabel={$N$}, xmode={log}, ymode={log}, legend cell align={left}, cycle list name={mark list*}, group style={group size={3 by 1}}]
    \nextgroupplot[]
    \addplot
        table[row sep={\\}]
        {
            \\
            128.0  0.01641736  \\
            256.0  0.018352284  \\
            1024.0  0.027898474  \\
            2048.0  0.026889006  \\
            8192.0  0.036481972  \\
            16384.0  0.067785835  \\
            65536.0  0.133567418  \\
            131072.0  0.503214514  \\
            262144.0  1.765784313  \\
            1.048576e6  8.670407272  \\
        }
        ;
    \addplot
        table[row sep={\\}]
        {
            \\
            128.0  0.03368338  \\
            256.0  0.034692052  \\
            1024.0  0.045759821  \\
            2048.0  0.045524327  \\
            8192.0  0.054495011  \\
            16384.0  0.116403498  \\
            65536.0  0.171331663  \\
            131072.0  0.589669178  \\
            262144.0  2.401287696  \\
            1.048576e6  11.986890966  \\
        }
        ;
    \addplot
        table[row sep={\\}]
        {
            \\
            128.0  0.05225485  \\
            256.0  0.059210613  \\
            1024.0  0.057853528  \\
            2048.0  0.067934309  \\
            8192.0  0.099882031  \\
            16384.0  0.121874484  \\
            65536.0  0.198053095  \\
            131072.0  0.642503399  \\
            262144.0  2.697758827  \\
            1.048576e6  11.857690128  \\
        }
        ;
    \addplot
        table[row sep={\\}]
        {
            \\
            128.0  0.011101017  \\
            256.0  0.010014863  \\
            1024.0  0.010530533  \\
            2048.0  0.011717069  \\
            8192.0  0.016730555  \\
            16384.0  0.04527034  \\
            65536.0  0.11996116  \\
            131072.0  0.341034974  \\
            262144.0  1.706501353  \\
            1.048576e6  8.622987773  \\
        }
        ;
    \addplot
        table[row sep={\\}]
        {
            \\
            128.0  0.010013508  \\
            256.0  0.010552498  \\
            1024.0  0.010566111  \\
            2048.0  0.011877867  \\
            8192.0  0.017115778  \\
            16384.0  0.320062267  \\
            65536.0  0.098502603  \\
            131072.0  0.443886008  \\
            262144.0  2.054343688  \\
            1.048576e6  10.366463651  \\
        }
        ;
    \addplot
        table[row sep={\\}]
        {
            \\
            128.0  0.01009407  \\
            256.0  0.011613476  \\
            1024.0  0.010469314  \\
            2048.0  0.011667276  \\
            8192.0  0.018048597  \\
            16384.0  0.054719619  \\
            65536.0  0.4428434  \\
            131072.0  0.468008155  \\
            262144.0  2.134378046  \\
            1.048576e6  9.601493537  \\
        }
        ;
    \addplot
        table[row sep={\\}]
        {
            \\
            128.0  0.019293235  \\
            256.0  0.019053361  \\
            1024.0  0.018807343  \\
            2048.0  0.020925962  \\
            8192.0  0.026537828  \\
            16384.0  0.079413595  \\
            65536.0  0.125417799  \\
            131072.0  0.448438626  \\
            262144.0  2.139283054  \\
            1.048576e6  10.851535373  \\
        }
        ;
    \addplot
        table[row sep={\\}]
        {
            \\
            128.0  0.019442165  \\
            256.0  0.020055778  \\
            1024.0  0.01967791  \\
            2048.0  0.021374131  \\
            8192.0  0.025653927  \\
            16384.0  0.0788598  \\
            65536.0  0.107473955  \\
            131072.0  0.481610282  \\
            262144.0  2.228423843  \\
            1.048576e6  10.9973603  \\
        }
        ;
    \addplot
        table[row sep={\\}]
        {
            \\
            128.0  0.020747355  \\
            256.0  0.023889096  \\
            1024.0  0.020442163  \\
            2048.0  0.021369848  \\
            8192.0  0.031200173  \\
            16384.0  0.068330119  \\
            65536.0  0.12474847  \\
            131072.0  0.536885218  \\
            262144.0  2.448845502  \\
            1.048576e6  12.807178525  \\
        }
        ;
    \addplot[style={black,dashed}]
        table[row sep={\\}]
        {
            \\
            128.0  0.00128  \\
            256.0  0.00256  \\
            1024.0  0.01024  \\
            2048.0  0.02048  \\
            8192.0  0.08192  \\
            16384.0  0.16384  \\
            65536.0  0.65536  \\
            131072.0  1.31072  \\
            262144.0  2.62144  \\
            1.048576e6  10.48576  \\
        }
        ;
    \legend{{},{},{},{},{},{},{},{},{},{$\mathcal O(N^{1.00})$}}
    \addplot[style={black,dotted}]
        table[row sep={\\}]
        {
            \\
            128.0  0.009899876  \\
            256.0  0.014340951  \\
            1024.0  0.020245233  \\
            2048.0  0.024358307  \\
            8192.0  0.041900199  \\
            16384.0  0.077871171  \\
            65536.0  0.126022159  \\
            131072.0  0.449578517  \\
            262144.0  1.673088098  \\
            1.048576e6  8.333767936  \\
        }
        ;
    \addplot[style={black,dotted}]
        table[row sep={\\}]
        {
            \\
            128.0  0.017142486  \\
            256.0  0.020231018  \\
            1024.0  0.026454036  \\
            2048.0  0.03096518  \\
            8192.0  0.050358526  \\
            16384.0  0.108172684  \\
            65536.0  0.170933755  \\
            131072.0  0.945229296  \\
            262144.0  2.430539623  \\
            1.048576e6  10.973877053  \\
        }
        ;
    \addplot[style={black,dotted}]
        table[row sep={\\}]
        {
            \\
            128.0  0.021897963  \\
            256.0  0.026700621  \\
            1024.0  0.038818468  \\
            2048.0  0.039193275  \\
            8192.0  0.062070302  \\
            16384.0  0.131094383  \\
            65536.0  0.23088579  \\
            131072.0  0.75621866  \\
            262144.0  3.561206865  \\
            1.048576e6  14.173685418  \\
        }
        ;
    \addplot[style={black,dotted}]
        table[row sep={\\}]
        {
            \\
            128.0  0.004602924  \\
            256.0  0.005238237  \\
            1024.0  0.005492982  \\
            2048.0  0.010109468  \\
            8192.0  0.021395398  \\
            16384.0  0.05945554  \\
            65536.0  0.109249105  \\
            131072.0  0.4262257  \\
            262144.0  2.187695273  \\
            1.048576e6  8.989004114  \\
        }
        ;
    \addplot[style={black,dotted}]
        table[row sep={\\}]
        {
            \\
            128.0  0.004211221  \\
            256.0  0.005893489  \\
            1024.0  0.005718999  \\
            2048.0  0.0099338  \\
            8192.0  0.021205718  \\
            16384.0  0.067279071  \\
            65536.0  0.137144721  \\
            131072.0  0.481626918  \\
            262144.0  1.913128743  \\
            1.048576e6  10.437178284  \\
        }
        ;
    \addplot[style={black,dotted}]
        table[row sep={\\}]
        {
            \\
            128.0  0.004613329  \\
            256.0  0.007711363  \\
            1024.0  0.00666258  \\
            2048.0  0.009419637  \\
            8192.0  0.02282761  \\
            16384.0  0.06986202  \\
            65536.0  0.129459413  \\
            131072.0  0.550871963  \\
            262144.0  2.18945405  \\
            1.048576e6  10.284367146  \\
        }
        ;
    \addplot[style={black,dotted}]
        table[row sep={\\}]
        {
            \\
            128.0  0.008390612  \\
            256.0  0.010172008  \\
            1024.0  0.009335014  \\
            2048.0  0.013630088  \\
            8192.0  0.023891777  \\
            16384.0  0.072475184  \\
            65536.0  0.131459053  \\
            131072.0  0.52763786  \\
            262144.0  2.183443629  \\
            1.048576e6  10.204141149  \\
        }
        ;
    \addplot[style={black,dotted}]
        table[row sep={\\}]
        {
            \\
            128.0  0.008686331  \\
            256.0  0.008958512  \\
            1024.0  0.009639375  \\
            2048.0  0.013123783  \\
            8192.0  0.032531769  \\
            16384.0  0.075033653  \\
            65536.0  0.129769062  \\
            131072.0  0.533436891  \\
            262144.0  2.243546931  \\
            1.048576e6  10.575803695  \\
        }
        ;
    \addplot[style={black,dotted}]
        table[row sep={\\}]
        {
            \\
            128.0  0.0085393  \\
            256.0  0.009585086  \\
            1024.0  0.010164779  \\
            2048.0  0.013699522  \\
            8192.0  0.029041157  \\
            16384.0  0.079779141  \\
            65536.0  0.148838902  \\
            131072.0  0.611762917  \\
            262144.0  2.888624073  \\
            1.048576e6  12.759567981  \\
        }
        ;
    \nextgroupplot[]
    \addplot
        table[row sep={\\}]
        {
            \\
            64.0  0.069172184  \\
            256.0  0.056223112  \\
            1024.0  0.077659512  \\
            4096.0  0.200994781  \\
            16384.0  1.061347391  \\
            65536.0  3.583179499  \\
            262144.0  15.437803005  \\
            1.048576e6  0.0  \\
        }
        ;
    \addplot
        table[row sep={\\}]
        {
            \\
            64.0  0.028147923  \\
            256.0  0.066569156  \\
            1024.0  0.207162927  \\
            4096.0  1.10022845  \\
            16384.0  5.120571052  \\
            65536.0  19.34280771  \\
            262144.0  85.392188779  \\
            1.048576e6  0.0  \\
        }
        ;
    \addplot
        table[row sep={\\}]
        {
            \\
            64.0  0.040967134  \\
            256.0  0.105294947  \\
            1024.0  0.662289958  \\
            4096.0  2.73809348  \\
            16384.0  14.73883701  \\
            65536.0  69.205879867  \\
            262144.0  206.817897228  \\
            1.048576e6  0.0  \\
        }
        ;
    \addplot
        table[row sep={\\}]
        {
            \\
            64.0  0.008370861  \\
            256.0  0.029689664  \\
            1024.0  0.055170232  \\
            4096.0  0.20068383  \\
            16384.0  1.014020454  \\
            65536.0  3.775359867  \\
            262144.0  14.159454617  \\
            1.048576e6  0.0  \\
        }
        ;
    \addplot
        table[row sep={\\}]
        {
            \\
            64.0  0.007768416  \\
            256.0  0.03179759  \\
            1024.0  0.08131809  \\
            4096.0  0.362737308  \\
            16384.0  1.848973761  \\
            65536.0  9.556342632  \\
            262144.0  44.954903727  \\
            1.048576e6  0.0  \\
        }
        ;
    \addplot
        table[row sep={\\}]
        {
            \\
            64.0  0.00799658  \\
            256.0  0.033874264  \\
            1024.0  0.09803432  \\
            4096.0  0.574383162  \\
            16384.0  3.624448367  \\
            65536.0  18.702717658  \\
            262144.0  92.555255567  \\
            1.048576e6  0.0  \\
        }
        ;
    \addplot
        table[row sep={\\}]
        {
            \\
            64.0  0.013459952  \\
            256.0  0.041978494  \\
            1024.0  0.129496866  \\
            4096.0  0.847446714  \\
            16384.0  3.694096035  \\
            65536.0  16.278107702  \\
            262144.0  75.599888947  \\
            1.048576e6  0.0  \\
        }
        ;
    \addplot
        table[row sep={\\}]
        {
            \\
            64.0  0.013236752  \\
            256.0  0.044441684  \\
            1024.0  0.142431565  \\
            4096.0  1.199708221  \\
            16384.0  7.426851056  \\
            65536.0  32.178535774  \\
            262144.0  160.518755883  \\
            1.048576e6  0.0  \\
        }
        ;
    \addplot
        table[row sep={\\}]
        {
            \\
            64.0  0.03234988  \\
            256.0  0.047645558  \\
            1024.0  0.217746727  \\
            4096.0  1.507180827  \\
            16384.0  10.005199867  \\
            65536.0  53.306766365  \\
            262144.0  266.500990847  \\
            1.048576e6  0.0  \\
        }
        ;
    \addplot[style={black,dashed}]
        table[row sep={\\}]
        {
            \\
            64.0  0.0064  \\
            256.0  0.0256  \\
            1024.0  0.1024  \\
            4096.0  0.4096  \\
            16384.0  1.6384  \\
            65536.0  6.5536  \\
            262144.0  26.2144  \\
        }
        ;
    \legend{{},{},{},{},{},{},{},{},{},{$\mathcal O(N^{1.00})$}}
    \addplot[style={black,dotted}]
        table[row sep={\\}]
        {
            \\
            64.0  0.019138978  \\
            256.0  0.098626629  \\
            1024.0  0.454422121  \\
            4096.0  3.123187186  \\
            16384.0  26.484014227  \\
            65536.0  255.193790657  \\
            262144.0  0.0  \\
            1.048576e6  0.0  \\
        }
        ;
    \addplot[style={black,dotted}]
        table[row sep={\\}]
        {
            \\
            64.0  0.030564859  \\
            256.0  0.240590857  \\
            1024.0  1.411568874  \\
            4096.0  9.581889099  \\
            16384.0  89.429581742  \\
            65536.0  923.462167195  \\
            262144.0  0.0  \\
            1.048576e6  0.0  \\
        }
        ;
    \addplot[style={black,dotted}]
        table[row sep={\\}]
        {
            \\
            64.0  0.043453888  \\
            256.0  0.33926656  \\
            1024.0  1.978699542  \\
            4096.0  13.137903145  \\
            16384.0  117.279186228  \\
            65536.0  1194.663278371  \\
            262144.0  0.0  \\
            1.048576e6  0.0  \\
        }
        ;
    \addplot[style={black,dotted}]
        table[row sep={\\}]
        {
            \\
            64.0  0.016160767  \\
            256.0  0.079943752  \\
            1024.0  0.449783366  \\
            4096.0  3.192330432  \\
            16384.0  25.410880038  \\
            65536.0  240.618940611  \\
            262144.0  0.0  \\
            1.048576e6  0.0  \\
        }
        ;
    \addplot[style={black,dotted}]
        table[row sep={\\}]
        {
            \\
            64.0  0.018713235  \\
            256.0  0.108382887  \\
            1024.0  0.580150583  \\
            4096.0  3.962469342  \\
            16384.0  29.919168844  \\
            65536.0  268.079083205  \\
            262144.0  0.0  \\
            1.048576e6  0.0  \\
        }
        ;
    \addplot[style={black,dotted}]
        table[row sep={\\}]
        {
            \\
            64.0  0.025871177  \\
            256.0  0.129458084  \\
            1024.0  0.755191669  \\
            4096.0  4.716746636  \\
            16384.0  35.295167234  \\
            65536.0  302.152513666  \\
            262144.0  0.0  \\
            1.048576e6  0.0  \\
        }
        ;
    \addplot[style={black,dotted}]
        table[row sep={\\}]
        {
            \\
            64.0  0.022967965  \\
            256.0  0.113595252  \\
            1024.0  0.661453073  \\
            4096.0  4.587808727  \\
            16384.0  37.995362847  \\
            65536.0  362.682190331  \\
            262144.0  0.0  \\
            1.048576e6  0.0  \\
        }
        ;
    \addplot[style={black,dotted}]
        table[row sep={\\}]
        {
            \\
            64.0  0.028516398  \\
            256.0  0.144397503  \\
            1024.0  0.836187631  \\
            4096.0  5.594721182  \\
            16384.0  43.725694077  \\
            65536.0  400.497684764  \\
            262144.0  0.0  \\
            1.048576e6  0.0  \\
        }
        ;
    \addplot[style={black,dotted}]
        table[row sep={\\}]
        {
            \\
            64.0  0.035963478  \\
            256.0  0.181783422  \\
            1024.0  1.054088342  \\
            4096.0  6.505950088  \\
            16384.0  50.38282855  \\
            65536.0  440.102833138  \\
            262144.0  0.0  \\
            1.048576e6  0.0  \\
        }
        ;
    \nextgroupplot[]
    \addplot
        table[row sep={\\}]
        {
            \\
            64.0  0.083294029  \\
            512.0  0.189600045  \\
            4096.0  8.40313058  \\
            32768.0  359.846001586  \\
        }
        ;
    \addplot
        table[row sep={\\}]
        {
            \\
            64.0  0.053029292  \\
            512.0  0.535061235  \\
            4096.0  28.537315942  \\
            32768.0  1381.035493565  \\
        }
        ;
    \addplot
        table[row sep={\\}]
        {
            \\
            64.0  0.071527756  \\
            512.0  1.015401753  \\
            4096.0  77.708344632  \\
            32768.0  2825.806047609  \\
        }
        ;
    \addplot
        table[row sep={\\}]
        {
            \\
            64.0  0.027427256  \\
            512.0  0.164480158  \\
            4096.0  7.478236207  \\
            32768.0  289.450580528  \\
        }
        ;
    \addplot
        table[row sep={\\}]
        {
            \\
            64.0  0.027178229  \\
            512.0  0.23056075  \\
            4096.0  12.250088554  \\
            32768.0  589.95178702  \\
        }
        ;
    \addplot
        table[row sep={\\}]
        {
            \\
            64.0  0.028104655  \\
            512.0  0.232629574  \\
            4096.0  15.768962018  \\
            32768.0  780.803343876  \\
        }
        ;
    \addplot
        table[row sep={\\}]
        {
            \\
            64.0  0.039465131  \\
            512.0  0.376563042  \\
            4096.0  21.178938308  \\
            32768.0  1889.255021745  \\
        }
        ;
    \addplot
        table[row sep={\\}]
        {
            \\
            64.0  0.043427167  \\
            512.0  0.707447736  \\
            4096.0  27.223738699  \\
            32768.0  1509.995750552  \\
        }
        ;
    \addplot
        table[row sep={\\}]
        {
            \\
            64.0  0.062854337  \\
            512.0  0.456949504  \\
            4096.0  32.519593768  \\
            32768.0  1834.989058737  \\
        }
        ;
    \addplot[style={black,dashed}]
        table[row sep={\\}]
        {
            \\
            64.0  0.064  \\
            512.0  0.512  \\
            4096.0  4.096  \\
            32768.0  32.768  \\
        }
        ;
    \legend{{},{},{},{},{},{},{},{},{},{$\mathcal O(N^{1.00})$}}
    \addplot[style={black,dotted}]
        table[row sep={\\}]
        {
            \\
            64.0  0.043804572  \\
            512.0  1.603994372  \\
            4096.0  40.521657436  \\
            32768.0  1567.190559478  \\
        }
        ;
    \addplot[style={black,dotted}]
        table[row sep={\\}]
        {
            \\
            64.0  0.0664589  \\
            512.0  2.809119305  \\
            4096.0  130.509693414  \\
            32768.0  6518.404201037  \\
        }
        ;
    \addplot[style={black,dotted}]
        table[row sep={\\}]
        {
            \\
            64.0  0.079444104  \\
            512.0  3.999944357  \\
            4096.0  175.417624062  \\
            32768.0  3393.167662544  \\
        }
        ;
    \addplot[style={black,dotted}]
        table[row sep={\\}]
        {
            \\
            64.0  0.040745216  \\
            512.0  1.509244403  \\
            4096.0  38.549607027  \\
            32768.0  1252.696788967  \\
        }
        ;
    \addplot[style={black,dotted}]
        table[row sep={\\}]
        {
            \\
            64.0  0.052778126  \\
            512.0  2.104878631  \\
            4096.0  51.968447306  \\
            32768.0  1480.905336022  \\
        }
        ;
    \addplot[style={black,dotted}]
        table[row sep={\\}]
        {
            \\
            64.0  0.064147449  \\
            512.0  3.084559341  \\
            4096.0  69.474925273  \\
            32768.0  1817.673329234  \\
        }
        ;
    \addplot[style={black,dotted}]
        table[row sep={\\}]
        {
            \\
            64.0  0.057773266  \\
            512.0  2.296781883  \\
            4096.0  106.453014747  \\
            32768.0  3931.600476457  \\
        }
        ;
    \addplot[style={black,dotted}]
        table[row sep={\\}]
        {
            \\
            64.0  0.070114044  \\
            512.0  3.27999242  \\
            4096.0  139.123403616  \\
            32768.0  4436.922176358  \\
        }
        ;
    \addplot[style={black,dotted}]
        table[row sep={\\}]
        {
            \\
            64.0  0.082277131  \\
            512.0  4.261419345  \\
            4096.0  176.5255296  \\
            32768.0  5725.115170135  \\
        }
        ;
\end{groupplot}
\end{tikzpicture}

%% file: img/AZAStimings1d-3d.tikz
\begin{tikzpicture}
\begin{groupplot}[legend style={fill opacity={0.6}, text opacity={1}}, legend pos={north west}, xlabel={$N$}, xmode={log}, ymode={log}, legend cell align={left}, cycle list name={mark list*}, group style={group size={3 by 1}}]
    \nextgroupplot[]
    \addplot
        table[row sep={\\}]
        {
            \\
            128.0  0.006559886  \\
            256.0  0.010815244  \\
            1024.0  0.018016971  \\
            2048.0  0.021587044  \\
            8192.0  0.07260329  \\
            16384.0  0.165589974  \\
            65536.0  0.326611388  \\
            131072.0  1.863352492  \\
            262144.0  10.41134604  \\
            1.048576e6  67.839123656  \\
        }
        ;
    \addplot
        table[row sep={\\}]
        {
            \\
            128.0  0.007521395  \\
            256.0  0.010992711  \\
            1024.0  0.018268995  \\
            2048.0  0.024074927  \\
            8192.0  0.057630148  \\
            16384.0  0.235928896  \\
            65536.0  0.49707117  \\
            131072.0  3.067578127  \\
            262144.0  16.791556921  \\
            1.048576e6  107.814356703  \\
        }
        ;
    \addplot
        table[row sep={\\}]
        {
            \\
            128.0  0.007567829  \\
            256.0  0.011993446  \\
            1024.0  0.016737899  \\
            2048.0  0.029736122  \\
            8192.0  0.074103907  \\
            16384.0  0.309268747  \\
            65536.0  0.656678735  \\
            131072.0  3.875361699  \\
            262144.0  24.087592542  \\
            1.048576e6  143.897491801  \\
        }
        ;
    \addplot
        table[row sep={\\}]
        {
            \\
            128.0  0.001298928  \\
            256.0  0.001364934  \\
            1024.0  0.001979252  \\
            2048.0  0.047040543  \\
            8192.0  0.023151072  \\
            16384.0  0.55269533  \\
            65536.0  0.28193101  \\
            131072.0  1.648317367  \\
            262144.0  9.179458347  \\
            1.048576e6  61.264547542  \\
        }
        ;
    \addplot
        table[row sep={\\}]
        {
            \\
            128.0  0.001425117  \\
            256.0  0.001473176  \\
            1024.0  0.002149916  \\
            2048.0  0.007189636  \\
            8192.0  0.034574971  \\
            16384.0  0.19921386  \\
            65536.0  0.513558377  \\
            131072.0  2.777212148  \\
            262144.0  15.189015064  \\
            1.048576e6  104.922311922  \\
        }
        ;
    \addplot
        table[row sep={\\}]
        {
            \\
            128.0  0.001023883  \\
            256.0  0.00291371  \\
            1024.0  0.00252236  \\
            2048.0  0.009656688  \\
            8192.0  0.045518704  \\
            16384.0  0.269014554  \\
            65536.0  0.864711771  \\
            131072.0  3.651483244  \\
            262144.0  22.478449118  \\
            1.048576e6  154.937746686  \\
        }
        ;
    \addplot
        table[row sep={\\}]
        {
            \\
            128.0  0.001339155  \\
            256.0  0.001940379  \\
            1024.0  0.003026792  \\
            2048.0  0.007782118  \\
            8192.0  0.034265335  \\
            16384.0  0.19576301  \\
            65536.0  0.472700758  \\
            131072.0  2.737454565  \\
            262144.0  15.071435837  \\
            1.048576e6  91.971108316  \\
        }
        ;
    \addplot
        table[row sep={\\}]
        {
            \\
            128.0  0.00131654  \\
            256.0  0.002289365  \\
            1024.0  0.002745027  \\
            2048.0  0.008857955  \\
            8192.0  0.045817817  \\
            16384.0  0.267981912  \\
            65536.0  0.611280159  \\
            131072.0  3.591692904  \\
            262144.0  20.165145417  \\
            1.048576e6  128.945657511  \\
        }
        ;
    \addplot
        table[row sep={\\}]
        {
            \\
            128.0  0.001497559  \\
            256.0  0.002818636  \\
            1024.0  0.003369996  \\
            2048.0  0.010706927  \\
            8192.0  0.056065697  \\
            16384.0  0.344293601  \\
            65536.0  0.770968906  \\
            131072.0  5.0012135  \\
            262144.0  28.328314468  \\
            1.048576e6  281.401492409  \\
        }
        ;
    \addplot[style={black,dashed}]
        table[row sep={\\}]
        {
            \\
            128.0  0.00128  \\
            256.0  0.00256  \\
            1024.0  0.01024  \\
            2048.0  0.02048  \\
            8192.0  0.08192  \\
            16384.0  0.16384  \\
            65536.0  0.65536  \\
            131072.0  1.31072  \\
            262144.0  2.62144  \\
            1.048576e6  10.48576  \\
        }
        ;
    \legend{{},{},{},{},{},{},{},{},{},{$\mathcal O(N^{1.00})$}}
    \addplot[style={black,dotted}]
        table[row sep={\\}]
        {
            \\
            128.0  0.01641736  \\
            256.0  0.018352284  \\
            1024.0  0.027898474  \\
            2048.0  0.026889006  \\
            8192.0  0.036481972  \\
            16384.0  0.067785835  \\
            65536.0  0.133567418  \\
            131072.0  0.503214514  \\
            262144.0  1.765784313  \\
            1.048576e6  8.670407272  \\
        }
        ;
    \addplot[style={black,dotted}]
        table[row sep={\\}]
        {
            \\
            128.0  0.03368338  \\
            256.0  0.034692052  \\
            1024.0  0.045759821  \\
            2048.0  0.045524327  \\
            8192.0  0.054495011  \\
            16384.0  0.116403498  \\
            65536.0  0.171331663  \\
            131072.0  0.589669178  \\
            262144.0  2.401287696  \\
            1.048576e6  11.986890966  \\
        }
        ;
    \addplot[style={black,dotted}]
        table[row sep={\\}]
        {
            \\
            128.0  0.05225485  \\
            256.0  0.059210613  \\
            1024.0  0.057853528  \\
            2048.0  0.067934309  \\
            8192.0  0.099882031  \\
            16384.0  0.121874484  \\
            65536.0  0.198053095  \\
            131072.0  0.642503399  \\
            262144.0  2.697758827  \\
            1.048576e6  11.857690128  \\
        }
        ;
    \addplot[style={black,dotted}]
        table[row sep={\\}]
        {
            \\
            128.0  0.011101017  \\
            256.0  0.010014863  \\
            1024.0  0.010530533  \\
            2048.0  0.011717069  \\
            8192.0  0.016730555  \\
            16384.0  0.04527034  \\
            65536.0  0.11996116  \\
            131072.0  0.341034974  \\
            262144.0  1.706501353  \\
            1.048576e6  8.622987773  \\
        }
        ;
    \addplot[style={black,dotted}]
        table[row sep={\\}]
        {
            \\
            128.0  0.010013508  \\
            256.0  0.010552498  \\
            1024.0  0.010566111  \\
            2048.0  0.011877867  \\
            8192.0  0.017115778  \\
            16384.0  0.320062267  \\
            65536.0  0.098502603  \\
            131072.0  0.443886008  \\
            262144.0  2.054343688  \\
            1.048576e6  10.366463651  \\
        }
        ;
    \addplot[style={black,dotted}]
        table[row sep={\\}]
        {
            \\
            128.0  0.01009407  \\
            256.0  0.011613476  \\
            1024.0  0.010469314  \\
            2048.0  0.011667276  \\
            8192.0  0.018048597  \\
            16384.0  0.054719619  \\
            65536.0  0.4428434  \\
            131072.0  0.468008155  \\
            262144.0  2.134378046  \\
            1.048576e6  9.601493537  \\
        }
        ;
    \addplot[style={black,dotted}]
        table[row sep={\\}]
        {
            \\
            128.0  0.019293235  \\
            256.0  0.019053361  \\
            1024.0  0.018807343  \\
            2048.0  0.020925962  \\
            8192.0  0.026537828  \\
            16384.0  0.079413595  \\
            65536.0  0.125417799  \\
            131072.0  0.448438626  \\
            262144.0  2.139283054  \\
            1.048576e6  10.851535373  \\
        }
        ;
    \addplot[style={black,dotted}]
        table[row sep={\\}]
        {
            \\
            128.0  0.019442165  \\
            256.0  0.020055778  \\
            1024.0  0.01967791  \\
            2048.0  0.021374131  \\
            8192.0  0.025653927  \\
            16384.0  0.0788598  \\
            65536.0  0.107473955  \\
            131072.0  0.481610282  \\
            262144.0  2.228423843  \\
            1.048576e6  10.9973603  \\
        }
        ;
    \addplot[style={black,dotted}]
        table[row sep={\\}]
        {
            \\
            128.0  0.020747355  \\
            256.0  0.023889096  \\
            1024.0  0.020442163  \\
            2048.0  0.021369848  \\
            8192.0  0.031200173  \\
            16384.0  0.068330119  \\
            65536.0  0.12474847  \\
            131072.0  0.536885218  \\
            262144.0  2.448845502  \\
            1.048576e6  12.807178525  \\
        }
        ;
    \nextgroupplot[]
    \addplot
        table[row sep={\\}]
        {
            \\
            64.0  0.073788436  \\
            256.0  0.016458177  \\
            1024.0  0.175881052  \\
            4096.0  1.634061625  \\
            16384.0  14.493813903  \\
            65536.0  171.49888932  \\
            262144.0  0.0  \\
            1.048576e6  0.0  \\
        }
        ;
    \addplot
        table[row sep={\\}]
        {
            \\
            64.0  0.003537906  \\
            256.0  0.043132394  \\
            1024.0  0.294348281  \\
            4096.0  3.843025294  \\
            16384.0  41.701734845  \\
            65536.0  0.0  \\
            262144.0  0.0  \\
            1.048576e6  0.0  \\
        }
        ;
    \addplot
        table[row sep={\\}]
        {
            \\
            64.0  0.003701912  \\
            256.0  0.039124683  \\
            1024.0  0.521440069  \\
            4096.0  5.385231175  \\
            16384.0  86.75837726  \\
            65536.0  0.0  \\
            262144.0  0.0  \\
            1.048576e6  0.0  \\
        }
        ;
    \addplot
        table[row sep={\\}]
        {
            \\
            64.0  0.002207973  \\
            256.0  0.017563084  \\
            1024.0  0.176119847  \\
            4096.0  1.326755743  \\
            16384.0  15.870128747  \\
            65536.0  0.0  \\
            262144.0  0.0  \\
            1.048576e6  0.0  \\
        }
        ;
    \addplot
        table[row sep={\\}]
        {
            \\
            64.0  0.002382077  \\
            256.0  0.024442399  \\
            1024.0  0.276461466  \\
            4096.0  3.896135244  \\
            16384.0  41.997363357  \\
            65536.0  0.0  \\
            262144.0  0.0  \\
            1.048576e6  0.0  \\
        }
        ;
    \addplot
        table[row sep={\\}]
        {
            \\
            64.0  0.002406884  \\
            256.0  0.037259319  \\
            1024.0  0.472249103  \\
            4096.0  5.419083452  \\
            16384.0  67.250483017  \\
            65536.0  0.0  \\
            262144.0  0.0  \\
            1.048576e6  0.0  \\
        }
        ;
    \addplot
        table[row sep={\\}]
        {
            \\
            64.0  0.002412947  \\
            256.0  0.029382742  \\
            1024.0  0.260025046  \\
            4096.0  3.346172675  \\
            16384.0  43.93722154  \\
            65536.0  0.0  \\
            262144.0  0.0  \\
            1.048576e6  0.0  \\
        }
        ;
    \addplot
        table[row sep={\\}]
        {
            \\
            64.0  0.002446731  \\
            256.0  0.028957977  \\
            1024.0  0.351941131  \\
            4096.0  5.25629096  \\
            16384.0  87.891955777  \\
            65536.0  0.0  \\
            262144.0  0.0  \\
            1.048576e6  0.0  \\
        }
        ;
    \addplot
        table[row sep={\\}]
        {
            \\
            64.0  0.00238433  \\
            256.0  0.040601137  \\
            1024.0  0.54819183  \\
            4096.0  6.398400552  \\
            16384.0  88.85393391  \\
            65536.0  0.0  \\
            262144.0  0.0  \\
            1.048576e6  0.0  \\
        }
        ;
    \addplot[style={black,dashed}]
        table[row sep={\\}]
        {
            \\
            64.0  0.0064  \\
            256.0  0.0256  \\
            1024.0  0.1024  \\
            4096.0  0.4096  \\
            16384.0  1.6384  \\
            65536.0  6.5536  \\
            262144.0  26.2144  \\
        }
        ;
    \legend{{},{},{},{},{},{},{},{},{},{$\mathcal O(N^{1.00})$}}
    \addplot[style={black,dotted}]
        table[row sep={\\}]
        {
            \\
            64.0  0.069172184  \\
            256.0  0.056223112  \\
            1024.0  0.077659512  \\
            4096.0  0.200994781  \\
            16384.0  1.061347391  \\
            65536.0  3.583179499  \\
            262144.0  15.437803005  \\
            1.048576e6  0.0  \\
        }
        ;
    \addplot[style={black,dotted}]
        table[row sep={\\}]
        {
            \\
            64.0  0.028147923  \\
            256.0  0.066569156  \\
            1024.0  0.207162927  \\
            4096.0  1.10022845  \\
            16384.0  5.120571052  \\
            65536.0  19.34280771  \\
            262144.0  85.392188779  \\
            1.048576e6  0.0  \\
        }
        ;
    \addplot[style={black,dotted}]
        table[row sep={\\}]
        {
            \\
            64.0  0.040967134  \\
            256.0  0.105294947  \\
            1024.0  0.662289958  \\
            4096.0  2.73809348  \\
            16384.0  14.73883701  \\
            65536.0  69.205879867  \\
            262144.0  206.817897228  \\
            1.048576e6  0.0  \\
        }
        ;
    \addplot[style={black,dotted}]
        table[row sep={\\}]
        {
            \\
            64.0  0.008370861  \\
            256.0  0.029689664  \\
            1024.0  0.055170232  \\
            4096.0  0.20068383  \\
            16384.0  1.014020454  \\
            65536.0  3.775359867  \\
            262144.0  14.159454617  \\
            1.048576e6  0.0  \\
        }
        ;
    \addplot[style={black,dotted}]
        table[row sep={\\}]
        {
            \\
            64.0  0.007768416  \\
            256.0  0.03179759  \\
            1024.0  0.08131809  \\
            4096.0  0.362737308  \\
            16384.0  1.848973761  \\
            65536.0  9.556342632  \\
            262144.0  44.954903727  \\
            1.048576e6  0.0  \\
        }
        ;
    \addplot[style={black,dotted}]
        table[row sep={\\}]
        {
            \\
            64.0  0.00799658  \\
            256.0  0.033874264  \\
            1024.0  0.09803432  \\
            4096.0  0.574383162  \\
            16384.0  3.624448367  \\
            65536.0  18.702717658  \\
            262144.0  92.555255567  \\
            1.048576e6  0.0  \\
        }
        ;
    \addplot[style={black,dotted}]
        table[row sep={\\}]
        {
            \\
            64.0  0.013459952  \\
            256.0  0.041978494  \\
            1024.0  0.129496866  \\
            4096.0  0.847446714  \\
            16384.0  3.694096035  \\
            65536.0  16.278107702  \\
            262144.0  75.599888947  \\
            1.048576e6  0.0  \\
        }
        ;
    \addplot[style={black,dotted}]
        table[row sep={\\}]
        {
            \\
            64.0  0.013236752  \\
            256.0  0.044441684  \\
            1024.0  0.142431565  \\
            4096.0  1.199708221  \\
            16384.0  7.426851056  \\
            65536.0  32.178535774  \\
            262144.0  160.518755883  \\
            1.048576e6  0.0  \\
        }
        ;
    \addplot[style={black,dotted}]
        table[row sep={\\}]
        {
            \\
            64.0  0.03234988  \\
            256.0  0.047645558  \\
            1024.0  0.217746727  \\
            4096.0  1.507180827  \\
            16384.0  10.005199867  \\
            65536.0  53.306766365  \\
            262144.0  266.500990847  \\
            1.048576e6  0.0  \\
        }
        ;
    \nextgroupplot[]
    \addplot
        table[row sep={\\}]
        {
            \\
            64.0  0.048795897  \\
            512.0  0.118993486  \\
            4096.0  8.779235549  \\
            32768.0  940.6559111  \\
        }
        ;
    \addplot
        table[row sep={\\}]
        {
            \\
            64.0  0.007740059  \\
            512.0  0.156440213  \\
            4096.0  13.550588782  \\
            32768.0  2384.096208943  \\
        }
        ;
    \addplot
        table[row sep={\\}]
        {
            \\
            64.0  0.006007689  \\
            512.0  0.156905906  \\
            4096.0  16.887705298  \\
            32768.0  1580.715220099  \\
        }
        ;
    \addplot
        table[row sep={\\}]
        {
            \\
            64.0  0.002179503  \\
            512.0  0.073188388  \\
            4096.0  9.807359329  \\
            32768.0  690.130884584  \\
        }
        ;
    \addplot
        table[row sep={\\}]
        {
            \\
            64.0  0.002434419  \\
            512.0  0.112182185  \\
            4096.0  14.836551047  \\
            32768.0  1995.150671764  \\
        }
        ;
    \addplot
        table[row sep={\\}]
        {
            \\
            64.0  0.002682028  \\
            512.0  0.110233997  \\
            4096.0  19.895062377  \\
            32768.0  2384.628860119  \\
        }
        ;
    \addplot
        table[row sep={\\}]
        {
            \\
            64.0  0.002766554  \\
            512.0  0.142957262  \\
            4096.0  15.000250406  \\
            32768.0  1817.898344863  \\
        }
        ;
    \addplot
        table[row sep={\\}]
        {
            \\
            64.0  0.00252855  \\
            512.0  0.150239627  \\
            4096.0  20.441463889  \\
            32768.0  2551.675423327  \\
        }
        ;
    \addplot
        table[row sep={\\}]
        {
            \\
            64.0  0.002897048  \\
            512.0  0.179874395  \\
            4096.0  23.707935529  \\
            32768.0  3476.654942793  \\
        }
        ;
    \addplot[style={black,dashed}]
        table[row sep={\\}]
        {
            \\
            64.0  0.064  \\
            512.0  0.512  \\
            4096.0  4.096  \\
            32768.0  32.768  \\
        }
        ;
    \legend{{},{},{},{},{},{},{},{},{},{$\mathcal O(N^{1.00})$}}
    \addplot[style={black,dotted}]
        table[row sep={\\}]
        {
            \\
            64.0  0.083294029  \\
            512.0  0.189600045  \\
            4096.0  8.40313058  \\
            32768.0  359.846001586  \\
        }
        ;
    \addplot[style={black,dotted}]
        table[row sep={\\}]
        {
            \\
            64.0  0.053029292  \\
            512.0  0.535061235  \\
            4096.0  28.537315942  \\
            32768.0  1381.035493565  \\
        }
        ;
    \addplot[style={black,dotted}]
        table[row sep={\\}]
        {
            \\
            64.0  0.071527756  \\
            512.0  1.015401753  \\
            4096.0  77.708344632  \\
            32768.0  2825.806047609  \\
        }
        ;
    \addplot[style={black,dotted}]
        table[row sep={\\}]
        {
            \\
            64.0  0.027427256  \\
            512.0  0.164480158  \\
            4096.0  7.478236207  \\
            32768.0  289.450580528  \\
        }
        ;
    \addplot[style={black,dotted}]
        table[row sep={\\}]
        {
            \\
            64.0  0.027178229  \\
            512.0  0.23056075  \\
            4096.0  12.250088554  \\
            32768.0  589.95178702  \\
        }
        ;
    \addplot[style={black,dotted}]
        table[row sep={\\}]
        {
            \\
            64.0  0.028104655  \\
            512.0  0.232629574  \\
            4096.0  15.768962018  \\
            32768.0  780.803343876  \\
        }
        ;
    \addplot[style={black,dotted}]
        table[row sep={\\}]
        {
            \\
            64.0  0.039465131  \\
            512.0  0.376563042  \\
            4096.0  21.178938308  \\
            32768.0  1889.255021745  \\
        }
        ;
    \addplot[style={black,dotted}]
        table[row sep={\\}]
        {
            \\
            64.0  0.043427167  \\
            512.0  0.707447736  \\
            4096.0  27.223738699  \\
            32768.0  1509.995750552  \\
        }
        ;
    \addplot[style={black,dotted}]
        table[row sep={\\}]
        {
            \\
            64.0  0.062854337  \\
            512.0  0.456949504  \\
            4096.0  32.519593768  \\
            32768.0  1834.989058737  \\
        }
        ;
\end{groupplot}
\end{tikzpicture}

%% file: img/smallresidual1d-3d_q4.tikz
\begin{tikzpicture}
\begin{groupplot}[xmode={log}, ymode={log}, legend cell align={left}, legend style={fill opacity={0.6}, text opacity={1}}, legend pos={south west}, width={0.5\textwidth}, height={0.4\textwidth}, ylabel={Residual}, cycle list name={mark list*}, group style={x descriptions at={edge bottom}, ylabels at={edge left}, vertical sep={1em}, group size={3 by 1}}]
    \nextgroupplot[legend to name={test-1.353689255797669}]
    \addplot
        table[row sep={\\}]
        {
            \\
            64.0  0.046487224092280836  \\
            256.0  0.03603680237782705  \\
            1024.0  0.020974441308179412  \\
            4096.0  0.011148565477848118  \\
            16384.0  0.005766941473513112  \\
            65536.0  0.002939118361841445  \\
            262144.0  0.0  \\
            1.048576e6  0.0  \\
        }
        ;
    \addlegendentry {$\texttt{db2}$}
    \addplot
        table[row sep={\\}]
        {
            \\
            64.0  0.06754613537666118  \\
            256.0  0.0025127317815469147  \\
            1024.0  0.0009915132340239182  \\
            4096.0  0.00029939468344496986  \\
            16384.0  8.240241711569973e-5  \\
            65536.0  2.1616361938933608e-5  \\
            262144.0  0.0  \\
            1.048576e6  0.0  \\
        }
        ;
    \addlegendentry {$\texttt{db3}$}
    \addplot
        table[row sep={\\}]
        {
            \\
            64.0  0.08586134687714009  \\
            256.0  0.00014382583579377665  \\
            1024.0  1.120496965866664e-5  \\
            4096.0  2.1560323725818357e-6  \\
            16384.0  4.1060915062359384e-7  \\
            65536.0  0.0  \\
            262144.0  0.0  \\
            1.048576e6  0.0  \\
        }
        ;
    \addlegendentry {$\texttt{db4}$}
    \addplot
        table[row sep={\\}]
        {
            \\
            64.0  0.031004855393026942  \\
            256.0  2.4009881395641202e-5  \\
            1024.0  6.119908971464633e-6  \\
            4096.0  1.5428583659497483e-6  \\
            16384.0  3.868351439973905e-7  \\
            65536.0  0.0  \\
            262144.0  0.0  \\
            1.048576e6  0.0  \\
        }
        ;
    \addlegendentry {$\texttt{cdf31}$}
    \addplot
        table[row sep={\\}]
        {
            \\
            64.0  0.031004855393024233  \\
            256.0  2.4009881395993832e-5  \\
            1024.0  6.119908972043627e-6  \\
            4096.0  1.542858371151035e-6  \\
            16384.0  3.868351423114853e-7  \\
            65536.0  0.0  \\
            262144.0  0.0  \\
            1.048576e6  0.0  \\
        }
        ;
    \addlegendentry {$\texttt{cdf33}$}
    \addplot
        table[row sep={\\}]
        {
            \\
            64.0  0.031004855393025603  \\
            256.0  2.4009881395257666e-5  \\
            1024.0  6.119908968873037e-6  \\
            4096.0  1.5428583654967953e-6  \\
            16384.0  3.8683514395411305e-7  \\
            65536.0  0.0  \\
            262144.0  0.0  \\
            1.048576e6  0.0  \\
        }
        ;
    \addlegendentry {$\texttt{cdf35}$}
    \addplot
        table[row sep={\\}]
        {
            \\
            64.0  0.05179424571729861  \\
            256.0  1.676032429867409e-6  \\
            1024.0  3.118176094512769e-7  \\
            4096.0  4.518268502703894e-8  \\
            16384.0  6.107218731117894e-9  \\
            65536.0  0.0  \\
            262144.0  0.0  \\
            1.048576e6  0.0  \\
        }
        ;
    \addlegendentry {$\texttt{cdf42}$}
    \addplot
        table[row sep={\\}]
        {
            \\
            64.0  0.051794245717548766  \\
            256.0  1.6760324978931215e-6  \\
            1024.0  3.118176991541471e-7  \\
            4096.0  4.5183049607342403e-8  \\
            16384.0  6.1073549183426105e-9  \\
            65536.0  0.0  \\
            262144.0  0.0  \\
            1.048576e6  0.0  \\
        }
        ;
    \addlegendentry {$\texttt{cdf44}$}
    \addplot
        table[row sep={\\}]
        {
            \\
            64.0  0.0517942457175684  \\
            256.0  1.6760323586603268e-6  \\
            1024.0  3.118175607421164e-7  \\
            4096.0  4.518265100123969e-8  \\
            16384.0  6.107358909784225e-9  \\
            65536.0  0.0  \\
            262144.0  0.0  \\
            1.048576e6  0.0  \\
        }
        ;
    \addlegendentry {$\texttt{cdf46}$}
    \nextgroupplot[legend to name={test0.35117786982356297}]
    \addplot
        table[row sep={\\}]
        {
            \\
            64.0  0.0565586969328909  \\
            256.0  0.0675569113533023  \\
            1024.0  0.09013299986406081  \\
            4096.0  0.1914339663285852  \\
            16384.0  537.1321856046962  \\
            65536.0  45.218597832483155  \\
            262144.0  0.0  \\
            1.048576e6  0.0  \\
        }
        ;
    \addlegendentry {$\texttt{db2}$}
    \addplot
        table[row sep={\\}]
        {
            \\
            64.0  0.06758892046940022  \\
            256.0  0.010538338762118784  \\
            1024.0  0.03771071097368857  \\
            4096.0  0.27169784869266744  \\
            16384.0  2.1884031965427906  \\
            65536.0  0.09691925857364207  \\
            262144.0  0.0  \\
            1.048576e6  0.0  \\
        }
        ;
    \addlegendentry {$\texttt{db3}$}
    \addplot
        table[row sep={\\}]
        {
            \\
            64.0  3.4017547251309434e12  \\
            256.0  5.384311455877342e9  \\
            1024.0  2.647635483113165e10  \\
            4096.0  2.770029876116731e10  \\
            16384.0  3.1523514384037308e10  \\
            65536.0  0.0  \\
            262144.0  0.0  \\
            1.048576e6  0.0  \\
        }
        ;
    \addlegendentry {$\texttt{db4}$}
    \addplot
        table[row sep={\\}]
        {
            \\
            64.0  0.031004855393022567  \\
            256.0  2.400988139750792e-5  \\
            1024.0  7.518958448119615e-5  \\
            4096.0  0.00010279831891480127  \\
            16384.0  0.0005706269967237781  \\
            65536.0  0.0  \\
            262144.0  0.0  \\
            1.048576e6  0.0  \\
        }
        ;
    \addlegendentry {$\texttt{cdf31}$}
    \addplot
        table[row sep={\\}]
        {
            \\
            64.0  0.031004855393024226  \\
            256.0  2.4009881414010616e-5  \\
            1024.0  6.314863723709179e-5  \\
            4096.0  0.0007183544714184973  \\
            16384.0  0.00012739370548240492  \\
            65536.0  0.0  \\
            262144.0  0.0  \\
            1.048576e6  0.0  \\
        }
        ;
    \addlegendentry {$\texttt{cdf33}$}
    \addplot
        table[row sep={\\}]
        {
            \\
            64.0  0.031004855393014105  \\
            256.0  2.4009881419798744e-5  \\
            1024.0  0.0003614664419997655  \\
            4096.0  0.009132485589515137  \\
            16384.0  0.005598970024863244  \\
            65536.0  0.0  \\
            262144.0  0.0  \\
            1.048576e6  0.0  \\
        }
        ;
    \addlegendentry {$\texttt{cdf35}$}
    \addplot
        table[row sep={\\}]
        {
            \\
            64.0  0.0591741250237381  \\
            256.0  3.0084229378966627e-6  \\
            1024.0  1.953118587159765e-5  \\
            4096.0  0.00010759226320696445  \\
            16384.0  6.900238750956133e-5  \\
            65536.0  0.0  \\
            262144.0  0.0  \\
            1.048576e6  0.0  \\
        }
        ;
    \addlegendentry {$\texttt{cdf42}$}
    \addplot
        table[row sep={\\}]
        {
            \\
            64.0  0.06762122308631785  \\
            256.0  2.945234481698714e-6  \\
            1024.0  4.4843470309636754e-6  \\
            4096.0  8.89253195639968e-5  \\
            16384.0  1.446308241953466e-5  \\
            65536.0  0.0  \\
            262144.0  0.0  \\
            1.048576e6  0.0  \\
        }
        ;
    \addlegendentry {$\texttt{cdf44}$}
    \addplot
        table[row sep={\\}]
        {
            \\
            64.0  0.07061632122919723  \\
            256.0  2.332532004930083e-6  \\
            1024.0  1.0541985657433788e-6  \\
            4096.0  0.00033250464257330903  \\
            16384.0  3.254927156412794e-5  \\
            65536.0  0.0  \\
            262144.0  0.0  \\
            1.048576e6  0.0  \\
        }
        ;
    \addlegendentry {$\texttt{cdf46}$}
    \nextgroupplot[legend to name={test-1.5502161325289747}]
    \addplot
        table[row sep={\\}]
        {
            \\
            64.0  0.025770956639057808  \\
            256.0  0.014772710578410612  \\
            1024.0  0.00788805668488226  \\
            4096.0  0.0040443909268980584  \\
            16384.0  0.002046174284823846  \\
            65536.0  0.001031457929327939  \\
            262144.0  0.0  \\
            1.048576e6  0.0  \\
        }
        ;
    \addlegendentry {$\texttt{db2}$}
    \addplot
        table[row sep={\\}]
        {
            \\
            64.0  0.062327753048890946  \\
            256.0  0.00023307945729319506  \\
            1024.0  6.598363326664734e-5  \\
            4096.0  1.8612231426029078e-5  \\
            16384.0  4.391455081210662e-6  \\
            65536.0  1.1184806163211453e-6  \\
            262144.0  0.0  \\
            1.048576e6  0.0  \\
        }
        ;
    \addlegendentry {$\texttt{db3}$}
    \addplot
        table[row sep={\\}]
        {
            \\
            64.0  0.08529502543416703  \\
            256.0  7.196811487789718e-5  \\
            1024.0  5.082333017678014e-7  \\
            4096.0  8.652586106306157e-8  \\
            16384.0  9.242374828000262e-9  \\
            65536.0  0.0  \\
            262144.0  0.0  \\
            1.048576e6  0.0  \\
        }
        ;
    \addlegendentry {$\texttt{db4}$}
    \addplot
        table[row sep={\\}]
        {
            \\
            64.0  0.027542771235698656  \\
            256.0  2.3937191331230974e-5  \\
            1024.0  6.113187127666321e-6  \\
            4096.0  1.5417306998500925e-6  \\
            16384.0  3.869041353648323e-7  \\
            65536.0  0.0  \\
            262144.0  0.0  \\
            1.048576e6  0.0  \\
        }
        ;
    \addlegendentry {$\texttt{cdf31}$}
    \addplot
        table[row sep={\\}]
        {
            \\
            64.0  0.027542771235697847  \\
            256.0  2.3937299332767325e-5  \\
            1024.0  6.121560885530305e-6  \\
            4096.0  1.543876907545998e-6  \\
            16384.0  3.870154546919204e-7  \\
            65536.0  0.0  \\
            262144.0  0.0  \\
            1.048576e6  0.0  \\
        }
        ;
    \addlegendentry {$\texttt{cdf33}$}
    \addplot
        table[row sep={\\}]
        {
            \\
            64.0  0.02754277123569811  \\
            256.0  2.393719133196891e-5  \\
            1024.0  6.183113750319968e-6  \\
            4096.0  1.5422957528001244e-6  \\
            16384.0  3.894552836210371e-7  \\
            65536.0  0.0  \\
            262144.0  0.0  \\
            1.048576e6  0.0  \\
        }
        ;
    \addlegendentry {$\texttt{cdf35}$}
    \addplot
        table[row sep={\\}]
        {
            \\
            64.0  0.04565684296829485  \\
            256.0  1.772591921794951e-7  \\
            1024.0  2.298270875270557e-8  \\
            4096.0  2.9255027416254398e-9  \\
            16384.0  3.6964958249083904e-10  \\
            65536.0  0.0  \\
            262144.0  0.0  \\
            1.048576e6  0.0  \\
        }
        ;
    \addlegendentry {$\texttt{cdf42}$}
    \addplot
        table[row sep={\\}]
        {
            \\
            64.0  0.04565684296829521  \\
            256.0  1.7726019877599978e-7  \\
            1024.0  2.3026613284410947e-8  \\
            4096.0  2.9828758316724957e-9  \\
            16384.0  3.7072362249916883e-10  \\
            65536.0  0.0  \\
            262144.0  0.0  \\
            1.048576e6  0.0  \\
        }
        ;
    \addlegendentry {$\texttt{cdf44}$}
    \addplot
        table[row sep={\\}]
        {
            \\
            64.0  0.04565684296829431  \\
            256.0  1.773071502668083e-7  \\
            1024.0  2.3049528633413223e-8  \\
            4096.0  2.9344457283516343e-9  \\
            16384.0  6.665612822133378e-10  \\
            65536.0  0.0  \\
            262144.0  0.0  \\
            1.048576e6  0.0  \\
        }
        ;
    \addlegendentry {$\texttt{cdf46}$}
\end{groupplot}
\end{tikzpicture}

%% file: img/smallcoefnorm1d-3d_q4.tikz
\begin{tikzpicture}
\begin{groupplot}[xmode={log}, ymode={log}, legend cell align={left}, legend style={fill opacity={0.6}, text opacity={1}}, legend pos={south west}, width={0.5\textwidth}, height={0.4\textwidth}, xlabel={$N$}, ylabel={$\|x\|$}, cycle list name={mark list*}, group style={x descriptions at={edge bottom}, ylabels at={edge left}, vertical sep={1em}, group size={3 by 1}}]
    \nextgroupplot[legend to name={test0.14457046171983484}]
    \addplot
        table[row sep={\\}]
        {
            \\
            64.0  1.6632380120406487  \\
            256.0  1.47306179076303  \\
            1024.0  1.446013305290509  \\
            4096.0  1.4234732923954663  \\
            16384.0  1.3968649553717496  \\
            65536.0  1.3974687395662626  \\
            262144.0  0.0  \\
            1.048576e6  0.0  \\
        }
        ;
    \addlegendentry {$\texttt{db2}$}
    \addplot
        table[row sep={\\}]
        {
            \\
            64.0  4.3260173746171775  \\
            256.0  1870.5910469610865  \\
            1024.0  23.6293847197484  \\
            4096.0  2.4638483405344367  \\
            16384.0  1.444044855193355  \\
            65536.0  1.3593743001379497  \\
            262144.0  0.0  \\
            1.048576e6  0.0  \\
        }
        ;
    \addlegendentry {$\texttt{db3}$}
    \addplot
        table[row sep={\\}]
        {
            \\
            64.0  5.718345255479511  \\
            256.0  3.918249247757657  \\
            1024.0  1.3371034947409355  \\
            4096.0  1.330462957342077  \\
            16384.0  1.3136542225570773  \\
            65536.0  0.0  \\
            262144.0  0.0  \\
            1.048576e6  0.0  \\
        }
        ;
    \addlegendentry {$\texttt{db4}$}
    \addplot
        table[row sep={\\}]
        {
            \\
            64.0  2.001161208676086  \\
            256.0  1.6922702959151563  \\
            1024.0  1.5689259710485661  \\
            4096.0  1.5338944799343528  \\
            16384.0  1.5007937268202265  \\
            65536.0  0.0  \\
            262144.0  0.0  \\
            1.048576e6  0.0  \\
        }
        ;
    \addlegendentry {$\texttt{cdf31}$}
    \addplot
        table[row sep={\\}]
        {
            \\
            64.0  2.2165054581284696  \\
            256.0  2.05496053796269  \\
            1024.0  1.8709727786241364  \\
            4096.0  1.7805885229651452  \\
            16384.0  1.730211182235585  \\
            65536.0  0.0  \\
            262144.0  0.0  \\
            1.048576e6  0.0  \\
        }
        ;
    \addlegendentry {$\texttt{cdf33}$}
    \addplot
        table[row sep={\\}]
        {
            \\
            64.0  2.3317919367449242  \\
            256.0  2.207671894268745  \\
            1024.0  2.012719037973462  \\
            4096.0  1.904597156133155  \\
            16384.0  1.8484190852852047  \\
            65536.0  0.0  \\
            262144.0  0.0  \\
            1.048576e6  0.0  \\
        }
        ;
    \addlegendentry {$\texttt{cdf35}$}
    \addplot
        table[row sep={\\}]
        {
            \\
            64.0  89.85782537645935  \\
            256.0  1.3879574467414015  \\
            1024.0  1.3877189213687473  \\
            4096.0  1.4442222895925962  \\
            16384.0  1.6596044146289515  \\
            65536.0  0.0  \\
            262144.0  0.0  \\
            1.048576e6  0.0  \\
        }
        ;
    \addlegendentry {$\texttt{cdf42}$}
    \addplot
        table[row sep={\\}]
        {
            \\
            64.0  76.54570966030556  \\
            256.0  1.4096815955209903  \\
            1024.0  1.3838709690404492  \\
            4096.0  1.4098757490897567  \\
            16384.0  1.4387727026781703  \\
            65536.0  0.0  \\
            262144.0  0.0  \\
            1.048576e6  0.0  \\
        }
        ;
    \addlegendentry {$\texttt{cdf44}$}
    \addplot
        table[row sep={\\}]
        {
            \\
            64.0  71.73651322259794  \\
            256.0  1.426126623300203  \\
            1024.0  1.3861244034266857  \\
            4096.0  1.4086931725070817  \\
            16384.0  1.440841006044435  \\
            65536.0  0.0  \\
            262144.0  0.0  \\
            1.048576e6  0.0  \\
        }
        ;
    \addlegendentry {$\texttt{cdf46}$}
    \nextgroupplot[legend to name={test-1.1883991397908722}]
    \addplot
        table[row sep={\\}]
        {
            \\
            64.0  3.19838683724031e12  \\
            256.0  8.087183482913797e13  \\
            1024.0  7.753359122815472e12  \\
            4096.0  1.5526403554090814e13  \\
            16384.0  5.98951686339739e12  \\
            65536.0  2.2228213237768797e14  \\
            262144.0  0.0  \\
            1.048576e6  0.0  \\
        }
        ;
    \addlegendentry {$\texttt{db2}$}
    \addplot
        table[row sep={\\}]
        {
            \\
            64.0  44.503914654932686  \\
            256.0  2.0537380940765388e9  \\
            1024.0  6.763780799843103e11  \\
            4096.0  1.779230031024187e12  \\
            16384.0  2.4983281507729748e10  \\
            65536.0  6.744493126871735e9  \\
            262144.0  0.0  \\
            1.048576e6  0.0  \\
        }
        ;
    \addlegendentry {$\texttt{db3}$}
    \addplot
        table[row sep={\\}]
        {
            \\
            64.0  3.321288153711472e12  \\
            256.0  1.1194723838255553e11  \\
            1024.0  1.2854106567924881e11  \\
            4096.0  2.1353211735468164e11  \\
            16384.0  1.9555491300570486e12  \\
            65536.0  0.0  \\
            262144.0  0.0  \\
            1.048576e6  0.0  \\
        }
        ;
    \addlegendentry {$\texttt{db4}$}
    \addplot
        table[row sep={\\}]
        {
            \\
            64.0  42.34601428042735  \\
            256.0  234.89624107812423  \\
            1024.0  7.639659131858768e8  \\
            4096.0  9.979185039153547e9  \\
            16384.0  4.619199781784456e9  \\
            65536.0  0.0  \\
            262144.0  0.0  \\
            1.048576e6  0.0  \\
        }
        ;
    \addlegendentry {$\texttt{cdf31}$}
    \addplot
        table[row sep={\\}]
        {
            \\
            64.0  36.10113839061419  \\
            256.0  17.89519053049507  \\
            1024.0  1.2370478013155836e10  \\
            4096.0  2.0505395347017392e11  \\
            16384.0  9.823785211460483e10  \\
            65536.0  0.0  \\
            262144.0  0.0  \\
            1.048576e6  0.0  \\
        }
        ;
    \addlegendentry {$\texttt{cdf33}$}
    \addplot
        table[row sep={\\}]
        {
            \\
            64.0  34.20346454114026  \\
            256.0  15.773038461204242  \\
            1024.0  5.363031540144155e11  \\
            4096.0  1.2361212533583079e12  \\
            16384.0  2.3390030352162375e13  \\
            65536.0  0.0  \\
            262144.0  0.0  \\
            1.048576e6  0.0  \\
        }
        ;
    \addlegendentry {$\texttt{cdf35}$}
    \addplot
        table[row sep={\\}]
        {
            \\
            64.0  2.8611295166861225e10  \\
            256.0  4.835260597741448e7  \\
            1024.0  1.8289009760409813e9  \\
            4096.0  6.208431855034355e8  \\
            16384.0  3.2918714607043386e7  \\
            65536.0  0.0  \\
            262144.0  0.0  \\
            1.048576e6  0.0  \\
        }
        ;
    \addlegendentry {$\texttt{cdf42}$}
    \addplot
        table[row sep={\\}]
        {
            \\
            64.0  1.7173553347484816e10  \\
            256.0  3.521964932706465e7  \\
            1024.0  1.2086757416925874e9  \\
            4096.0  5.657865709810683e9  \\
            16384.0  1.0510943736062614e9  \\
            65536.0  0.0  \\
            262144.0  0.0  \\
            1.048576e6  0.0  \\
        }
        ;
    \addlegendentry {$\texttt{cdf44}$}
    \addplot
        table[row sep={\\}]
        {
            \\
            64.0  1.5274398839418066e10  \\
            256.0  5.07846574471037e7  \\
            1024.0  1.167179577083757e7  \\
            4096.0  1.6625282205944054e10  \\
            16384.0  1.1998995006844022e9  \\
            65536.0  0.0  \\
            262144.0  0.0  \\
            1.048576e6  0.0  \\
        }
        ;
    \addlegendentry {$\texttt{cdf46}$}
    \nextgroupplot[legend to name={test0.7829703463259083}]
    \addplot
        table[row sep={\\}]
        {
            \\
            64.0  2.1363357687630757  \\
            256.0  35.110743515911935  \\
            1024.0  3.3294318646464466e10  \\
            4096.0  2.5848389383749728e11  \\
            16384.0  4.8331794279768805e11  \\
            65536.0  1.535098118254908e12  \\
            262144.0  0.0  \\
            1.048576e6  0.0  \\
        }
        ;
    \addlegendentry {$\texttt{db2}$}
    \addplot
        table[row sep={\\}]
        {
            \\
            64.0  11.358892678801036  \\
            256.0  4.244576616161855e9  \\
            1024.0  5.556597487349728e12  \\
            4096.0  7.154145399728834e11  \\
            16384.0  5.843050394195754e11  \\
            65536.0  1.641687978926454e12  \\
            262144.0  0.0  \\
            1.048576e6  0.0  \\
        }
        ;
    \addlegendentry {$\texttt{db3}$}
    \addplot
        table[row sep={\\}]
        {
            \\
            64.0  4.525514521901371  \\
            256.0  2.282687611445861e6  \\
            1024.0  2.3462941836335425e9  \\
            4096.0  4.483462600759125e10  \\
            16384.0  2.0538056229878778e9  \\
            65536.0  0.0  \\
            262144.0  0.0  \\
            1.048576e6  0.0  \\
        }
        ;
    \addlegendentry {$\texttt{db4}$}
    \addplot
        table[row sep={\\}]
        {
            \\
            64.0  3.7170615258154522  \\
            256.0  2.480090389012968  \\
            1024.0  2.471530481418408e8  \\
            4096.0  1.2073491514046166e9  \\
            16384.0  2.2594902167306885e8  \\
            65536.0  0.0  \\
            262144.0  0.0  \\
            1.048576e6  0.0  \\
        }
        ;
    \addlegendentry {$\texttt{cdf31}$}
    \addplot
        table[row sep={\\}]
        {
            \\
            64.0  3.2744947044546024  \\
            256.0  1.5217049282831132e8  \\
            1024.0  1.0222167501751938e10  \\
            4096.0  2.8650012968286755e10  \\
            16384.0  4.141605269147019e10  \\
            65536.0  0.0  \\
            262144.0  0.0  \\
            1.048576e6  0.0  \\
        }
        ;
    \addlegendentry {$\texttt{cdf33}$}
    \addplot
        table[row sep={\\}]
        {
            \\
            64.0  3.148812743688125  \\
            256.0  4.018161548527892  \\
            1024.0  1.9839341579972024e12  \\
            4096.0  1.3869123058941425e11  \\
            16384.0  2.460888288662572e12  \\
            65536.0  0.0  \\
            262144.0  0.0  \\
            1.048576e6  0.0  \\
        }
        ;
    \addlegendentry {$\texttt{cdf35}$}
    \addplot
        table[row sep={\\}]
        {
            \\
            64.0  25.78662617295695  \\
            256.0  1.4717402480057685e6  \\
            1024.0  8.302364227935669e6  \\
            4096.0  4.546548419505867e7  \\
            16384.0  6.117612017323895e7  \\
            65536.0  0.0  \\
            262144.0  0.0  \\
            1.048576e6  0.0  \\
        }
        ;
    \addlegendentry {$\texttt{cdf42}$}
    \addplot
        table[row sep={\\}]
        {
            \\
            64.0  17.871170932322087  \\
            256.0  625321.9047771601  \\
            1024.0  6.159327830241685e7  \\
            4096.0  8.172131642306546e8  \\
            16384.0  7.589027309816933e7  \\
            65536.0  0.0  \\
            262144.0  0.0  \\
            1.048576e6  0.0  \\
        }
        ;
    \addlegendentry {$\texttt{cdf44}$}
    \addplot
        table[row sep={\\}]
        {
            \\
            64.0  15.148712616904433  \\
            256.0  2.494284436999029e6  \\
            1024.0  7.467136591012032e7  \\
            4096.0  1.331860520205291e10  \\
            16384.0  3.0236591368478703e9  \\
            65536.0  0.0  \\
            262144.0  0.0  \\
            1.048576e6  0.0  \\
        }
        ;
    \addlegendentry {$\texttt{cdf46}$}
\end{groupplot}
\end{tikzpicture}

%% file: img/1dadaptivity.tikz
\begin{tikzpicture}
\begin{groupplot}[group style={group size={2 by 1}, horizontal sep={4em}}]
    \nextgroupplot[width={.5\textwidth}, height={.25\textwidth}, no marks, cycle list name={color list}]
    \addplot+
        table[row sep={\\}]
        {
            \\
            0.0  1.0000000000000002  \\
            0.0078125  1.0078430972046153  \\
            0.015625  1.0157477085848383  \\
            0.0234375  1.023714316600496  \\
            0.03125  1.0317434074972263  \\
            0.0390625  1.0398354713343392  \\
            0.046875  1.047991002014727  \\
            0.0546875  1.056210497315011  \\
            0.0625  1.0644944589159229  \\
            0.0703125  1.072843392432926  \\
            0.078125  1.0812578074470727  \\
            0.0859375  1.0897382175361112  \\
            0.09375  1.0982851403058278  \\
            0.1015625  1.1068990974216435  \\
            0.109375  1.1155806146404508  \\
            0.1171875  1.1243302218427047  \\
            0.125  1.1331484530647646  \\
            0.1328125  1.142035846531488  \\
            0.140625  1.1509929446890825  \\
            0.1484375  1.160020294238215  \\
            0.15625  1.1691184461673776  \\
            0.1640625  1.1782879557865202  \\
            0.171875  1.1875293827609408  \\
            0.1796875  1.196843291145448  \\
            0.1875  1.2062302494187866  \\
            0.1953125  1.2156908305183365  \\
            0.203125  1.225225611875079  \\
            0.2109375  1.2348351754488454  \\
            0.21875  1.2445201077638322  \\
            0.2265625  1.2542809999444027  \\
            0.234375  1.2641184477511678  \\
            0.2421875  1.2740330516173441  \\
            0.25  1.284025416685406  \\
            0.2578125  1.2940961528440196  \\
            0.265625  1.3042458747652659  \\
            0.2734375  1.3144752019421588  \\
            0.28125  1.324784758726456  \\
            0.2890625  1.3351751743667684  \\
            0.296875  1.3456470830469631  \\
            0.3046875  1.3562011239248737  \\
            0.3125  1.3668379411713096  \\
            0.3203125  1.3775581840093776  \\
            0.328125  1.3883625067541006  \\
            0.3359375  1.3992515688523608  \\
            0.34375  1.410226034923145  \\
            0.3515625  1.421286574798111  \\
            0.359375  1.4324338635624723  \\
            0.3671875  1.4436685815962007  \\
            0.375  1.4549914146155545  \\
            0.3828125  1.4664030537149317  \\
            0.390625  1.47790419540905  \\
            0.3984375  1.4894955416754612  \\
            0.40625  1.5011777999973923  \\
            0.4140625  1.5129516834069334  \\
            0.421875  1.5248179105285533  \\
            0.4296875  1.5367772056229618  \\
            0.4375  1.5488302986313165  \\
            0.4453125  1.5609779252197735  \\
            0.453125  1.573220826824392  \\
            0.4609375  1.5855597506963837  \\
            0.46875  1.5979954499477271  \\
            0.4765625  1.6105286835971286  \\
            0.484375  1.6231602166163537  \\
            0.4921875  1.6358908199769153  \\
            0.5  1.64872127069713  \\
            0.5078125  1.661652351889546  \\
            0.515625  1.6746848528087386  \\
            0.5234375  1.687819568899484  \\
            0.53125  1.7010573018453075  \\
            0.5390625  1.714398859617418  \\
            0.546875  1.7278450565240209  \\
            0.5546875  1.7413967132600194  \\
            0.5625  1.7550546569571066  \\
            0.5703125  1.7688197212342507  \\
            0.578125  1.7826927462485722  \\
            0.5859375  1.7966745787466296  \\
            0.59375  1.810766072116094  \\
            0.6015625  1.823286344483741  \\
            0.609375  1.8138221579226226  \\
            0.6171875  1.8032096695776887  \\
            0.625  1.791215778333505  \\
            0.6328125  1.7778510738271747  \\
            0.640625  1.7631198881829657  \\
            0.6484375  1.7470222214008782  \\
            0.65625  1.7295580734809124  \\
            0.6640625  1.7107274444230853  \\
            0.671875  1.6905303342274043  \\
            0.6796875  1.6689667428938688  \\
            0.6875  1.6460366704225302  \\
            0.6953125  1.6217401168145167  \\
            0.703125  1.5960770820702885  \\
            0.7109375  1.5690475661898464  \\
            0.71875  1.5406515691731892  \\
            0.7265625  1.5108890910203094  \\
            0.734375  1.4797601317312017  \\
            0.7421875  1.4472646913058669  \\
            0.75  1.4134863120573657  \\
            0.7578125  1.3802629248730238  \\
            0.765625  1.348346410570384  \\
            0.7734375  1.3177367691494468  \\
            0.78125  1.288434000610212  \\
            0.7890625  1.2604381049526792  \\
            0.796875  1.2337490821768486  \\
            0.8046875  1.2083669322827202  \\
            0.8125  1.1842916552702627  \\
            0.8203125  1.161523251138785  \\
            0.828125  1.1400617198880043  \\
            0.8359375  1.1199070615179203  \\
            0.84375  1.1010592760285336  \\
            0.8515625  1.083518363419844  \\
            0.859375  1.067284323691851  \\
            0.8671875  1.0523571568445553  \\
            0.875  1.0387372161121269  \\
            0.8828125  1.0264322726463186  \\
            0.890625  1.0154455055546654  \\
            0.8984375  1.0057769148371678  \\
            0.90625  0.9974265004938236  \\
            0.9140625  0.9903942625245987  \\
            0.921875  0.9846802009294793  \\
            0.9296875  0.9802843157084653  \\
            0.9375  0.9772066068616608  \\
            0.9453125  0.9754470743913657  \\
            0.953125  0.9750057182985201  \\
            0.9609375  0.975882538583122  \\
            0.96875  0.9780775352451802  \\
            0.9765625  0.98159070828491  \\
            0.984375  0.986422057702401  \\
            0.9921875  0.992571583497667  \\
        }
        ;
    \addplot+
        table[row sep={\\}]
        {
            \\
            0.0  0.9999999999999991  \\
            0.0078125  1.0078430972074937  \\
            0.015625  1.015747708587635  \\
            0.0234375  1.0237143166032943  \\
            0.03125  1.0317434075000422  \\
            0.0390625  1.0398354713371765  \\
            0.046875  1.0479910020175869  \\
            0.0546875  1.056210497317894  \\
            0.0625  1.0644944589188274  \\
            0.0703125  1.0728433924358531  \\
            0.078125  1.081257807450022  \\
            0.0859375  1.089738217539083  \\
            0.09375  1.0982851403088234  \\
            0.1015625  1.1068990974246633  \\
            0.109375  1.1155806146434948  \\
            0.1171875  1.1243302218457727  \\
            0.125  1.1331484530678564  \\
            0.1328125  1.142035846534604  \\
            0.140625  1.150992944692223  \\
            0.1484375  1.1600202942413798  \\
            0.15625  1.1691184461705673  \\
            0.1640625  1.178287955789734  \\
            0.171875  1.1875293827641806  \\
            0.1796875  1.1968432911487128  \\
            0.1875  1.206230249422078  \\
            0.1953125  1.2156908305216536  \\
            0.203125  1.2252256118784226  \\
            0.2109375  1.2348351754522147  \\
            0.21875  1.2445201077672279  \\
            0.2265625  1.2542809999478255  \\
            0.234375  1.2641184477546166  \\
            0.2421875  1.2740330516208198  \\
            0.25  1.28402541668891  \\
            0.2578125  1.2940961528475508  \\
            0.265625  1.3042458747688244  \\
            0.2734375  1.3144752019457453  \\
            0.28125  1.3247847587300712  \\
            0.2890625  1.3351751743704123  \\
            0.296875  1.3456470830506353  \\
            0.3046875  1.356201123928574  \\
            0.3125  1.3668379411750398  \\
            0.3203125  1.3775581840131366  \\
            0.328125  1.3883625067578893  \\
            0.3359375  1.3992515688561793  \\
            0.34375  1.4102260349269933  \\
            0.3515625  1.4212865748019887  \\
            0.359375  1.4324338635663807  \\
            0.3671875  1.4436685816001398  \\
            0.375  1.4549914146195237  \\
            0.3828125  1.4664030537189325  \\
            0.390625  1.4779041954130823  \\
            0.3984375  1.4894955416795246  \\
            0.40625  1.5011778000014884  \\
            0.4140625  1.5129516834110621  \\
            0.421875  1.5248179105327138  \\
            0.4296875  1.536777205627155  \\
            0.4375  1.5488302986355422  \\
            0.4453125  1.560977925224032  \\
            0.453125  1.573220826828684  \\
            0.4609375  1.5855597507007095  \\
            0.46875  1.597995449952086  \\
            0.4765625  1.610528683601522  \\
            0.484375  1.6231602166207821  \\
            0.4921875  1.6358908199813782  \\
            0.5  1.648721270701628  \\
            0.5078125  1.6616523518940793  \\
            0.515625  1.674684852813308  \\
            0.5234375  1.6878195689040885  \\
            0.53125  1.7010573018499477  \\
            0.5390625  1.7143988596220958  \\
            0.546875  1.7278450565287353  \\
            0.5546875  1.741396713264771  \\
            0.5625  1.7550546569618954  \\
            0.5703125  1.7688197212390786  \\
            0.578125  1.782692746253447  \\
            0.5859375  1.7966745787515785  \\
            0.59375  1.810766072121263  \\
            0.6015625  1.7954681825137702  \\
            0.609375  1.5791786663151715  \\
            0.6171875  1.5447517281434033  \\
            0.625  1.5069717974622796  \\
            0.6328125  1.5123303772542895  \\
            0.640625  1.5258852889757075  \\
            0.6484375  1.5171841494333695  \\
            0.65625  1.4774560019461884  \\
            0.6640625  1.4296861429110372  \\
            0.671875  1.3832327985083057  \\
            0.6796875  1.3371092771346995  \\
            0.6875  1.2911081353010256  \\
            0.6953125  1.2495458406746773  \\
            0.703125  1.2141882209377703  \\
            0.7109375  1.1850352760903045  \\
            0.71875  1.1619834916123186  \\
            0.7265625  1.142755548064665  \\
            0.734375  1.1264198147676923  \\
            0.7421875  1.1129762917214014  \\
            0.75  1.1023913435953188  \\
            0.7578125  1.0939249931190371  \\
            0.765625  1.0872745223182987  \\
            0.7734375  1.0824399311931039  \\
            0.78125  1.0794212197434523  \\
            0.7890625  1.0782183879693443  \\
            0.796875  1.07883143587078  \\
            0.8046875  1.0812603634477589  \\
            0.8125  1.0854242796794538  \\
            0.8203125  1.089543582107658  \\
            0.828125  1.0928902515449237  \\
            0.8359375  1.0954642879912506  \\
            0.84375  1.0972656914466388  \\
            0.8515625  1.0982944619110884  \\
            0.859375  1.098550599384599  \\
            0.8671875  1.0980341038671715  \\
            0.875  1.0967600686200678  \\
            0.8828125  1.0950605453910753  \\
            0.890625  1.0930713735315614  \\
            0.8984375  1.0907925530415257  \\
            0.90625  1.0882091506010352  \\
            0.9140625  1.0849926331715571  \\
            0.921875  1.0810086008736928  \\
            0.9296875  1.0762570537074414  \\
            0.9375  1.0707679548537463  \\
            0.9453125  1.0652004942933573  \\
            0.953125  1.0597500292742825  \\
            0.9609375  1.0527817094259588  \\
            0.96875  1.0437787782769612  \\
            0.9765625  1.0357656598170655  \\
            0.984375  1.0237715808674408  \\
            0.9921875  1.0120112527924106  \\
        }
        ;
    \addplot+
        table[row sep={\\}]
        {
            \\
            0.0  1.0000000000000002  \\
            0.0078125  1.0078430972046153  \\
            0.015625  1.0157477085848383  \\
            0.0234375  1.0237143166004958  \\
            0.03125  1.0317434074972263  \\
            0.0390625  1.0398354713343392  \\
            0.046875  1.0479910020147272  \\
            0.0546875  1.0562104973150113  \\
            0.0625  1.064494458915923  \\
            0.0703125  1.0728433924329261  \\
            0.078125  1.0812578074470727  \\
            0.0859375  1.0897382175361108  \\
            0.09375  1.098285140305828  \\
            0.1015625  1.1068990974216444  \\
            0.109375  1.1155806146404517  \\
            0.1171875  1.1243302218427056  \\
            0.125  1.1331484530647655  \\
            0.1328125  1.1420358465314884  \\
            0.140625  1.1509929446890828  \\
            0.1484375  1.1600202942382154  \\
            0.15625  1.1691184461673778  \\
            0.1640625  1.1782879557865198  \\
            0.171875  1.18752938276094  \\
            0.1796875  1.1968432911454476  \\
            0.1875  1.2062302494187866  \\
            0.1953125  1.2156908305183356  \\
            0.203125  1.2252256118750786  \\
            0.2109375  1.2348351754488442  \\
            0.21875  1.244520107763831  \\
            0.2265625  1.2542809999444018  \\
            0.234375  1.2641184477511669  \\
            0.2421875  1.274033051617343  \\
            0.25  1.2840254166854053  \\
            0.2578125  1.2940961528440194  \\
            0.265625  1.3042458747652652  \\
            0.2734375  1.3144752019421582  \\
            0.28125  1.324784758726456  \\
            0.2890625  1.3351751743667686  \\
            0.296875  1.345647083046963  \\
            0.3046875  1.356201123924874  \\
            0.3125  1.3668379411713099  \\
            0.3203125  1.3775581840093774  \\
            0.328125  1.3883625067541012  \\
            0.3359375  1.3992515688523608  \\
            0.34375  1.4102260349231455  \\
            0.3515625  1.4212865747981112  \\
            0.359375  1.4324338635624725  \\
            0.3671875  1.4436685815962003  \\
            0.375  1.454991414615554  \\
            0.3828125  1.4664030537149313  \\
            0.390625  1.47790419540905  \\
            0.3984375  1.4894955416754607  \\
            0.40625  1.501177799997392  \\
            0.4140625  1.5129516834069334  \\
            0.421875  1.5248179105285529  \\
            0.4296875  1.5367772056229612  \\
            0.4375  1.5488302986313158  \\
            0.4453125  1.5609779252197733  \\
            0.453125  1.573220826824391  \\
            0.4609375  1.5855597506963832  \\
            0.46875  1.597995449947726  \\
            0.4765625  1.6105286835971282  \\
            0.484375  1.6231602166163528  \\
            0.4921875  1.6358908199769142  \\
            0.5  1.6487212706971293  \\
            0.5078125  1.6616523518895454  \\
            0.515625  1.6746848528087384  \\
            0.5234375  1.6878195688994833  \\
            0.53125  1.7010573018453068  \\
            0.5390625  1.7143988596174173  \\
            0.546875  1.7278450565240204  \\
            0.5546875  1.7413967132600183  \\
            0.5625  1.7550546569571062  \\
            0.5703125  1.7688197212342498  \\
            0.578125  1.7826927462485729  \\
            0.5859375  1.7966745787466305  \\
            0.59375  1.8107660721160936  \\
            0.6015625  1.7487205352321953  \\
            0.609375  1.2903856016642876  \\
            0.6171875  1.2371915789325763  \\
            0.625  1.1665488147188374  \\
            0.6328125  1.2103139753475105  \\
            0.640625  1.2190080596699557  \\
            0.6484375  1.1749439611026826  \\
            0.65625  1.0809365471202175  \\
            0.6640625  0.9853810354560785  \\
            0.671875  0.9079422013442966  \\
            0.6796875  0.8456881403318899  \\
            0.6875  0.7972386490004871  \\
            0.6953125  0.7586163922227658  \\
            0.703125  0.7281942783557304  \\
            0.7109375  0.7059723073993807  \\
            0.71875  0.6917518582041485  \\
            0.7265625  0.681163265479532  \\
            0.734375  0.6724189388794162  \\
            0.7421875  0.6655188784038013  \\
            0.75  0.6604453341488162  \\
            0.7578125  0.6568078082292931  \\
            0.765625  0.6544465515103909  \\
            0.7734375  0.6533615639921095  \\
            0.78125  0.6535528456744492  \\
            0.7890625  0.6550203965574097  \\
            0.796875  0.657764216640991  \\
            0.8046875  0.6617843059251934  \\
            0.8125  0.6670195415502216  \\
            0.8203125  0.6721252206005841  \\
            0.828125  0.6765512373381258  \\
            0.8359375  0.6802975917628462  \\
            0.84375  0.6833642838747456  \\
            0.8515625  0.685751313673824  \\
            0.859375  0.6874586811600812  \\
            0.8671875  0.6884863863335173  \\
            0.875  0.6888120477911027  \\
            0.8828125  0.6879432746661829  \\
            0.890625  0.68567863433149  \\
            0.8984375  0.682018126787024  \\
            0.90625  0.6771476811266997  \\
            0.9140625  0.6751577374166309  \\
            0.921875  0.6777216575020465  \\
            0.9296875  0.6848394413829468  \\
            0.9375  0.6961481545814249  \\
            0.9453125  0.7036632385835365  \\
            0.953125  0.7048583257134483  \\
            0.9609375  0.7160143537283179  \\
            0.96875  0.7423524800423974  \\
            0.9765625  0.7551405084534637  \\
            0.984375  0.8005789136310848  \\
            0.9921875  0.8631232591497748  \\
        }
        ;
    \nextgroupplot[width={.5\textwidth}, height={.25\textwidth}, ymode={log}, cycle list name={color list}]
    \addplot+
        table[row sep={\\}]
        {
            \\
            1.0  1.3449688624240168  \\
            2.0  0.3004958140915016  \\
            3.0  0.059685797270837115  \\
            4.0  0.7483123019395072  \\
            5.0  0.0005124472046710536  \\
            6.0  0.024690996949841768  \\
            7.0  0.022463068701592863  \\
            8.0  0.0039393876778139394  \\
            9.0  6.83362504416672e-5  \\
            10.0  7.349827383105358e-5  \\
            11.0  8.328043011545604e-5  \\
            12.0  0.0011253377621449393  \\
            13.0  0.06844155725772046  \\
            14.0  0.002364220924305391  \\
            15.0  0.00021532141547024963  \\
            16.0  3.6197502528274013e-6  \\
            17.0  8.291181800497885e-6  \\
            18.0  8.365238992933603e-6  \\
            19.0  8.902261659024177e-6  \\
            20.0  9.476408207673842e-6  \\
            21.0  1.0087584027634183e-5  \\
            22.0  1.0738177301308771e-5  \\
            23.0  1.14307302366487e-5  \\
            24.0  1.2552878225716916e-5  \\
            25.0  0.00020522478412194096  \\
            26.0  0.027843580864479474  \\
            27.0  0.004456296082194264  \\
            28.0  1.6416456936757783e-5  \\
            29.0  0.0  \\
            30.0  0.0  \\
            31.0  2.360515015756285e-9  \\
            32.0  4.429678584963597e-7  \\
            33.0  1.0186806532896752e-6  \\
            34.0  9.977318268270862e-7  \\
            35.0  1.0291557662033685e-6  \\
            36.0  1.0618246775884882e-6  \\
            37.0  1.0955306110694077e-6  \\
            38.0  1.1303064845792632e-6  \\
            39.0  1.1661862641884624e-6  \\
            40.0  1.203204989985024e-6  \\
            41.0  1.2413988158127554e-6  \\
            42.0  1.280805045078548e-6  \\
            43.0  1.3214621619392353e-6  \\
            44.0  1.3634098738196654e-6  \\
            45.0  1.4066891491359979e-6  \\
            46.0  1.4513422560008538e-6  \\
            47.0  1.4974128044326074e-6  \\
            48.0  1.5449457890850948e-6  \\
            49.0  1.5939876320814061e-6  \\
            50.0  9.85640974352576e-6  \\
            51.0  0.0021776921387844247  \\
            52.0  0.013070000742817566  \\
            53.0  0.00035021774908674847  \\
            54.0  8.432307310377414e-17  \\
            55.0  0.0  \\
            56.0  0.0  \\
            57.0  0.0  \\
            58.0  0.0  \\
            59.0  0.0  \\
            60.0  0.0  \\
            61.0  0.0  \\
            62.0  0.0  \\
            63.0  2.405833450160362e-10  \\
            64.0  5.235825606288711e-8  \\
            65.0  1.250320885552402e-7  \\
            66.0  1.2179878038295597e-7  \\
            67.0  1.237161788643533e-7  \\
            68.0  1.256644252990843e-7  \\
            69.0  1.2764335189340004e-7  \\
            70.0  1.2965344189843492e-7  \\
            71.0  1.3169518615560594e-7  \\
            72.0  1.3376908437195967e-7  \\
            73.0  1.3587564074400704e-7  \\
            74.0  1.3801537033369948e-7  \\
            75.0  1.4018879689048358e-7  \\
            76.0  1.4239644855612577e-7  \\
            77.0  1.446388662954684e-7  \\
            78.0  1.4691659667998014e-7  \\
            79.0  1.4923019691498451e-7  \\
            80.0  1.5158023041958457e-7  \\
            81.0  1.5396727218522366e-7  \\
            82.0  1.563919034448802e-7  \\
            83.0  1.588547169015242e-7  \\
            84.0  1.6135631515645787e-7  \\
            85.0  1.6389730693455729e-7  \\
            86.0  1.6647831450888884e-7  \\
            87.0  1.6909996692487939e-7  \\
            88.0  1.7176290412895812e-7  \\
            89.0  1.7446777612895037e-7  \\
            90.0  1.7721524356306695e-7  \\
            91.0  1.8000597776929306e-7  \\
            92.0  1.8284065945312067e-7  \\
            93.0  1.857199808733001e-7  \\
            94.0  1.8864464523121738e-7  \\
            95.0  1.9161536623374387e-7  \\
            96.0  1.9463286832915871e-7  \\
            97.0  1.9769789041945707e-7  \\
            98.0  2.0081117904635626e-7  \\
            99.0  2.0397349467793369e-7  \\
            100.0  2.0718560993920887e-7  \\
            101.0  2.1044830882598692e-7  \\
            102.0  0.000175399330649641  \\
            103.0  0.004694430710915491  \\
            104.0  0.0014598891229744488  \\
            105.0  2.815843308992858e-18  \\
            106.0  0.0  \\
            107.0  0.0  \\
            108.0  0.0  \\
            109.0  0.0  \\
            110.0  0.0  \\
            111.0  0.0  \\
            112.0  0.0  \\
            113.0  0.0  \\
            114.0  0.0  \\
            115.0  0.0  \\
            116.0  0.0  \\
            117.0  0.0  \\
            118.0  0.0  \\
            119.0  0.0  \\
            120.0  0.0  \\
            121.0  0.0  \\
            122.0  0.0  \\
            123.0  0.0  \\
            124.0  0.0  \\
            125.0  0.0  \\
            126.0  0.0  \\
            127.0  5.08362962667908e-15  \\
            128.0  5.1324141881759624e-9  \\
            129.0  1.4944025134703553e-8  \\
            130.0  1.5047393191597918e-8  \\
            131.0  1.5165410927053813e-8  \\
            132.0  1.528435517039739e-8  \\
            133.0  1.540423189005935e-8  \\
            134.0  1.552504916377956e-8  \\
            135.0  1.5646813482893263e-8  \\
            136.0  1.576953309579375e-8  \\
            137.0  1.5893214661000243e-8  \\
            138.0  1.60178667812233e-8  \\
            139.0  1.61434964038136e-8  \\
            140.0  1.6270111341532e-8  \\
            141.0  1.6397719265759392e-8  \\
            142.0  1.652632794502118e-8  \\
            143.0  1.665594528185016e-8  \\
            144.0  1.6786580040936694e-8  \\
            145.0  1.6918239004182212e-8  \\
            146.0  1.7050930473105907e-8  \\
            147.0  1.718466279693187e-8  \\
            148.0  1.7319443912020005e-8  \\
            149.0  1.7455281107678355e-8  \\
            150.0  1.7592184885922846e-8  \\
            151.0  1.7730161949724987e-8  \\
            152.0  1.7869221269339874e-8  \\
            153.0  1.8009371438587607e-8  \\
            154.0  1.8150620811896445e-8  \\
            155.0  1.829297799349483e-8  \\
            156.0  1.8436451403730514e-8  \\
            157.0  1.8581050639093766e-8  \\
            158.0  1.8726783727884833e-8  \\
            159.0  1.8873659961282652e-8  \\
            160.0  1.9021687169828994e-8  \\
            161.0  1.917087547043117e-8  \\
            162.0  1.9321235056324326e-8  \\
            163.0  1.947277322722485e-8  \\
            164.0  1.9625500446984745e-8  \\
            165.0  1.9779425319832455e-8  \\
            166.0  1.9934556672041026e-8  \\
            167.0  2.009090534910163e-8  \\
            168.0  2.0248481107099092e-8  \\
            169.0  2.0407291766166846e-8  \\
            170.0  2.0567347429334415e-8  \\
            171.0  2.072865892127629e-8  \\
            172.0  2.089123602583287e-8  \\
            173.0  2.105508796479416e-8  \\
            174.0  2.1220224952211986e-8  \\
            175.0  2.1386657514388396e-8  \\
            176.0  2.1554395643330615e-8  \\
            177.0  2.1723448748178775e-8  \\
            178.0  2.1893828319741182e-8  \\
            179.0  2.206554387818027e-8  \\
            180.0  2.2238606248170534e-8  \\
            181.0  2.2413024866607678e-8  \\
            182.0  2.258881222350073e-8  \\
            183.0  2.2765978616168248e-8  \\
            184.0  2.294453503581817e-8  \\
            185.0  2.3124491502213296e-8  \\
            186.0  2.3305859089828296e-8  \\
            187.0  2.3488649206204748e-8  \\
            188.0  2.3672872787039445e-8  \\
            189.0  2.3858541045584936e-8  \\
            190.0  2.4045666083272188e-8  \\
            191.0  2.4234258891309146e-8  \\
            192.0  2.4424330599681632e-8  \\
            193.0  2.4615893254309462e-8  \\
            194.0  2.480895791578952e-8  \\
            195.0  2.5003535783496567e-8  \\
            196.0  2.5199641061346423e-8  \\
            197.0  2.5397284907080486e-8  \\
            198.0  2.5596478013534263e-8  \\
            199.0  2.5797233765834093e-8  \\
            200.0  2.5999563779688373e-8  \\
            201.0  2.6203480524289446e-8  \\
            202.0  2.6408997690074987e-8  \\
            203.0  2.6616125938655477e-8  \\
            204.0  2.6824878797451054e-8  \\
            205.0  2.7038776329748298e-8  \\
            206.0  0.0037373193994677547  \\
            207.0  9.2943092638286e-20  \\
            208.0  0.0  \\
            209.0  0.0  \\
            210.0  0.0  \\
            211.0  0.0  \\
            212.0  0.0  \\
            213.0  0.0  \\
            214.0  0.0  \\
            215.0  0.0  \\
            216.0  0.0  \\
            217.0  0.0  \\
            218.0  0.0  \\
            219.0  0.0  \\
            220.0  0.0  \\
            221.0  0.0  \\
            222.0  0.0  \\
            223.0  0.0  \\
            224.0  0.0  \\
            225.0  0.0  \\
            226.0  0.0  \\
            227.0  0.0  \\
            228.0  0.0  \\
            229.0  0.0  \\
            230.0  0.0  \\
            231.0  0.0  \\
            232.0  0.0  \\
            233.0  0.0  \\
            234.0  0.0  \\
            235.0  0.0  \\
            236.0  0.0  \\
            237.0  0.0  \\
            238.0  0.0  \\
            239.0  0.0  \\
            240.0  0.0  \\
            241.0  0.0  \\
            242.0  0.0  \\
            243.0  0.0  \\
            244.0  0.0  \\
            245.0  0.0  \\
            246.0  0.0  \\
            247.0  0.0  \\
            248.0  0.0  \\
            249.0  0.0  \\
            250.0  0.0  \\
            251.0  0.0  \\
            252.0  0.0  \\
            253.0  0.0  \\
            254.0  0.0  \\
            255.0  9.633042886946583e-19  \\
            256.0  1.0721366649557947e-13  \\
        }
        ;
    \addplot+
        table[row sep={\\}]
        {
            \\
            1.0  1.2950333303580197  \\
            2.0  0.05784323992601773  \\
            3.0  0.11430733339185853  \\
            4.0  0.4325954921873454  \\
            5.0  0.06236092009031331  \\
            6.0  0.04193434573877511  \\
            7.0  0.35164136178984406  \\
            8.0  0.05994545106247551  \\
            9.0  0.03369208306244365  \\
            10.0  0.000630058643889244  \\
            11.0  8.379669400946557e-5  \\
            12.0  0.0003449751004710297  \\
            13.0  0.09989508259916541  \\
            14.0  0.06148898076728966  \\
            15.0  0.036011601382095006  \\
            16.0  0.01877160136852461  \\
            17.0  0.017941050199963844  \\
            18.0  0.0003444500174569718  \\
            19.0  8.902261659818792e-6  \\
            20.0  9.476408208269867e-6  \\
            21.0  1.0087584027598183e-5  \\
            22.0  1.073817730149895e-5  \\
            23.0  1.1430730236118643e-5  \\
            24.0  2.356650794616185e-5  \\
            25.0  0.008755089262481407  \\
            26.0  0.06292275815356091  \\
            27.0  0.026382276322527765  \\
            28.0  0.01460907113684242  \\
            29.0  0.0  \\
            30.0  0.0  \\
            31.0  0.0009516028353332508  \\
            32.0  0.012008153895611634  \\
            33.0  0.008776754094938347  \\
            34.0  0.00014912794473451727  \\
            35.0  1.0291557663089801e-6  \\
            36.0  1.0618246768535154e-6  \\
            37.0  1.095530610567978e-6  \\
            38.0  1.1303064852666145e-6  \\
            39.0  1.1661862644472627e-6  \\
            40.0  1.203204990590296e-6  \\
            41.0  1.2413988167134252e-6  \\
            42.0  1.280805045358549e-6  \\
            43.0  1.3214621613833098e-6  \\
            44.0  1.3634098735252609e-6  \\
            45.0  1.4066891485738178e-6  \\
            46.0  1.451342255848604e-6  \\
            47.0  1.4974128031274567e-6  \\
            48.0  1.5449457884459426e-6  \\
            49.0  1.5939876328759214e-6  \\
            50.0  0.00024481384377668977  \\
            51.0  0.0020212854439974406  \\
            52.0  0.1574840155121907  \\
            53.0  0.043716409812914296  \\
            54.0  0.0017663849681346475  \\
            55.0  0.0  \\
            56.0  0.0  \\
            57.0  0.0  \\
            58.0  0.0  \\
            59.0  0.0  \\
            60.0  0.0  \\
            61.0  0.0  \\
            62.0  0.0  \\
            63.0  0.0014570081891225291  \\
            64.0  0.008335323568021453  \\
            65.0  0.0032025560713410013  \\
            66.0  1.217983708996942e-7  \\
            67.0  1.2371617841912983e-7  \\
            68.0  1.2566442454460964e-7  \\
            69.0  1.2764335081381337e-7  \\
            70.0  1.296534416974496e-7  \\
            71.0  1.316951866026569e-7  \\
            72.0  1.3376908407236028e-7  \\
            73.0  1.3587564045505226e-7  \\
            74.0  1.3801537057918374e-7  \\
            75.0  1.4018879659421057e-7  \\
            76.0  1.4239644893236165e-7  \\
            77.0  1.4463886672383263e-7  \\
            78.0  1.469165971600328e-7  \\
            79.0  1.4923019634887658e-7  \\
            80.0  1.5158023007934388e-7  \\
            81.0  1.5396727189042202e-7  \\
            82.0  1.563919037526048e-7  \\
            83.0  1.58854717571816e-7  \\
            84.0  1.6135631548505618e-7  \\
            85.0  1.6389730785132283e-7  \\
            86.0  1.6647831496145783e-7  \\
            87.0  1.6909996696760066e-7  \\
            88.0  1.717629033661581e-7  \\
            89.0  1.744677757168563e-7  \\
            90.0  1.7721524387334504e-7  \\
            91.0  1.800059778241817e-7  \\
            92.0  1.828406591433911e-7  \\
            93.0  1.8571998088040573e-7  \\
            94.0  1.886446446467721e-7  \\
            95.0  1.9161536537255883e-7  \\
            96.0  1.9463286833538308e-7  \\
            97.0  1.976978903147758e-7  \\
            98.0  2.0081117892595855e-7  \\
            99.0  2.0397349414574733e-7  \\
            100.0  2.0718561006247922e-7  \\
            101.0  2.1044831285440692e-7  \\
            102.0  0.005187824590040085  \\
            103.0  0.03442144255242336  \\
            104.0  0.01569372998643822  \\
            105.0  0.00706553987253859  \\
            106.0  0.0  \\
            107.0  0.0  \\
            108.0  0.0  \\
            109.0  0.0  \\
            110.0  0.0  \\
            111.0  0.0  \\
            112.0  0.0  \\
            113.0  0.0  \\
            114.0  0.0  \\
            115.0  0.0  \\
            116.0  0.0  \\
            117.0  0.0  \\
            118.0  0.0  \\
            119.0  0.0  \\
            120.0  0.0  \\
            121.0  0.0  \\
            122.0  0.0  \\
            123.0  0.0  \\
            124.0  0.0  \\
            125.0  0.0  \\
            126.0  0.0  \\
            127.0  0.0015285056277785916  \\
            128.0  0.006909514141664819  \\
            129.0  1.4937497450295876e-8  \\
            130.0  1.504792896443238e-8  \\
            131.0  1.516551523425356e-8  \\
            132.0  1.5284374210552945e-8  \\
            133.0  1.540423552211834e-8  \\
            134.0  1.552504944762702e-8  \\
            135.0  1.5646813823399154e-8  \\
            136.0  1.5769533869088267e-8  \\
            137.0  1.5893214922549257e-8  \\
            138.0  1.6017866520880964e-8  \\
            139.0  1.6143495301205485e-8  \\
            140.0  1.6270110814442466e-8  \\
            141.0  1.639771941924162e-8  \\
            142.0  1.6526329055126325e-8  \\
            143.0  1.665594640470885e-8  \\
            144.0  1.6786580403220613e-8  \\
            145.0  1.691823897231087e-8  \\
            146.0  1.7050930478838104e-8  \\
            147.0  1.7184662186036807e-8  \\
            148.0  1.73194436446955e-8  \\
            149.0  1.7455281510086262e-8  \\
            150.0  1.759218522171377e-8  \\
            151.0  1.77301627441285e-8  \\
            152.0  1.786922201116195e-8  \\
            153.0  1.8009372135842653e-8  \\
            154.0  1.8150621266390565e-8  \\
            155.0  1.8292977856986286e-8  \\
            156.0  1.8436451115113483e-8  \\
            157.0  1.858105060892589e-8  \\
            158.0  1.872678371843108e-8  \\
            159.0  1.8873659376585347e-8  \\
            160.0  1.902168664995393e-8  \\
            161.0  1.917087561389214e-8  \\
            162.0  1.9321234397647614e-8  \\
            163.0  1.9472773023522646e-8  \\
            164.0  1.9625500124308877e-8  \\
            165.0  1.9779425046864068e-8  \\
            166.0  1.993455743827486e-8  \\
            167.0  2.0090906106850405e-8  \\
            168.0  2.0248481081192388e-8  \\
            169.0  2.0407291753761867e-8  \\
            170.0  2.056734785057844e-8  \\
            171.0  2.0728659673392034e-8  \\
            172.0  2.0891237106413266e-8  \\
            173.0  2.1055089281566203e-8  \\
            174.0  2.1220225578952207e-8  \\
            175.0  2.138665838477705e-8  \\
            176.0  2.155439618152258e-8  \\
            177.0  2.172344897943474e-8  \\
            178.0  2.1893827802345843e-8  \\
            179.0  2.2065543458831104e-8  \\
            180.0  2.223860569554928e-8  \\
            181.0  2.2413025227757142e-8  \\
            182.0  2.2588812834310725e-8  \\
            183.0  2.276597922773288e-8  \\
            184.0  2.2944535130622227e-8  \\
            185.0  2.3124491302147422e-8  \\
            186.0  2.3305858716133106e-8  \\
            187.0  2.3488648776601555e-8  \\
            188.0  2.3672871889386704e-8  \\
            189.0  2.385854020830826e-8  \\
            190.0  2.4045665397815295e-8  \\
            191.0  2.4234259080098288e-8  \\
            192.0  2.4424330910349544e-8  \\
            193.0  2.461589286820278e-8  \\
            194.0  2.4808957459546486e-8  \\
            195.0  2.500353669926299e-8  \\
            196.0  2.5199641157485264e-8  \\
            197.0  2.5397284521475358e-8  \\
            198.0  2.559647797600633e-8  \\
            199.0  2.579723376375267e-8  \\
            200.0  2.599956446673349e-8  \\
            201.0  2.62034840218263e-8  \\
            202.0  2.6409014868931207e-8  \\
            203.0  2.661621111892843e-8  \\
            204.0  2.6825305259763903e-8  \\
            205.0  2.7037738503199734e-8  \\
            206.0  0.11066905824859176  \\
            207.0  0.02826215949015436  \\
            208.0  0.0  \\
            209.0  0.0  \\
            210.0  0.0  \\
            211.0  0.0  \\
            212.0  0.0  \\
            213.0  0.0  \\
            214.0  0.0  \\
            215.0  0.0  \\
            216.0  0.0  \\
            217.0  0.0  \\
            218.0  0.0  \\
            219.0  0.0  \\
            220.0  0.0  \\
            221.0  0.0  \\
            222.0  0.0  \\
            223.0  0.0  \\
            224.0  0.0  \\
            225.0  0.0  \\
            226.0  0.0  \\
            227.0  0.0  \\
            228.0  0.0  \\
            229.0  0.0  \\
            230.0  0.0  \\
            231.0  0.0  \\
            232.0  0.0  \\
            233.0  0.0  \\
            234.0  0.0  \\
            235.0  0.0  \\
            236.0  0.0  \\
            237.0  0.0  \\
            238.0  0.0  \\
            239.0  0.0  \\
            240.0  0.0  \\
            241.0  0.0  \\
            242.0  0.0  \\
            243.0  0.0  \\
            244.0  0.0  \\
            245.0  0.0  \\
            246.0  0.0  \\
            247.0  0.0  \\
            248.0  0.0  \\
            249.0  0.0  \\
            250.0  0.0  \\
            251.0  0.0  \\
            252.0  0.0  \\
            253.0  0.0  \\
            254.0  0.0  \\
            255.0  0.0016931647601229043  \\
            256.0  0.00799568382324315  \\
        }
        ;
    \addplot+
        table[row sep={\\}]
        {
            \\
            1.0  1.1393800347936223  \\
            2.0  0.8877760445192617  \\
            3.0  0.07279437812768463  \\
            4.0  0.7013651893037076  \\
            5.0  0.17289910979258766  \\
            6.0  0.07674358608090745  \\
            7.0  0.7452450817661707  \\
            8.0  0.18589172469584564  \\
            9.0  0.1406400497195439  \\
            10.0  0.0038862138161061942  \\
            11.0  8.480088196331044e-5  \\
            12.0  0.0006160311875218344  \\
            13.0  0.1731896285708262  \\
            14.0  0.0016546971269310595  \\
            15.0  0.022979341700889142  \\
            16.0  0.09980587691195765  \\
            17.0  0.10178867122450974  \\
            18.0  0.002461459055891466  \\
            19.0  8.902261659143873e-6  \\
            20.0  9.476408207281795e-6  \\
            21.0  1.0087584027179686e-5  \\
            22.0  1.0738177300823049e-5  \\
            23.0  1.1430730236183795e-5  \\
            24.0  4.498918432724146e-5  \\
            25.0  0.02604601002767968  \\
            26.0  0.20016010085407235  \\
            27.0  0.05652436193441094  \\
            28.0  0.026268895179663528  \\
            29.0  0.0  \\
            30.0  0.0  \\
            31.0  0.01763205822831764  \\
            32.0  0.04890066377360025  \\
            33.0  0.06214098818733421  \\
            34.0  0.0011778739013687672  \\
            35.0  1.0291557663022477e-6  \\
            36.0  1.0618246769900086e-6  \\
            37.0  1.0955306103954676e-6  \\
            38.0  1.1303064851031497e-6  \\
            39.0  1.1661862648337795e-6  \\
            40.0  1.203204990724016e-6  \\
            41.0  1.2413988163817447e-6  \\
            42.0  1.2808050451618147e-6  \\
            43.0  1.3214621614188182e-6  \\
            44.0  1.3634098736947653e-6  \\
            45.0  1.406689148733542e-6  \\
            46.0  1.4513422561396316e-6  \\
            47.0  1.497412804654652e-6  \\
            48.0  1.5449457895222451e-6  \\
            49.0  1.5939876329140734e-6  \\
            50.0  0.000701830939915174  \\
            51.0  0.028548939645518814  \\
            52.0  0.28244235201682083  \\
            53.0  0.06177431187577999  \\
            54.0  0.005248724055688161  \\
            55.0  0.0  \\
            56.0  0.0  \\
            57.0  0.0  \\
            58.0  0.0  \\
            59.0  0.0  \\
            60.0  0.0  \\
            61.0  0.0  \\
            62.0  0.0  \\
            63.0  0.01617788333429564  \\
            64.0  0.003255890874605788  \\
            65.0  0.025106816648982407  \\
            66.0  1.217987805491455e-7  \\
            67.0  1.2371617905126975e-7  \\
            68.0  1.2566442532391253e-7  \\
            69.0  1.2764335204085153e-7  \\
            70.0  1.2965344214910246e-7  \\
            71.0  1.3169518631173105e-7  \\
            72.0  1.3376908417073174e-7  \\
            73.0  1.3587564087758075e-7  \\
            74.0  1.3801537041002732e-7  \\
            75.0  1.40188796404761e-7  \\
            76.0  1.423964481606088e-7  \\
            77.0  1.446388659398501e-7  \\
            78.0  1.4691659708243598e-7  \\
            79.0  1.4923019706070129e-7  \\
            80.0  1.5158023020100941e-7  \\
            81.0  1.539672718486873e-7  \\
            82.0  1.563919040659112e-7  \\
            83.0  1.5885471760235248e-7  \\
            84.0  1.6135631483726876e-7  \\
            85.0  1.638973075868133e-7  \\
            86.0  1.6647831532767832e-7  \\
            87.0  1.6909996730651855e-7  \\
            88.0  1.71762904135897e-7  \\
            89.0  1.744677761983393e-7  \\
            90.0  1.7721524354918916e-7  \\
            91.0  1.800059777345986e-7  \\
            92.0  1.8284065914087044e-7  \\
            93.0  1.8571998060268324e-7  \\
            94.0  1.8864464512713397e-7  \\
            95.0  1.9161536554679337e-7  \\
            96.0  1.9463286829446425e-7  \\
            97.0  1.9769789033619034e-7  \\
            98.0  2.0081117893533396e-7  \\
            99.0  2.0397349526080077e-7  \\
            100.0  2.0718561041799255e-7  \\
            101.0  2.1044830918366675e-7  \\
            102.0  0.014937522640975108  \\
            103.0  0.16968978614101946  \\
            104.0  0.05813427439513724  \\
            105.0  0.020994896222752645  \\
            106.0  0.0  \\
            107.0  0.0  \\
            108.0  0.0  \\
            109.0  0.0  \\
            110.0  0.0  \\
            111.0  0.0  \\
            112.0  0.0  \\
            113.0  0.0  \\
            114.0  0.0  \\
            115.0  0.0  \\
            116.0  0.0  \\
            117.0  0.0  \\
            118.0  0.0  \\
            119.0  0.0  \\
            120.0  0.0  \\
            121.0  0.0  \\
            122.0  0.0  \\
            123.0  0.0  \\
            124.0  0.0  \\
            125.0  0.0  \\
            126.0  0.0  \\
            127.0  0.014808425545828285  \\
            128.0  0.03667070820855846  \\
            129.0  1.494402528112726e-8  \\
            130.0  1.5047393100639517e-8  \\
            131.0  1.5165411176853993e-8  \\
            132.0  1.5284355378564207e-8  \\
            133.0  1.540423182067041e-8  \\
            134.0  1.5525048843723077e-8  \\
            135.0  1.5646813275593807e-8  \\
            136.0  1.576953290172156e-8  \\
            137.0  1.589321461936688e-8  \\
            138.0  1.6017866782958023e-8  \\
            139.0  1.6143496095466503e-8  \\
            140.0  1.6270110994587306e-8  \\
            141.0  1.6397719238003816e-8  \\
            142.0  1.652632833359924e-8  \\
            143.0  1.6655945678234474e-8  \\
            144.0  1.6786580070426993e-8  \\
            145.0  1.6918238851526546e-8  \\
            146.0  1.705093020942794e-8  \\
            147.0  1.7184662567948372e-8  \\
            148.0  1.7319443657015654e-8  \\
            149.0  1.745528116492423e-8  \\
            150.0  1.759218470898105e-8  \\
            151.0  1.773016225503632e-8  \\
            152.0  1.78692212415843e-8  \\
            153.0  1.8009371521854334e-8  \\
            154.0  1.8150620756385294e-8  \\
            155.0  1.829297749389447e-8  \\
            156.0  1.8436451365566597e-8  \\
            157.0  1.8581050614807637e-8  \\
            158.0  1.872678326297894e-8  \\
            159.0  1.8873659579643487e-8  \\
            160.0  1.902168730860687e-8  \\
            161.0  1.9170875997787107e-8  \\
            162.0  1.932123527836893e-8  \\
            163.0  1.9472773393758303e-8  \\
            164.0  1.962550066902935e-8  \\
            165.0  1.9779425458610334e-8  \\
            166.0  1.993455703286351e-8  \\
            167.0  2.0090906022174337e-8  \\
            168.0  2.024848128057144e-8  \\
            169.0  2.0407291745350165e-8  \\
            170.0  2.0567348046895972e-8  \\
            171.0  2.0728659344548817e-8  \\
            172.0  2.0891236157671855e-8  \\
            173.0  2.1055088561539037e-8  \\
            174.0  2.1220225326912256e-8  \\
            175.0  2.138665780582194e-8  \\
            176.0  2.1554395157608042e-8  \\
            177.0  2.1723448137556112e-8  \\
            178.0  2.1893827972796487e-8  \\
            179.0  2.206554325367982e-8  \\
            180.0  2.2238606040003717e-8  \\
            181.0  2.2413025130285646e-8  \\
            182.0  2.2588812625956578e-8  \\
            183.0  2.276597890760179e-8  \\
            184.0  2.2944534799895777e-8  \\
            185.0  2.3124491377313205e-8  \\
            186.0  2.3305858798394752e-8  \\
            187.0  2.348864849843757e-8  \\
            188.0  2.3672872800917233e-8  \\
            189.0  2.3858540726395816e-8  \\
            190.0  2.404566576408307e-8  \\
            191.0  2.4234258322319846e-8  \\
            192.0  2.442433011395906e-8  \\
            193.0  2.461589282409804e-8  \\
            194.0  2.4808957443944735e-8  \\
            195.0  2.5003536102685686e-8  \\
            196.0  2.5199641262574346e-8  \\
            197.0  2.5397284802997078e-8  \\
            198.0  2.5596477770672976e-8  \\
            199.0  2.5797233738078518e-8  \\
            200.0  2.599956372417722e-8  \\
            201.0  2.620348107940096e-8  \\
            202.0  2.640899776640282e-8  \\
            203.0  2.6616125747835895e-8  \\
            204.0  2.68248787765879e-8  \\
            205.0  2.703877577392433e-8  \\
            206.0  0.31866261668639373  \\
            207.0  0.08397958489101055  \\
            208.0  0.0  \\
            209.0  0.0  \\
            210.0  0.0  \\
            211.0  0.0  \\
            212.0  0.0  \\
            213.0  0.0  \\
            214.0  0.0  \\
            215.0  0.0  \\
            216.0  0.0  \\
            217.0  0.0  \\
            218.0  0.0  \\
            219.0  0.0  \\
            220.0  0.0  \\
            221.0  0.0  \\
            222.0  0.0  \\
            223.0  0.0  \\
            224.0  0.0  \\
            225.0  0.0  \\
            226.0  0.0  \\
            227.0  0.0  \\
            228.0  0.0  \\
            229.0  0.0  \\
            230.0  0.0  \\
            231.0  0.0  \\
            232.0  0.0  \\
            233.0  0.0  \\
            234.0  0.0  \\
            235.0  0.0  \\
            236.0  0.0  \\
            237.0  0.0  \\
            238.0  0.0  \\
            239.0  0.0  \\
            240.0  0.0  \\
            241.0  0.0  \\
            242.0  0.0  \\
            243.0  0.0  \\
            244.0  0.0  \\
            245.0  0.0  \\
            246.0  0.0  \\
            247.0  0.0  \\
            248.0  0.0  \\
            249.0  0.0  \\
            250.0  0.0  \\
            251.0  0.0  \\
            252.0  0.0  \\
            253.0  0.0  \\
            254.0  0.0  \\
            255.0  0.014610429738475426  \\
            256.0  0.06251987897073587  \\
        }
        ;
\end{groupplot}
\end{tikzpicture}

%% file: weightedtabular.tikz
\begin{tabular}{l|l|l|l|l|}
&\multicolumn{2}{l}{wavelet extension}&\multicolumn{2}{|l|}{spline extension}\\
&$\|x\|$  & $\|Ax-b\|$  & $\|x\|$  & $\|\hat Ax-b\|$\\\hline
Reduced AZ&1.58&1.14e-06&47.00&1.14e-06\\
Weighted reduced AZ&3.37&1.14e-06&&\\
Pivoted QR&1.56&1.14e-06&47.00&1.14e-06\\\hline
Sparse AZ&1.58&1.14e-06&47.00&1.14e-06\\
Weighted sparse AZ&3.38&1.14e-06&&
\end{tabular}

%% file: verticaldualtabulardb2.tikz
\begin{tabular}{|l|l|l|} 
 \hline
k&db2&\\\hline
1&2.6389584337646843&\\\hline
2&3.8637033051562737&0.12940952255126012\\\hline
3&&0.7071067811865474\\\hline
4&-1.0352761804100834&-0.48296291314453427\\\hline
5&0.18946869098150604&\\\hline
\end{tabular} 

%% file: verticaldualtabularcdf31.tikz
\begin{tabular}{|l|l|l|} 
 \hline
k&cdf31&\\\hline
-2&0.3535533905932738&\\\hline
-1&1.4142135623730951&-0.17677669529663703\\\hline
0&2.121320343559643&0.7071067811865476\\\hline
1&1.4142135623730951&-0.17677669529663687\\\hline
2&0.3535533905932738&\\\hline
\end{tabular} 

%% file: verticaldualtabulardb3.tikz
\begin{tabular}{|l|l|l|} 
 \hline
k&db3&\\\hline
1&1.711703195299721&\\\hline
2&3.63830500187592&-0.0012134914677774413\\\hline
3&1.247682791903474&0.02435873888065617\\\hline
4&-1.0913117263517993&0.16899390552758248\\\hline
5&-0.04234322675143443&0.7980018678049798\\\hline
6&0.269457311225099&2.4036661942863304\\\hline
7&-0.0892122768840828&-5.741389143657221\\\hline
8&0.011976537996970554&2.701135319218721\\\hline
9&0.0005966411785128591&\\\hline
\end{tabular} 

%% file: verticaldualtabularcdf42.tikz
\begin{tabular}{|l|l|l|} 
 \hline
k&cdf42&\\\hline
-3&0.05892556509887896&\\\hline
-2&0.4714045207910316&0.058925565098878696\\\hline
-1&1.355287997274216&-0.4714045207910281\\\hline
0&1.8856180831641276&1.1785113019775755\\\hline
1&1.3552879972742171&-0.4714045207910299\\\hline
2&0.47140452079103&0.05892556509887848\\\hline
3&0.05892556509887875&\\\hline
\end{tabular} 

%% file: verticaldualtabulardb4.tikz
\begin{tabular}{|l|l|l|} 
 \hline
k&db4&\\\hline
1&0.9281184683951778&-0.0009679452571134971\\\hline
2&2.848706884229112&0.00734700187024821\\\hline
3&2.481365240821777&0.07373328922533529\\\hline
4&-0.09570535866183345&0.10889244990221354\\\hline
5&-0.6845007058653815&0.42109222379000844\\\hline
6&0.11203530716940972&2.166676570753769\\\hline
7&0.09701428436973923&1.8834043099303606\\\hline
8&-0.0332746298543078&-16.336082959145102\\\hline
9&0.006087195587662892&17.749143869466998\\\hline
10&-0.003388335759322769&-5.51513309994731\\\hline
11&0.0003434396086452516&-0.3034016016000655\\\hline
12&5.325762313281169e-5&0.09884928187755511\\\hline
13&-7.981714296518615e-7&\\\hline
\end{tabular} 